\documentclass[a4paper,oneside,english]{amsart}
\usepackage[T1]{fontenc}
\usepackage[latin9]{inputenc}
\usepackage{babel}
\usepackage{mathrsfs}
\usepackage{amsthm}
\usepackage{amsbsy}
\usepackage{amstext}
\usepackage{amssymb}
\usepackage[unicode=true,
 bookmarks=true,bookmarksnumbered=false,bookmarksopen=false,
 breaklinks=false,pdfborder={0 0 1},backref=page,colorlinks=false]
 {hyperref}
\usepackage{breakurl}

\makeatletter

\numberwithin{equation}{section}
\numberwithin{figure}{section}
\newenvironment{lyxlist}[1]
{\begin{list}{}
{\settowidth{\labelwidth}{#1}
 \setlength{\leftmargin}{\labelwidth}
 \addtolength{\leftmargin}{\labelsep}
 }}
{\end{list}}
  \theoremstyle{remark}
  \newtheorem*{rem*}{\protect\remarkname}
  \theoremstyle{definition}
  \newtheorem*{example*}{\protect\examplename}

\usepackage{calc}
\usepackage{upgreek}

\AtBeginDocument{
  
}

\makeatother

  \providecommand{\examplename}{Example}
  \providecommand{\remarkname}{Remark}

\begin{document}
\title[Classical Tensors]{Interpretations and Representations of\\Classical Tensors}

\author{Dan Jonsson}

\date{January 5, 2014}

\address{Department of Sociology and Work Science, University of Gothenburg,
Box 720, SE 405 30 Göteborg, Sweden.}

\thanks{\emph{E-mail}: dan.jonsson@gu.se. }
\begin{abstract}
Classical tensors, the familiar mathematical objects denoted by symbols
such as $t_{i}$, $t^{ij}$ and $t_{k}^{ij}$, are usually interpreted
either as 'coordinatizable objects' with coordinates changing in a
specific way under a change of coordinate system or as elements of
tensor spaces of the form $V^{\otimes n}\otimes\left(V^{*}\right)^{\otimes m}$.
An alternative interpretation of classical tensors as linear tensor
maps of the form $V^{\otimes m}\rightarrow V^{\otimes n}$ is presented
here. In this interpretation, tensor multiplication is seen as generalized
function composition. Representations of classical tensors by means
of arrays are also considered.
\end{abstract}
\maketitle
\tableofcontents{}

\section{Introduction}

\emph{Classical tensors}, associated with Ricci and Levi-Civita \cite{key-6}
and denoted by symbols such as $t_{i}$, $t^{ij}$ and $t_{k}^{ij}$,
have traditionally been defined in a basis-dependent way, using double-index
notation. For example, $t_{k}^{ij}$ can be seen as a system of scalars
collectively representing a tensor relative to a basis for a vector
space; these scalars change in a specific way under a change of basis,
or equivalently a change of coordinate map. Tensors are thus defined
\emph{only} \emph{indirectly}, as certain \emph{unspecified coordinatizable
objects} which are \emph{represented} by systems of scalars. \newpage{}

This definition and the accompanying notation are problematic in several
ways:
\begin{lyxlist}{00.00.0000}
\item [{$\qquad$(a)}] There is a void in the definition; it says how tensors
can be recognized -- namely, by the behavior of their scalar representations
under a change of basis or coordinate map -- but not what they \emph{are}.
\item [{$\qquad$(b)}] The definition and notation are not coordinate-free,
not manifestly basis-independent; the situation is the same as if
we use tuples of coordinates to refer to vectors, writing, for example,
$\left(\ldots u_{i}\ldots\right)+\left(\ldots v_{i}\ldots\right)=\left(\ldots w_{i}\ldots\right)$
instead of $\mathbf{u}+\mathbf{v}=\mathbf{w}$.
\item [{$\qquad$(c)}] In many cases, indices may make formulas somewhat
difficult to read and understand. For example, $\boldsymbol{g}\!\left(\mathbf{u},\mathbf{v}\right)$
without indices is less cluttered and maybe more transparent than
$g_{ij}v^{i}u^{j}$.
\end{lyxlist}
Addressing primarily the first two problems, several 'intrinsic' definitions
of classical tensors have been proposed. For example, in one common
interpretation, they are regarded as elements of tensor spaces of
the form $V^{\otimes n}\otimes\left(V^{*}\right)^{\otimes m}$. This
works, formally, but is quite non-intuitive.

There have also been attempts since the early days of tensor calculus
to develop an index-free, 'direct' notation for tensors \cite{key-9,key-11},
but no comprehensive, generally accepted index-free notation has yet
emerged. Apparently, many still share Weyl's \cite{key-12} opinion:
\begin{quote}
Various attempts have been made to set up a standard terminology in
this branch of mathematics involving only the vectors themselves and
not their components, analogous to that of vectors in vector analysis.
This is highly expedient in the latter, but very cumbersome for the
much more complicated framework of tensor calculus. In trying to avoid
continual reference to the components we are obliged to adopt an endless
profusion of symbols in addition to an intricate set of rules for
carrying out calculations, so that the balance of advantage is considerably
on the negative side. (p. 54).
\end{quote}
The crux of the matter would seem to be that indices have a conceptual
as well as a computational function in classical tensor notation.
Indices refer to scalars (components) which are subject to arithmetic
operations, but in addition the arrangement of indices shows what
kind of tensor one is dealing with. The form of symbols such as $t_{i}$,
$t^{i}$ and $t_{k}^{ij}$ contains information about invariant properties
of the denoted tensor, and how indices are matched in expressions
such as $t_{i}^{ij}$, $s_{i}t^{i}$ and $s_{i}^{j}t^{ik}$ indicates
what kind of operation on tensors that is being considered. This suggests
that a \emph{coordinate-free} definition of classical tensors should
be combined with a notation which is not necessarily \emph{index-free}.

Recognizing that in some situations indices help more than hinder,
Penrose \cite{key-5} proposed an \emph{abstract index notation} for
tensors. In this notation, indices are retained, but used to distinguish
between different types of tensors, not to identify their components
relative to a particular basis. Abstract index notation thus presupposes
a coordinate-free notion of classical tensors, but it is not based
on any specific 'intrinsic' definition.

In this article, the three objections against the traditional definition
of and notation for tensors are addressed mainly on the basis of an
interpretation of classical tensors as (linear) \emph{tensor maps},
defined in a basis-independent manner. The corresponding notation
has two forms, double-index notation and index-free notation. This
means, for example, that we can write $g_{ij}u^{i}v^{j}$ as $g_{\mathrm{ab}}\circ u^{\mathrm{a}}\circ v^{\mathrm{b}}$
or $\boldsymbol{g}\circ\boldsymbol{u}\otimes\boldsymbol{v}$. In other
words, both a form of abstract index notation and an index-free notation
are available.

It should be emphasized that what is of prime importance here is the
suggested interpretation of classical tensors. The proposed notation
is adapted to that interpretation rather than the other way around,
contrary to Penrose's emphasis on the notation itself.

An important consequence of viewing classical tensors as tensor maps
is that tensor multiplication can then be viewed as generalized function
composition. This means, for example, that we can write $FG=H$, where
$F$, $G$ and $H$ are linear operators on $V$, as $f_{\mathrm{a}}^{\mathrm{b}}\circ g_{\mathrm{c}}^{\mathrm{a}}=h_{\mathrm{c}}^{\mathrm{b}}$,
where the match between the subscript of $f_{\mathrm{a}}^{\mathrm{b}}$
and the superscript of $g_{\mathrm{c}}^{\mathrm{a}}$ means that the
'output' from $g_{\mathrm{c}}^{\mathrm{a}}$ is accepted as 'input'
to $f_{\mathrm{a}}^{\mathrm{b}}$. In general, we can think of subscripts
as 'slots' for arguments of tensors as functions, and superscripts
as 'slots' for values produced by tensors as functions. 

On the other hand, $f_{\mathrm{a}}^{\mathrm{b}}\circ g_{\mathrm{c}}^{\mathrm{d}}=h_{\mathrm{ac}}^{\mathrm{bd}}$
is an example of ordinary tensor multiplication (without contraction).
General products such as $f_{\mathrm{a}}^{\mathrm{b}}\circ g^{\mathrm{ac}}$
and $f_{\mathrm{ab}}\circ g_{\mathrm{d}}^{\mathrm{abc}}$, where some
but not all indices match, can also be formed, and have meanings related
to function composition. Symbols such as $v^{a}$, $g_{\mathrm{a}}$
or $h$ represent functions with scalar arguments and/or scalar values.
In particular,
\[
v^{\mathrm{a}}:K\rightarrow V,\qquad\eta\mapsto\eta\mathbf{v}
\]
is a tensor map representing the vector $\mathbf{v}$, and the composite
function $g_{\mathrm{a}}\circ v^{\mathrm{a}}$ represents the scalar
$g_{\mathrm{a}}\!\left(\mathbf{v}\right)=\boldsymbol{g}\!\left(\mathbf{v}\right)$.
Composition of tensor maps can thus sometimes be interpreted as function
application.

It is important to make a clear distinction between a tensor and its
scalar representation. To emphasize this distinction, special notation
is used; tensors are represented by generalized matrices, \emph{arrays}
such as $\left[t_{i}\right]$, $\left[t^{i}\right]$ and $\left[t_{k}^{ij}\right]$.
For example, a vector-like tensor map $v^{\mathrm{a}}$ is represented
as an array $\left[v^{i}\right]$ relative to some basis; the distinction
between a tensor and its representation is not obscured. Arrays can
also be other objects than scalar representations of tensors; for
example, $\left[\mathbf{e}_{i}\right]$ or $\left[e_{i}^{\mathrm{a}}\right]$
can denote an ordered basis.

While array notation (restricted to scalar arrays) is equivalent to
conventional indicial notation, array notation is conceptually cleaner
and extends notational conventions of elementary linear algebra instead
of introducing a different notational system.

The interpretation of classical tensors as tensor maps, the corresponding
notation and definitions of central notions corresponding to multiplication
and contraction of tensors, permutation of indices, and raising and
lowering of indices are presented in Sections 7--10. Section 6 discusses,
as a background, the conventional basis-independent (coordinate-free)
approach, where tensors are defined as elements of tensors spaces
of the form $V^{\otimes n}\otimes\left(V^{*}\right)^{\otimes m}$
or $\left(V^{*}\right)^{\otimes m}\otimes V^{\otimes n}$. The idea
that multiplication of classical tensors can be regarded as generalized
composition of linear functions is elaborated in Sections 8--9.

The representation of classical tensors as arrays is discussed in
Sections 11--13.

Observations on the equivalence between the present approach to classical
tensors and more conventional approaches are scattered throughout
Sections 6--13.

Finally, Sections 2--5 describe some notation and background, including
elements of algebraic tensor theory. Deviating somewhat from the usual
way of introducing tensors, I present a definition of a unique, strictly
associative tensor product, defining the tensor product in terms of
the tensor algebra rather than vice versa.

Readers familiar with tensors can skip Sections 2--5 -- except the
notational conventions in Section 5 and perhaps Section 4.3 -- without
losing too much context.

\newpage{}

\section{Preliminaries}

\subsection{Notation}

$ $\\
Throughout this article, $U,V,V_{i},W,T$ will denote finite-dimensional
vector spaces over a field $K$ (typically the real numbers) unless
otherwise indicated. Some further notation is shown below:\medskip{}

Scalar: $\text{\ensuremath{\eta}},\xi,\ldots,f_{i},v^{i},\ldots,a_{i}^{j},\ldots$

Vector; tensor as vector: $\mathbf{u},\mathbf{v},\ldots,\mathbf{e}_{i},\ldots;\mathbf{s},\mathbf{t},\ldots$

Linear form: $\boldsymbol{f},\boldsymbol{g},\ldots,\boldsymbol{f}^{i},\ldots;f_{\mathrm{a}},\ldots,f_{\mathrm{a}}^{i},\ldots$

Bilateral tensor: $\boldsymbol{\mathsf{s}},\boldsymbol{\mathsf{t}},\ldots;\mathsf{t},\mathsf{t}_{\mathrm{a}},\mathsf{t}^{\mathrm{a}},\mathsf{t}_{\mathrm{a}}^{\mathrm{b}},\ldots$

Linear tensor map: $\boldsymbol{s},\boldsymbol{t},\ldots,\boldsymbol{e}_{i},\ldots;t,f_{\mathrm{a}},v^{\mathrm{a}},t_{\mathrm{a}}^{\mathrm{bc}},g_{\mathrm{ab}},\ldots,e_{i}^{\mathrm{a}},\ldots$

Separately linear tensor map: $\boldsymbol{s},\boldsymbol{t},\ldots;t,t_{\mathrm{a,b}},t_{\mathrm{a,b}}^{\mathrm{cd}},\ldots$

Array: $\left[f_{i}\right],\left[v^{i}\right],\left[a_{i}^{j}\right],\ldots,\left[\mathbf{e}{}_{i}\right],\ldots,\left[\boldsymbol{e}_{i}\right],\ldots,\left[e_{i}^{\mathrm{a}}\right],\ldots$

Matrix; index-free array: $\boldsymbol{\mathfrak{a}},\boldsymbol{\mathfrak{a}}^{\mathrm{T}},\boldsymbol{\mathfrak{A}},\ldots;\overset{n}{\underset{m}{\boldsymbol{\mathfrak{A}}}},\ldots$

\subsection{Concatenation of tuples and Cartesian products of tuple sets}

$ $\\
\emph{(1)}. A \emph{tuple} or $N$\emph{-tuple} $\left(x_{1},\ldots,x_{N}\right)$
is a list of $N\geq0$ not necessarily distinct elements. Let the
\emph{concatenation} $\boldsymbol{x},\boldsymbol{y}$ of $\boldsymbol{x}=\left(x_{1},\ldots,x_{m}\right)$
and $\boldsymbol{y}=\left(y_{1},\ldots,y_{n}\right)$ be 
\[
\left(x_{1},\ldots,x_{m},y_{1},\ldots,y_{n}\right)
\]
rather than the nested tuple 
\[
\left(\left(x_{1},\ldots,x_{m}\right),\left(y_{1},\ldots,y_{n}\right)\right)
\]
associated with the usual Cartesian product. Note that with concatenation
of tuples defined in this way, the parentheses surrounding a tuple
are no longer essential, although they can be used to set a tuple
apart from surrounding text. Since parentheses are no longer part
of the tuple construct, we can regard a tuple with only one element
as that element\emph{.}

Parentheses still have the usual function of indicating in what order
(concatenation) operations are performed. For example, $\left(x,y\right),z$
is the tuple $x,y$ concatenated with the element $z$, while $x,\left(y,z\right)$
is the element $x$ concatenated with the tuple $y,z$. On the other
hand, concatenation of tuples is clearly an associative operation;
for arbitrary tuples $\boldsymbol{x}$, $\boldsymbol{y}$ and $\boldsymbol{z}$
we have
\[
\boldsymbol{x},\left(\boldsymbol{y},\boldsymbol{z}\right)=\boldsymbol{x},\boldsymbol{y},\boldsymbol{z}=\left(\boldsymbol{x},\boldsymbol{y}\right),\boldsymbol{z}.
\]

\emph{(2)}. A \emph{tuple set} is a set of $N$-tuples. Let $X^{m}$
be a set of $m$-tuples and $Y^{n}$ a set of $n$-tuples; the corresponding
\emph{Cartesian product of tuple sets} is
\[
X^{m}\times Y^{m}=\left\{ \boldsymbol{x},\boldsymbol{y}\mid\boldsymbol{x}\in X^{m},\boldsymbol{y}\in Y^{n}\right\} .
\]
When $X^{m}=X$ and $Y^{m}=Y$ are sets of $1$-tuples, $X^{m}\times Y^{m}$
can be regarded as the usual Cartesian product
\[
X\times Y=\left\{ \left(x,y\right)\mid x\in X,y\in Y\right\} .
\]
Note, though, that while the usual Cartesian product is not associative
because $\left(\boldsymbol{x},\left(\boldsymbol{y},\boldsymbol{z}\right)\right)\neq\left(\left(\boldsymbol{x},\boldsymbol{y}\right),\boldsymbol{z}\right)$,
the Cartesian product of tuple sets is associative because concatenation
of tuples is associative, 

The Cartesian product will be understood as the Cartesian product
of tuple sets in this article, and this product will be denoted in
the same way as the usual Cartesian product.

\specialsection*{\textbf{\textsl{A. Tensor products and related notions}}}

\section{On multiplication of vectors}

\subsection{Algebras, quotient algebras and monoid algebras}

$ $\\
\emph{(1)}. An \emph{algebra} $A$ over a field $K$ can be defined
as a vector space over $K$ equipped with a \emph{$K$-bilinear function}
called \emph{multiplication} of vectors\emph{,} 
\[
A\times A\rightarrow A,\qquad\left(x,y\right)\mapsto xy,
\]
so that $\left(x+x'\right)\! y=xy+x'y$, $x\!\left(y+y'\right)=xy+xy'$,
and $\left(kx\right)\! y=x\!\left(ky\right)=k\!\left(xy\right)$ for
every $x,x',y,y'\in A$, $k\in K$. $0_{A}$ denotes the zero element
in $A$.

Only \emph{unital associative algebras} will be considered here; these
are algebras which are associative under multiplication and have a
unit element $1_{\! A}$.

\emph{(2)}. An \emph{ideal} $I$ in $A$ is a subspace of $A$ (as
a vector space) such that $\iota x,x\iota\in I$ for any $\iota\in I,x\in A$.
We denote the set $\left\{ x+\iota\mid\iota\in I\right\} $ by $x+I$
or $\left[x\right]$. 

$I$ is a subgroup of the additive group of $A$, so $A/I=\left\{ \left[x\right]\mid x\in A\right\} $
is a partition of $A$. As $\lambda\left(x+I\right)=\lambda x+\lambda I=\lambda x+I$,
$\left(x+I\right)+\left(y+I\right)=x+y+I+I=x+y+I$, and $\left(x+I\right)\left(y+I\right)=xy+xI+Iy+II=xy+I$,
$A/I$ can be made into a unital associative algebra with operations
defined by 
\[
\lambda\left[x\right]=\left[\lambda x\right],\quad\left[x\right]+\left[y\right]=\left[x+y\right],\quad\left[x\right]\left[y\right]=\left[xy\right].
\]
$A/I$ is said to be a \emph{quotient algebra}. $\left[0_{A}\right]$
is obviously the zero element in $A/I$, while $\left[1_{A}\right]$
is the unit element.

For any set $S\subset A$ there is a unique \emph{smallest ideal}
$I\!\left(S\right)\supset S$, and a corresponding quotient algebra
$A/I\!\left(S\right)$. Alternatively, we can define a quotient algebra
by specifying a set $R$ of \emph{relations} of the form $x=y$, where
$x,y\in A$. For each relation $x=y$ there is a corresponding element
$x-y\in S$, and for each $z\in S$ there is a corresponding relation
$z=0$. In other words, there is a one-to-one-correspondence $S\leftrightarrow R$,
and we can write $A/I\!\left(S\right)$ as $A/I\!\left(R\right)$. 

Note that if $x-y\in I$ then $\left[x\right]-\left[y\right]=\left[x-y\right]=0_{A}+x-y+I=\left[0_{A}\right]$.
Hence, the relation $x=y$ is reflected by the identity $\left[x\right]=\left[y\right]$
in $A/I\!\left(R\right)$. Expressed in another way, $\left[x\right]=\left[y\right]$
if and only if $x$ can be rewritten as $y$ by using relations in
$R$ $n_{R}\geq0$ times and identities in $A$ $n_{A}\geq0$ times.

\emph{(3)}. Recall that a \emph{monoid} $M$ is a set with associative
multiplication $\left(x,y\right)\mapsto xy$ and an identity element.
Let $K$ be a field (such as the real numbers) and consider the set
$K\!\left[M\right]$ of functions $\phi:M\rightarrow K$ such that
$\left\{ x\mid\phi\left(x\right)\neq0\right\} $ is a (possibly empty)
finite set. Let addition in $K\!\left[M\right]$ be defined by $\left(\alpha+\beta\right)\!\left(x\right)=\alpha\!\left(x\right)+\beta\!\left(x\right)$
and scalar multiplication by $\left(k\phi\right)\!\left(x\right)=k\!\left(\phi\!\left(x\right)\right)$,
where $k\in K$. $\left\{ x\mid\left(\alpha+\beta\right)\!\left(x\right)\neq0\right\} $
and $\left\{ x\mid\left(k\phi\right)\!\left(x\right)\neq0\right\} $
are clearly finite sets, and $K\!\left[M\right]$ is a vector space
over $K$. The zero element $0_{K\left[M\right]}$ in $K\!\left[M\right]$
is the function $x\mapsto0$.

As $\left\{ \left(x,y\right)\mid\alpha\!\left(x\right)\beta\!\left(y\right)\neq0\right\} $
is finite, we can define multiplication in $K\!\left[M\right]$ by
\[
\alpha\beta\!\left(z\right)=\sum_{xy=z}\alpha\!\left(x\right)\beta\!\left(y\right);
\]
the map $\left(\alpha,\beta\right)\mapsto\alpha\beta$ is clearly
bilinear. Let $e$ be the identity element in $M$, and define a function
$1_{K\left[M\right]}:M\rightharpoondown K$ by $1_{K\left[M\right]}\!\left(x\right)=1$
if $x=e$ and $1_{K\left[M\right]}\!\left(x\right)=0$ if $x\neq e$.
Then
\[
1_{K\left[M\right]}\,\phi\!\left(z\right)=\sum_{xy=z}1_{K\left[M\right]}\!\left(x\right)\,\phi\!\left(y\right)=1_{K\left[M\right]}\!\left(e\right)\,\phi\!\left(z\right)=\phi\!\left(z\right),
\]
so $1_{K\left[M\right]}\,\phi=\phi$ and similarly $\phi\,1_{K\left[M\right]}=\phi$,
so $1_{K\left[M\right]}$ is indeed the identity in $K\!\left[M\right]$.
Furthermore,
\begin{gather*}
\left(\alpha\beta\right)\!\gamma\!\left(z\right)=\sum_{xy=z}\left(\alpha\beta\right)\!\left(x\right)\gamma\!\left(y\right)=\sum_{xy=z}\left(\sum_{uv=x}\alpha\!\left(u\right)\beta\!\left(v\right)\right)\gamma\!\left(y\right)\\
=\sum_{xy=z}\sum_{uv=x}\alpha\!\left(u\right)\beta\!\left(v\right)\gamma\!\left(y\right)=\sum_{\left(uv\right)y=z}\alpha\!\left(u\right)\beta\!\left(v\right)\gamma\!\left(y\right),
\end{gather*}
and similarly 
\[
\alpha\!\left(\beta\gamma\right)\!\left(z\right)=\sum_{u\left(vy\right)=z}\alpha\!\left(u\right)\beta\!\left(v\right)\gamma\!\left(y\right),
\]
so $\left(\alpha\beta\right)\!\gamma=\alpha\!\left(\beta\gamma\right)$
since $\left(uv\right)y=u\left(vy\right)$ in a monoid. Thus, $K\!\left[M\right]$
is a unital associative algebra, called a \emph{monoid algebra}.

\emph{(4)}. Let $\triangle\!\left(x\right):M\rightarrow K$ be a function
such that $\triangle\!\left(x\right)\!\left(y\right)=1$ if $x=y$
and $\triangle\!\left(x\right)\!\left(y\right)=0$ if $x\neq y$.
By definition, $\triangle\!\left(e\right)=1_{K\left[M\right]}$. Any
non-zero $\phi\in K\!\left[M\right]$ clearly has a unique non-empty
finite-sum expansion of the form 
\begin{equation}
\phi=\sum_{k_{x}\neq0}k_{x}\,\triangle\!\left(x\right),\quad\mathrm{where\,}\; k_{x}=\phi\left(x\right),\label{eq:expansionMonAlg}
\end{equation}
and $0_{K\!\left[M\right]}$ is uniquely represented by the empty
sum $\textrm{Ø}$, since $0_{K\!\left[M\right]}\!\left(x\right)=0$
for all $x\!\in\! M$. $\left\{ \triangle\!\left(x\right)\mid x\in M\right\} $
is thus a basis for $K\!\left[M\right]$, and in terms of this basis
we have 
\begin{gather*}
\lambda\phi=\lambda\left(\sum_{k_{x}\neq0}k_{x}\,\triangle\!\left(x\right)\right)=\sum_{\lambda k_{x}\neq0}\lambda k_{x}\,\triangle\!\left(x\right),\\
\alpha+\beta=\sum_{a_{x}\neq0}a_{x}\,\triangle\!\left(x\right)+\sum_{b_{x}\neq0}b_{x}\,\triangle\!\left(x\right)=\sum_{\left(a_{x}+b_{x}\right)\neq0}\left(a_{x}+b_{x}\right)\,\triangle\!\left(x\right),
\end{gather*}
for any $\lambda\in K$ and $\phi,\alpha,\beta\in K\!\left[M\right]$.

\subsection{The monoid algebras \textmd{$K\!\left[M\!\left(V_{1},\ldots,V_{n}\right)\right]$}
and\textmd{ $K\!\left[M\!\left(V_{1},\ldots,V_{n}\right)\right]/I\!\left(\mathsf{L}\right)$} }

$ $\\
\emph{(1)}. Consider a vector space $V$ over $K$ and the set $M\!\left(V\right)$
of all tuples of vectors in $V$. Let multiplication of elements of
$M\!\left(V\right)$ be concatenation of tuples, $\left(\mathbf{u}_{1},\ldots,\mathbf{u}_{m}\right)\left(\mathbf{v}_{1},\ldots,\mathbf{v}_{n}\right)=\left(\mathbf{u}_{1},\ldots,\mathbf{u}_{m},\mathbf{v}_{1},\ldots,\mathbf{v}_{n}\right)$.
Concatenation is associative, so $M\!\left(V\right)$ is a monoid
with the empty tuple $\left(\right)$ as the identity element. 

Given $M\!\left(V\right)$, we can construct the monoid algebra $K\!\left[M\!\left(V\right)\right]$,
with tuples of the form $\triangle\!\left(\mathbf{v}_{1},\ldots,\mathbf{v}_{k}\right)$,
or $\triangle\!\left(\right)$ if $k=0$. By the definitions of multiplication
in $K\!\left[M\right]$ and $M\!\left(V\right)$, 
\[
\triangle\!\left(\mathbf{u}\right)\!\triangle\!\left(\mathbf{v}\right)\!\left(\mathbf{x},\mathbf{y}\right)=\left.\sum\right._{\left(\mathbf{s}\right)\left(\mathbf{t}\right)=\left(\mathbf{x},\mathbf{y}\right)}\triangle\!\left(\mathbf{u}\right)\!\left(\mathbf{s}\right)\triangle\!\left(\mathbf{v}\right)\!\left(\mathbf{t}\right)=\triangle\!\left(\mathbf{u}\right)\!\left(\mathbf{x}\right)\triangle\!\left(\mathbf{v}\right)\!\left(\mathbf{y}\right),
\]
and clearly $\triangle\!\left(\mathbf{u}\right)\!\left(\mathbf{x}\right)\triangle\!\left(\mathbf{v}\right)\!\left(\mathbf{y}\right)=\triangle\!\left(\mathbf{u},\mathbf{v}\right)\!\left(\mathbf{x},\mathbf{y}\right)$,
so $\triangle\!\left(\mathbf{u}\right)\triangle\!\left(\mathbf{v}\right)=\triangle\!\left(\mathbf{u},\mathbf{v}\right)$.
More generally, we have 
\[
\triangle\!\left(\mathbf{u}_{1},\ldots,\mathbf{u}_{m}\right)\,\triangle\!\left(\mathbf{v}_{1},\ldots,\mathbf{v}_{n}\right)=\triangle\!\left(\mathbf{u}_{1},\ldots,\mathbf{u}_{m},\mathbf{v}_{1},\ldots,\mathbf{v}_{n}\right).
\]

By (\ref{eq:expansionMonAlg}), every $\phi\in K\!\left[M\!\left(V\right)\right]$
has a unique possibly empty finite-sum expansion
\[
\phi=\left.\sum\right._{i}k_{i}\triangle\!\left(\mathbf{v}_{1i},\ldots,\mathbf{v}_{n_{i}i}\right)=\left.\sum\right._{i}k_{i}\triangle\!\left(\mathbf{v}_{1i}\right)\ldots\triangle\!\left(\mathbf{v}_{n_{i}i}\right),
\]
where $k_{i}\neq0,n_{i}\geq0$ and $\mathbf{v}_{1i},\ldots,\mathbf{v}_{n_{i}i}\in V$.
Collecting terms, $\phi$ has a unique expansion as a sum of the form
\begin{equation}
\phi=\left.\sum\right._{n\in\mathbb{N}_{0}}\left.\sum\right._{i=1}^{m_{n}}k_{ni}\triangle\!\left(\mathbf{v}_{1i}\right)\ldots\triangle\!\left(\mathbf{v}_{ni}\right),\label{eq:expanFi1}
\end{equation}
where $\mathbb{N}_{0}$ is a possibly empty finite set of non-negative
integers, $m_{0}\!=\!1$,\linebreak{}
 $k_{ni}\!\neq\!0,\mathbf{v}_{1i},\ldots,\mathbf{v}_{ni}\in V$ and
$\sum_{i=1}^{m_{0}}k_{0i}\triangle\!\left(\mathbf{v}_{1i}\right)\ldots\triangle\!\left(\mathbf{v}_{0i}\right)=k\triangle\!\left(\right)$.
This means that
\[
\left\{ \triangle\!\left(\right)\right\} \cup\left\{ \triangle\!\left(\mathbf{v}_{i}\right)\right\} \cup\left\{ \triangle\!\left(\mathbf{v}_{i}\right)\triangle\!\left(\mathbf{v}_{j}\right)\right\} \cup\ldots\qquad\left(\mathbf{v}_{i},\mathbf{v}_{j},\ldots\in V\right)
\]
is a basis for $K\!\left[M\!\left(V\right)\right]$. Expressed in
terms of this basis, the zero element $0_{K\!\left[M\!\left(V\right)\right]}$
is the empty sum $\textrm{Ø}$, while the identity element $1_{K\!\left[M\!\left(V\right)\right]}$
is the one-term sum $\triangle\!\left(\right)$. 

\emph{(2)}. Consider the set of relations in $K\!\left[M\!\left(V\right)\right]$:
\[
\mathsf{L}=\left\{ \,\triangle\!\left(\lambda\mathbf{v}\right)=\lambda\triangle\!\left(\mathbf{v}\right),\quad\triangle\!\left(\mathbf{v}+\mathbf{v}'\right)=\triangle\!\left(\mathbf{v}\right)+\triangle\!\left(\mathbf{v}'\right)\,\right\} .
\]
Because multiplication in $K\!\left[M\!\left(V\right)\right]$ is
bilinear, these two relations imply that 
\begin{gather*}
\triangle\!\left(\lambda\mathbf{u}\right)\triangle\!\left(\mathbf{v}\right)=\triangle\!\left(\mathbf{u}\right)\triangle\!\left(\lambda\mathbf{v}\right)=\lambda\left(\triangle\!\left(\mathbf{u}\right)\triangle\!\left(\mathbf{v}\right)\right),\\
\triangle\!\left(\mathbf{u}+\mathbf{u}'\right)\triangle\!\left(\mathbf{v}\right)=\triangle\!\left(\mathbf{u}\right)\triangle\!\left(\mathbf{v}\right)+\triangle\!\left(\mathbf{u}'\right)\triangle\!\left(\mathbf{v}\right),\\
\triangle\!\left(\mathbf{u}\right)\triangle\!\left(\mathbf{v}+\mathbf{v}'\right)=\triangle\!\left(\mathbf{u}\right)\triangle\!\left(\mathbf{v}\right)+\triangle\!\left(\mathbf{u}\right)\triangle\!\left(\mathbf{v}'\right),
\end{gather*}
and analogously for functions of the form $\triangle\!\left(\mathbf{v}_{1}\right)\ldots\triangle\!\left(\mathbf{v}_{n}\right)$.

Corresponding to $\mathsf{L}$, there is a quotient algebra $K\!\left[M\!\left(V\right)\right]/I\!\left(\mathsf{L}\right)$;
its elements are equivalence classes of the form $\left[\phi\right]$,
where $\phi\in K\!\left[M\!\left(V\right)\right]$. $K\!\left[M\!\left(V\right)\right]/I\!\left(\mathsf{L}\right)$
is a unital associative algebra with identity element $\left[1_{K\!\left[M\!\left(V\right)\right]}\right]=\left[\triangle\!\left(\right)\right]$.

\emph{(3)}. For any $\mathbf{v}_{1},\ldots,\mathbf{v}_{n}\in V$,
we have 
\[
\left[\left.\sum\right._{i}\lambda_{i}\triangle\!\left(\mathbf{v}_{i}\right)\right]=\left[\triangle\!\left(\left.\sum\right._{i}\lambda_{i}\mathbf{v}_{i}\right)\right]
\]
in view of the relations in $\mathsf{L}$, and hence
\[
\left.\sum\right._{i}\alpha_{i}\mathbf{v}_{i}=\left.\sum\right._{i}\beta_{i}\mathbf{v}_{i}\quad\Longrightarrow\quad\left[\left.\sum\right._{i}\alpha_{i}\triangle\!\left(\mathbf{v}_{i}\right)\right]=\left[\left.\sum\right._{i}\beta_{i}\triangle\!\left(\mathbf{v}_{i}\right)\right].
\]
Conversely, reflection on the relations in $\mathsf{L}$ leads to
the conclusion that we can have $\left[\left.\sum\right._{i}\alpha_{i}\triangle\!\left(\mathbf{v}_{i}\right)\right]=\left[\left.\sum\right._{i}\beta_{i}\triangle\!\left(\mathbf{v}_{i}\right)\right]$
only if (a) $\left.\sum\right._{i}\alpha_{i}\triangle\!\left(\mathbf{v}_{i}\right)=\left.\sum\right._{i}\beta_{i}\triangle\!\left(\mathbf{v}_{i}\right)$
or (b) $\triangle\!\left(\left.\sum\right._{i}\alpha_{i}\mathbf{v}_{i}\right)=\triangle\!\left(\left.\sum\right._{i}\beta_{i}\mathbf{v}_{i}\right)$.
Let $\mathbf{v}_{1},\ldots,\mathbf{v}_{n}$ be \emph{distinct} vectors
in $V$. In case (a), $\alpha_{i}=\beta_{i}$ for all $i$, since
expansions of the form $\varphi=\left.\sum\right._{i}\lambda_{i}\triangle\!\left(\mathbf{v}_{i}\right)$
are unique. In case (b), $\left.\sum\right._{i}\alpha_{i}\mathbf{v}_{i}=\left.\sum\right._{i}\beta_{i}\mathbf{v}_{i}$,
because $\triangle\!\left(\mathbf{u}\right)=\triangle\!\left(\mathbf{v}\right)$
clearly implies $\mathbf{u}=\mathbf{v}$. Hence,
\[
\left[\left.\sum\right._{i}\alpha_{i}\triangle\!\left(\mathbf{v}_{i}\right)\right]=\left[\left.\sum\right._{i}\beta_{i}\triangle\!\left(\mathbf{v}_{i}\right)\right]\quad\Longrightarrow\quad\left.\sum\right._{i}\alpha_{i}\mathbf{v}_{i}=\left.\sum\right._{i}\beta_{i}\mathbf{v}_{i}.
\]

As $\lambda\left[x\right]=\left[\lambda x\right]$ and $\left[x\right]+\left[y\right]=\left[x+y\right]$,
we also have
\[
\left.\sum\right._{i}\lambda_{i}\left[\triangle\!\left(\mathbf{v}_{i}\right)\right]=\left[\left.\sum\right._{i}\lambda_{i}\triangle\!\left(\mathbf{v}_{i}\right)\right].
\]

Together, these facts imply that if $\left\{ \mathbf{e}_{i}\right\} $
is a basis for $V$ and $\left.\sum\right._{i}\alpha{}_{i}\left[\triangle\!\left(\mathbf{e}_{i}\right)\right]=\left.\sum\right._{i}\beta_{i}\left[\triangle\!\left(\mathbf{e}_{i}\right)\right]$
then $\alpha_{i}=\beta_{i}$ for every $i$. Thus, any $\left[\varphi\right]=\left[\left.\sum\right._{i}k_{i}\triangle\!\left(\mathbf{v}_{i}\right)\right]=\left[\triangle\!\left(\mathbf{v}\right)\right]=\left[\left.\sum\right._{i}\kappa{}_{i}\triangle\!\left(\mathbf{e}_{i}\right)\right]$
has a unique possibly empty finite-sum expansion
\[
\left[\varphi\right]=\left.\sum\right._{i}\kappa{}_{i}\left[\triangle\!\left(\mathbf{e}_{i}\right)\right]\qquad\left(\kappa_{i}\neq0\right).
\]

\emph{(4)}. As every $\phi\in K\!\left[M\!\left(V\right)\right]$
has an expansion of the form (\ref{eq:expanFi1}) and\linebreak{}
 $\left[\triangle\!\left(\mathbf{v}_{1}\right)\ldots\triangle\!\left(\mathbf{v}_{n}\right)\right]=\left[\triangle\!\left(\mathbf{v}_{1}\right)\right]\ldots\left[\triangle\!\left(\mathbf{v}_{n}\right)\right]$,
every $\left[\phi\right]\in K\!\left[M\!\left(V\right)\right]/I\!\left(\mathsf{L}\right)$
has a corresponding expansion of the form 
\begin{equation}
\left[\phi\right]=\left.\sum\right._{n\in\mathbb{N}_{0}}\left.\sum\right._{i=1}^{m_{n}}k_{ni}\left[\triangle\!\left(\mathbf{v}_{1i}\right)\right]\ldots\left[\triangle\!\left(\mathbf{v}_{ni}\right)\right].\label{eq:expanFi2}
\end{equation}
Let $\left\{ \mathbf{e}_{i}\right\} $ be a basis for $V$ and consider
the set
\[
\mathcal{B}=\left\{ \left[\triangle\!\left(\right)\right]\right\} \cup\left\{ \left[\triangle\!\left(\mathbf{e}_{i}\right)\right]\right\} \cup\left\{ \left[\triangle\!\left(\mathbf{e}_{i}\right)\right]\left[\triangle\!\left(\mathbf{e}_{j}\right)\right]\right\} \cup\ldots\;.
\]
Since 
\begin{gather*}
\left[\triangle\!\left(\mathbf{v}_{1i}\right)\right]\ldots\left[\triangle\!\left(\mathbf{v}_{ni}\right)\right]=\left[\triangle\!\left(\left.\sum\right._{j_{1}}\lambda_{1ij_{1}}\mathbf{e}_{j_{1}}\right)\right]\ldots\left[\triangle\!\left(\left.\sum\right._{j_{n}}\lambda_{nij_{n}}\mathbf{e}_{j_{n}}\right)\right]\\
=\left.\sum\right._{j_{1},\ldots,j_{n}}\lambda_{1ij_{1}}\ldots\lambda_{nij_{n}}\left[\triangle\!\left(\mathbf{e}_{j_{1}}\right)\right]\ldots\left[\triangle\!\left(\mathbf{e}_{j_{n}}\right)\right],
\end{gather*}
$\left[\phi\right]$ has an expansion
\[
\left[\phi\right]=\left.\sum\right._{n\in\mathbb{N}_{0}}\left.\sum\right._{i=1}^{m_{n}}\kappa_{ni}\left[\triangle\!\left(\mathbf{e}_{1i}\right)\right]\ldots\left[\triangle\!\left(\mathbf{e}_{ni}\right)\right]
\]
in terms of $\mathcal{B}$. It can be shown by an argument similar
to that in the previous subsection that this expansion is unique,
so $\mathcal{B}$ is a basis for $K\!\left[M\!\left(V\right)\right]/I\!\left(\mathsf{L}\right)$. 

\emph{(5)}. Suppose that we start from tuples of the form
\[
\left(\mathbf{v}_{1},\ldots,\mathbf{v}_{m}\right),\quad\mathrm{where}\quad\mathbf{v}_{i}\in V_{j_{i}}\in\left\{ V_{1},\ldots,V_{n}\right\} 
\]
instead of tuples of the form $\left(\mathbf{v}_{1},\ldots,\mathbf{v}_{m}\right)$,
where $\mathbf{v}_{i}\in V$. Using such tuples, we can define unital
associative algebras 
\[
K\!\left[M\!\left(V_{1},\ldots,V_{n}\right)\right]\quad\mathrm{and}\quad K\!\left[M\!\left(V_{1},\ldots,V_{n}\right)\right]/I\!\left(\mathsf{L}\right)
\]
in essentially the same way as $K\!\left[M\!\left(V\right)\right]$
and $K\!\left[M\!\left(V\right)\right]/I\!\left(\mathsf{L}\right)$,
and then obtain generalizations of all results.

In particular, it can be shown that if $\left\{ \mathbf{e}_{i}\right\} $
is a basis for $U$ and $\left\{ \mathbf{f}_{i}\right\} $ is a basis
for $V$, then
\begin{gather*}
\left\{ \left[\triangle\!\left(\right)\right]\right\} \cup\left\{ \left[\triangle\!\left(\mathbf{e}_{i}\right)\right]\right\} \cup\left\{ \left[\triangle\!\left(\mathbf{f}_{i}\right)\right]\right\} \cup\\
\left\{ \left[\triangle\!\left(\mathbf{e}_{i}\right)\right]\left[\triangle\!\left(\mathbf{e}_{j}\right)\right]\right\} \cup\left\{ \left[\triangle\!\left(\mathbf{e}_{i}\right)\right]\left[\triangle\!\left(\mathbf{f}_{j}\right)\right]\right\} \cup\left\{ \left[\triangle\!\left(\mathbf{f}_{i}\right)\right]\left[\triangle\!\left(\mathbf{e}_{j}\right)\right]\right\} \cup\left\{ \left[\triangle\!\left(\mathbf{f}_{i}\right)\right]\left[\triangle\!\left(\mathbf{f}_{j}\right)\right]\right\} \cup\ldots
\end{gather*}
is a basis for $K\!\left[M\!\left(U,V\right)\right]/I\!\left(\mathsf{L}\right)$.
This generalizes in an obvious (but complicated) way to the general
case with $n$ vector spaces involved.

\subsection{Free unital associative algebras on vector spaces}

\textbf{$ $}\\
\emph{(1)}. Let $V$ be a vector space over $K$. A \emph{free} unital
associative algebra $\mathcal{A}\!\left(V\right)$ on $V$ is defined
here as a unital associative algebra over $K$ which includes a copy
of $V$ and has the property that ($\star$) if $\left\{ \mathbf{e}_{i}\right\} $
is a basis for $V$ then the infinite union 
\[
\mathcal{E}=\left\{ 1\right\} \cup\left\{ \mathbf{e}_{i}\right\} \cup\left\{ \mathbf{e}_{i}\mathbf{e}_{j}\right\} \cup\left\{ \mathbf{e}_{i}\mathbf{e}_{j}\mathbf{e}_{k}\right\} \cup\ldots\qquad\left(1\in K\right),
\]
where distinct products $\mathbf{e}_{i_{1}}\ldots\mathbf{e}_{i_{n}}$
denote distinct vectors, is a basis for $\mathcal{A}\!\left(V\right)$.
The elements of $\mathcal{A}\!\left(V\right)$ can thus be written
as finite sums (polynomials) of the form 
\[
t+\sum_{i}t^{i}\mathbf{e}_{i}+\sum_{i,j}t^{ij}\mathbf{e}_{i}\mathbf{e}_{j}+\ldots+\sum_{i_{1},\ldots,i_{n}}t^{i_{1}\ldots i_{n}}\mathbf{e}_{i_{1}}\ldots\mathbf{e}_{i_{n}}.
\]

To construct a free unital associative algebra on $V$, we perform
'surgery' on $K\!\left[M\!\left(V\right)\right]/I\!\left(\mathsf{L}\right)$.
Set $\mathcal{S}\!\left(V\right)=\left\{ x\mid x\in K\!\left[M\!\left(V\right)\right]/I\!\left(\mathsf{L}\right)\right\} $,
and let $\mathcal{S}'\!\left(V\right)$ be the same set except that
$\left[\triangle\!\left(\right)\right]$ is replaced by $1\in K$
and $\left[\triangle\!\left(\mathbf{v}\right)\right]$ by $\mathbf{v}$
for all $\mathbf{v}\in V$. Let $\boldsymbol{\iota}$ be the bijection
$\mathcal{S}\!\left(V\right)\rightarrow\mathcal{S}'\!\left(V\right)$
which is equal to the identity map on $\mathcal{S}\!\left(V\right)$
except that $\boldsymbol{\iota}\!\left(\left[\triangle\!\left(\right)\right]\right)=1$
and $\boldsymbol{\iota}\!\left(\left[\triangle\!\left(\mathbf{v}\right)\right]\right)=\mathbf{v}$
for all $\mathbf{v}\in V$. Regard $\boldsymbol{\iota}$ as a mapping
$K\!\left[M\!\left(V\right)\right]/I\!\left(\mathsf{L}\right)\rightarrow\mathcal{S}'\!\left(V\right)$,
and define scalar multiplication, addition and multiplication of elements
$\mathbf{s},\mathbf{t}$ of $S'\!\left(V\right)$ by $\lambda\mathbf{s}=\boldsymbol{\iota}\!\left(\lambda\boldsymbol{\iota^{-1}}\!\mathbf{\left(\mathbf{s}\right)}\right)$,
$\mathbf{s}+\mathbf{t}=\boldsymbol{\iota}\!\left(\boldsymbol{\iota}^{-1}\!\mathbf{\left(\mathbf{s}\right)}+\boldsymbol{\iota}^{-1}\!\mathbf{\left(\mathbf{t}\right)}\right)$
and $\mathbf{s}\mathbf{t}=\boldsymbol{\iota}\!\left(\boldsymbol{\iota}^{-1}\!\mathbf{\left(\mathbf{s}\right)}\,\boldsymbol{\iota}^{-1}\!\mathbf{\left(\mathbf{t}\right)}\right)$;
let $\mathcal{S}^{*}\!\left(V\right)$ denote $S'\!\left(V\right)$
with these operations. By design, $\boldsymbol{\iota}$ regarded as
a mapping $K\!\left[M\!\left(V\right)\right]/I\!\left(\mathsf{L}\right)\rightarrow S{}^{*}\!\left(V\right)$
is an isomorphism, and we recover the vector space operations in $V$
since $\boldsymbol{\iota}\!\left(\lambda\boldsymbol{\iota^{-1}}\!\mathbf{\left(\mathbf{v}\right)}\right)=\boldsymbol{\iota}\!\left(\lambda\left[\triangle\!\left(\mathbf{v}\right)\right]\right)=\boldsymbol{\iota}\!\left(\left[\triangle\!\left(\lambda\mathbf{v}\right)\right]\right)=\lambda\mathbf{v}$
and $\boldsymbol{\iota}\!\left(\boldsymbol{\iota}^{-1}\!\mathbf{\left(\mathbf{u}\right)}+\boldsymbol{\iota}^{-1}\!\mathbf{\left(\mathbf{v}\right)}\right)=\boldsymbol{\iota}\!\left(\left[\triangle\!\left(\mathbf{u}\right)\right]+\left[\triangle\!\left(\mathbf{v}\right)\right]\right)=\boldsymbol{\iota}\!\left(\left[\triangle\!\left(\mathbf{u}+\mathbf{v}\right)\right]\right)=\mathbf{u}+\mathbf{v}$.
Together with $\boldsymbol{\iota}$, the main result of Subsection
3.2(4) implies that condition ($\star$) is satisfied. We denote $S^{*}\!\left(V\right)$
by $\mathcal{A}\!\left(V\right)$ and call it \emph{the} free unital
associative algebra on $V$.

\emph{(2)}. In view of the close analogy between the algebras $K\!\left[M\!\left(V\right)\right]/I\!\left(\mathsf{L}\right)$
and $K\!\left[M\!\left(V_{1},\ldots,V_{n}\right)\right]/I\!\left(\mathsf{L}\right)$
noted in Subsection 3.2(5), the construction of the free associative
algebra $\mathcal{A}\!\left(V\right)$ can be generalized to a construction
of the free unital associative algebra on two or more vector spaces,
$\mathcal{A}\!\left(V_{1},\ldots,V_{n}\right)$. In particular, the
free unital associative algebra on $U$ and $V$, denoted $\mathcal{A}\!\left(U,V\right)$,
is a unital associative algebra such that if $\left\{ \mathbf{e}_{i}\right\} $
is a basis for $U$ and $\left\{ \mathbf{f}_{j}\right\} $ is a basis
for $V$ then
\[
\left\{ 1\right\} \cup\left\{ \mathbf{e}_{i}\right\} \cup\left\{ \mathbf{f}_{j}\right\} \cup\left\{ \mathbf{e}_{i}\mathbf{e}_{j}\right\} \cup\left\{ \mathbf{e}_{i}\mathbf{f}_{j}\right\} \cup\left\{ \mathbf{f}_{i}\mathbf{e}_{j}\right\} \cup\left\{ \mathbf{f}_{i}\mathbf{f}_{j}\right\} \cup\ldots\,,
\]
where distinct expressions denote distinct vectors, is a basis for
for $\mathcal{A}\left(U,V\right)$.

\section{Tensor products of two vectors or vector spaces}

\subsection{Tensor product maps and related tensor products}

$ $\\
\emph{(1)}.~~Let $U,V$ be vector spaces, and suppose that there
exists a bilinear\emph{ }map 
\[
\mu_{\otimes}\left(U,V\right):U\times V\rightarrow W,\quad\left(\mathbf{u},\mathbf{v}\right)\mapsto\mathbf{u}\otimes\mathbf{v}
\]
such that ($\otimes$) if $\left\{ \mathbf{e}_{1},\ldots,\mathbf{e}_{m}\right\} $
is a basis for $U$ and $\left\{ \mathbf{f}{}_{1},\ldots,\mathbf{f}{}_{n}\right\} $
a basis for $V$ then the map $\left(\mathbf{e}_{i},\mathbf{f}{}_{j}\right)\mapsto\mathbf{e}_{i}\otimes\mathbf{f}{}_{j}$
is injective and $\left\{ \mathbf{e}_{i}\otimes\mathbf{f}{}_{j}\mid i=1,\ldots,m;\, j=1,\ldots,n\right\} $
is a basis for $W$. We call $\mu_{\otimes}\left(U,V\right)$ a\emph{
tensor product map} and $W$ \emph{a} \emph{tensor product} of \emph{$U$}
and \emph{$V$ }or \emph{the tensor product} of $U$ and $V$ \emph{for
$\mu_{\otimes}\left(U,V\right)$}. $W$ is usually denoted by $U\otimes V$,
and the elements of \emph{$U\otimes V$ }are called\emph{ tensors.} 

Clearly, $\dim\left(U\otimes V\right)=\dim\left(U\right)\dim\left(V\right)$.
As $\mu_{\otimes}$ is bilinear, we have
\begin{gather*}
\eta\mathbf{u}\otimes\mathbf{v}=\eta\left(\mathbf{u}\otimes\mathbf{v}\right),\quad\left(\mathbf{u}_{1}+\mathbf{u}_{2}\right)\otimes\mathsf{v}=\mathbf{u}_{1}\otimes\mathbf{v}+\mathbf{u}_{2}\otimes\mathbf{v},\\
\mathbf{u}\otimes\eta\mathbf{v}=\eta\left(\mathbf{u}\otimes\mathbf{v}\right),\quad\mathbf{u}\otimes\left(\mathbf{v}_{1}+\mathbf{v}_{2}\right)=\mathbf{u}\otimes\mathbf{v}_{1}+\mathbf{u}\otimes\mathbf{v}_{2}.
\end{gather*}

One should not assume that $\mathbf{u}_{1}\otimes\mathbf{v}_{1}=\mathbf{u}_{2}\otimes\mathbf{v}_{2}$
implies $\left(\mathbf{u}_{1},\mathbf{v}_{1}\right)=\left(\mathbf{u}_{2},\mathbf{v}_{2}\right)$,
or that every element of $U\otimes V$ can be written in the form
$\mathbf{u}\otimes\mathbf{v}$, where $\mathbf{u}\in U,\mathbf{v}\in V$.
In general, a tensor product map is neither injective nor surjective.

An element of $U\otimes V$ which can be written as $\mathbf{u}\otimes\mathbf{v}$
is said to be a \emph{simple tensor}. We denote the set of simple
tensors in $U\otimes V$ by $\left.U\otimes V\right|_{\mathcal{S}}$.
\begin{rem*}
There is some abuse of notation here, since '$\otimes$' does not
refer to a specific binary operation. This symbol is not used in the
same sense in $\mathbf{u}\otimes\mathbf{v}$ as in $U\otimes V$,
and it can be associated with different tensor product maps in different
contexts.
\end{rem*}
\emph{(2)}.~Recall that $K$ is a vector space over itself with bases
of the form $\left\{ \eta\right\} $, where $\eta\neq0$. Let $V$
be a vector space over $K$ and consider the bilinear functions
\begin{gather*}
\mu_{\otimes}\left(K,V\right):K\times V\rightarrow V,\qquad\left(\alpha,\mathbf{v}\right)\mapsto\alpha\otimes\mathbf{v}=\alpha\mathbf{v},\\
\mu_{\otimes}\left(V,K\right):V\times K\rightarrow V,\qquad\left(\mathbf{v},\beta\right)\mapsto\mathbf{v}\otimes\beta=\beta\mathbf{v}.
\end{gather*}
If $\left\{ \eta\right\} $ is a basis for $K$, so that $\eta\neq0$,
and $\left\{ \mathbf{e}_{i}\right\} $ a basis for $V$, then \linebreak{}
$\left(\eta,\mathbf{e}_{i}\right)\mapsto\eta\otimes\mathbf{e}_{i}=\eta\mathbf{e}_{i}$
is injective since $\mathbf{e}_{i}\neq\mathbf{e}_{j}$ implies $\eta\mathbf{e}_{i}\neq\eta\mathbf{e}_{j}$,
and $\left\{ \eta\mathbf{e}_{i}\right\} $ is a basis for $V$, so
$\mu_{\otimes}\left(K,V\right)$ satisfies ($\otimes$). Similarly,
$\mu_{\otimes}\left(V,K\right)$ is a tensor product map, so $K\otimes V=V\otimes K=V$.
In particular, we have a tensor product map 
\[
\mu_{\otimes}\left(K,K\right):K\times K\rightarrow V,\qquad\left(\alpha,\beta\right)\mapsto\alpha\otimes\beta=\alpha\beta,
\]
so $K\otimes K=K$. The usual scalar product is thus a tensor product,
and the usual product of scalars is also a tensor product.

\emph{(3)}. Although we have constructed tensor product maps $\mu_{\otimes}\left(K,V\right)$,
$\mu_{\otimes}\left(V,K\right)$ and $\mu_{\otimes}\left(K,K\right)$
for any field $K$ and any vector space $V$ over $K$, we have not
yet shown that a tensor product map $\mu_{\otimes}\left(U,V\right)$
exists for any $U,V$. This follows immediately from the results in
Section 3, however. We define the \emph{distinguished tensor product
map} $\underline{\mu}\phantom{}_{\otimes}\left(U,V\right)$ by setting
\[
\mathbf{u}\otimes\mathbf{v}=\underline{\mu}\phantom{}_{\otimes}\left(U,V\right)\left(\mathbf{u},\mathbf{v}\right)=\mathbf{u}\mathbf{v},
\]
where $\mathbf{u}\mathbf{v}$ is the product of $\mathbf{u}\in U$
and $\mathbf{v}\in V$ in the free unital associative algebra $\mathcal{A}\left(U,V\right)$.
Vector spaces of the form $U\otimes V$ can also be defined in terms
of the distinguished tensor product map $\underline{\mu}\phantom{}_{\otimes}\left(U,V\right)$.
Specifically, $U\otimes V$ is the subspace of $\mathcal{A}\left(U,V\right)$
spanned by all products $\mathbf{u}\otimes\mathbf{v}$, or $\mathbf{u}\mathbf{v}$,
where $\mathbf{u}\in U,\mathbf{v}\in V$. We call $\mathbf{u}\otimes\mathbf{v}$
($U\otimes V$) \emph{the} tensor product of $\mathbf{u}$ and $\mathbf{v}$
($U$ and $V$). 

It is easy to verify that the definitions of the tensor product maps
$\mu_{\otimes}\left(K,V\right)$, $\mu_{\otimes}\left(V,K\right)$
and $\mu_{\otimes}\left(K,K\right)$ just given are consistent with
the definitions of the corresponding distinguished tensor product
maps $\underline{\mu}\phantom{}_{\otimes}\left(K,V\right)$, $\underline{\mu}\phantom{}_{\otimes}\left(V,K\right)$
and $\underline{\mu}\phantom{}_{\otimes}\left(K,K\right)$. With $U\otimes V$
defined by $\underline{\mu}\phantom{}_{\otimes}\left(U,V\right)$,
we recover the identities $K\otimes V=V\otimes K=V$.

\emph{(4)}. We are not yet done, because we want to show that the
distinguished tensor product map is associative, and the notion of
associative multiplication requires that we consider at least three
factors simultaneously. The key to understanding the situation is
that the multiplication operation in $\mathcal{A}\left(V_{1},\ldots,V_{n}\right)$
defines a \emph{set} of distinguished tensor product maps
\[
\left\{ \underline{\mu}\phantom{}_{\otimes}\left(X,Y\right)\mid X,Y\in\mathscr{V}\right\} ,
\]
where $\mathscr{V}$ is the smallest set of subspaces of $\mathcal{A}\left(V_{1},\ldots,V_{n}\right)$
such that
\begin{lyxlist}{00.00.0000}
\item [{$\qquad$(a)}] $K,V_{1},\ldots,V_{n}\in\mathscr{V}$;
\item [{$\qquad$(b)}] if $X,Y\in\mathscr{V}$ then $X\otimes Y\in\mathscr{V}.$ 
\end{lyxlist}
Since multiplication in $\mathcal{A}\left(V_{1},\ldots,V_{n}\right)$
is associative, we have
\begin{gather}
\underline{\mu}\phantom{}_{\otimes}\left(U\otimes V,W\right)\left(\mathbf{u}\otimes\mathbf{v},\mathbf{w}\right)=\underline{\mu}\phantom{}_{\otimes}\left(U,V\otimes W\right)\left(\mathbf{u},\mathbf{v}\otimes\mathbf{w}\right),\label{eq:gentensmap-1}\\
\forall\mathbf{u}\in U,\:\forall\mathbf{v}\in V,\:\forall\mathbf{w}\in W,\nonumber 
\end{gather}
so that we can write 
\[
\left(\mathbf{u}\otimes\mathbf{v}\right)\otimes\mathbf{w}=\mathbf{u}\otimes\left(\mathbf{v}\otimes\mathbf{w}\right).
\]
As $U\otimes V$ has a basis with elements of the form $\mathbf{u}\otimes\mathbf{v}$,
and $V\otimes W$ a basis with elements of the form $\mathbf{v}\otimes\mathbf{w}$,
$\left(U\otimes V\right)\otimes W$ is generated by tensors of the
form $\left(\mathbf{u}\otimes\mathbf{v}\right)\otimes\mathbf{w}$
and $U\otimes\left(V\otimes W\right)$ by tensors of the form $\mathbf{u}\otimes\left(\mathbf{v}\otimes\mathbf{w}\right)$,
so
\[
\left(U\otimes V\right)\otimes W=U\otimes\left(V\otimes W\right).
\]

\emph{(5)}. Although $\mu_{\otimes}\left(U,V\right)$ depends on $U$
and $V$, $\mu_{\otimes}\left(U,V\right)$ will be written as $\mu_{\otimes}$
below to simplify the notation when this is not likely to lead to
any misunderstanding. Furthermore, we shall use the same symbol $\mu{}_{\otimes}$
for an arbitrary tensor product map and the distinguished product
map $\underline{\mu}\phantom{}_{\otimes}$ defined above. In cases
where this distinction is important, the intended interpretation will
hopefully be clear from the context.

\subsection{Equivalence of linear and bilinear maps}

$ $\\
\emph{(1)}. Consider vector spaces $U$, $V$, $T$ and a tensor product
map $\mu_{\otimes}:U\times V\rightarrow U\otimes V$. For any linear
map $\lambda:U\otimes V\rightarrow T$, $\lambda\circ\mu_{\otimes}$
is a bilinear map $\mu:U\times V\rightarrow T$, since $\mu_{\otimes}$
is bilinear. Thus, $\mu_{\otimes}$ defines a mapping from linear
to bilinear maps
\[
\mu_{\otimes}\mapsto\left(\lambda\mapsto\mu=\lambda\circ\mu_{\otimes}\right).
\]

Let us show that $\mu_{\otimes}$ also defines a mapping in the opposite
direction
\[
\mu_{\otimes}\mapsto\left(\mu\mapsto\lambda\right),
\]
because the equation $\mu=\lambda\circ\mu_{\otimes}$ has a unique
solution $\lambda$ for given $\mu_{\otimes}$ and $\mu$. Choose
bases $\left\{ \mathbf{e}_{i}\right\} $ and $\left\{ \mathbf{f}_{j}\right\} $
for $U$ and $V$. Since $\left(\mathbf{e}_{i},\mathbf{f}_{j}\right)\mapsto\mathbf{e}_{i}\otimes\mathbf{f}_{j}$
is injective, there is a unique mapping $\ell:\left\{ \mathbf{e}_{i}\otimes\mathbf{f}_{j}\right\} \rightarrow T$
such that $\ell\left(\mathbf{e}_{i}\otimes\mathbf{f}_{j}\right)=\mu\left(\mathbf{e}_{i},\mathbf{f}_{j}\right)$
for all $\mathbf{e}_{i}$ and $\mathbf{f}_{j}$, and since $\left\{ \mathbf{e}_{i}\otimes\mathbf{f}_{j}\right\} $
is a basis for $U\otimes V$, $\ell$ can be extended to a unique
linear map $\lambda:U\otimes V\rightarrow T$ such that $\lambda\left(\mathbf{e}_{i}\otimes\mathbf{f}_{j}\right)=\ell\left(\mathbf{e}_{i}\otimes\mathbf{f}_{j}\right)=\mu\left(\mathbf{e}_{i},\mathbf{f}_{j}\right)$.
Thus, $\lambda\circ\mu_{\otimes}\!\left(\mathbf{e}_{i},\mathbf{f}_{j}\right)=\mu\!\left(\mathbf{e}_{i},\mathbf{f}_{j}\right)$
for all $\mathbf{e}_{i}$ and $\mathbf{f}_{j}$, and since $\lambda\circ\mu_{\otimes}$
and $\mu$ are bilinear, this implies that $\lambda\circ\mu_{\otimes}\left(\mathbf{u},\mathbf{v}\right)=\mu\left(\mathbf{u},\mathbf{v}\right)$
for all $\mathbf{u}\in U$ and $\mathbf{v}\in V$, because
\begin{gather*}
\lambda\!\circ\!\mu_{\otimes}\!\left(\left.\sum\right._{i}\!\alpha_{i}\mathbf{e}_{i},\left.\sum\right._{j}\!\beta_{j}\mathbf{f}_{j}\right)=\left.\sum\right._{i,j}\!\alpha_{i}\beta_{j}\,\lambda\!\circ\!\mu_{\otimes}\!\left(\mathbf{e}_{i},\mathbf{f}_{j}\right)\\
=\left.\sum\right._{i,j}\!\alpha_{i}\beta_{j}\,\mu\!\left(\mathbf{e}_{i},\mathbf{f}_{j}\right)=\mu\!\left(\left.\sum\right._{i,}\!\alpha_{i}\mathbf{e}_{i},\left.\sum\right._{j}\!\beta_{j}\mathbf{f}_{j}\right).
\end{gather*}

We conclude that $\mu_{\otimes}$ induces a one-to-one correspondence
$\lambda\leftrightarrow\mu_{\lambda}$ or $\mu\leftrightarrow\lambda_{\mu}$.

\emph{(2)}. As an example, consider a function 
\[
\mu_{\mathcal{C}}:U\times V\rightarrow V\otimes U,\qquad\left(\mathbf{u},\mathbf{v}\right)\mapsto\mathbf{v}\otimes\mathbf{u}.
\]
 For the first argument we have
\begin{gather*}
\mu_{\mathcal{C}}\!\left(\eta\mathbf{u},\mathbf{v}\right)=\mathbf{v}\otimes\eta\mathbf{u}=\eta\left(\mathbf{v}\otimes\mathbf{u}\right)=\eta\,\mu_{\mathcal{C}}\!\left(\mathbf{u},\mathbf{v}\right),\\
\mu_{\mathcal{C}}\!\left(\left(\mathbf{u}_{1}+\mathbf{u}_{2}\right),\mathbf{v}\right)=\mathbf{v}\otimes\left(\mathbf{u}_{1}+\mathbf{u}_{2}\right)=\mathbf{v}\otimes\mathbf{u}_{1}+\mathbf{v}\otimes\mathbf{u}_{2}=\mu_{\mathcal{C}}\!\left(\mathbf{u}_{1},\mathbf{v}\right)+\mu_{\mathcal{C}}\!\left(\mathbf{u}_{2},\mathbf{v}\right),
\end{gather*}
and similarly for the second argument, so $\mu_{\mathcal{C}}$ is
bilinear and gives a unique linear map 
\begin{equation}
\Lambda_{C}:U\otimes V\rightarrow V\otimes U.\label{eq:isomcom}
\end{equation}
By definition, $\Lambda_{C}\left(\mathbf{u}\otimes\mathbf{v}\right)=\mu_{\mathcal{C}}\!\left(\mathbf{u},\mathbf{v}\right)=\mathbf{v}\otimes\mathbf{u}$,
$ $so $\Lambda_{C}$ maps a basis $\left\{ \mathbf{e}_{i}\otimes\mathbf{f}_{j}\right\} $
for $U\otimes V$ bijectively to a basis $\left\{ \mathbf{f}_{j}\otimes\mathbf{e}_{i}\right\} $
for $V\otimes U$, so $\Lambda_{C}$ is a canonical isomorphism.

\emph{(3)}. $\Lambda_{C}$ is also given by the \emph{internally linear}
map
\[
\widetilde{\lambda}_{C}:\left.U\otimes V\right|_{\mathcal{S}}\rightarrow V\otimes U,\qquad\mathbf{u}\otimes\mathbf{v}\mapsto\mathbf{v}\otimes\mathbf{u},
\]
because $\mu_{\mathcal{C}}$ can be obtained from $\widetilde{\lambda}_{C}$
by setting $\mu_{\mathcal{C}}\left(\mathbf{u},\mathbf{v}\right)=\widetilde{\lambda}_{C}\left(\mathbf{u}\otimes\mathbf{v}\right)$. 

In the general case, $\widetilde{\lambda}_{\mu}$ is the restriction
of $\lambda_{\mu}$ to the set of simple tensors in $U\otimes V$,
and $\lambda_{\mu}$ can be recovered from $\widetilde{\lambda}_{\mu}$
through\emph{ }extension by linearity\emph{,} setting
\[
\lambda_{\mu}\left(\sum_{i,j}t^{ij}\left(\mathbf{u}_{i}\otimes\mathbf{v}_{j}\right)\right)=\sum_{i,j}t^{ij}\widetilde{\lambda}_{\mu}\left(\mathbf{u}_{i}\otimes\mathbf{v}_{j}\right).
\]

\begin{rem*}
$\widetilde{\lambda}_{\mu}$ is 'internally linear' in the sense that
if $\mathbf{u},\mathbf{v},\alpha\mathbf{u}+\beta\mathbf{v}\in\left.U\otimes V\right|_{\mathcal{S}}$
then $\widetilde{\lambda}_{\mu}\left(\alpha\mathbf{u}+\beta\mathbf{v}\right)=\alpha\widetilde{\lambda}_{\mu}\left(\mathbf{u}\right)+\beta\widetilde{\lambda}_{\mu}\left(\mathbf{v}\right)$.
Note, though, that $\widetilde{\lambda}_{\mu}$ is not a linear map,
since $\left.U\otimes V\right|_{\mathcal{S}}$ is not a closed subset
of $U\otimes V$ under addition.
\end{rem*}

\subsection{Uniqueness of the tensor product}

\emph{$ $}\\
If $\mu_{\otimes}:U\times V\rightarrow W$ is a tensor product map
and $i:W\rightarrow W'$ is an isomorphism then $\mu_{\otimes}'=i\circ\mu_{\otimes}$
is clearly a tensor product map $U\times V\rightarrow W'$, so any
vector space isomorphic to a tensor product $U\otimes V$ is itself
a tensor product $U\otimes'V$. 

Conversely, if $\mu_{\otimes}:U\times V\rightarrow W$ and $\mu_{\otimes}':U\times V\rightarrow W'$
are tensor product maps, $\left\{ \mathbf{e}_{i}\right\} $ is a basis
for $U$ and $\left\{ \mathbf{f}_{j}\right\} $ a basis for $V$,
then $\mu_{\otimes}\left(\mathbf{e}_{i},\mathbf{f}_{j}\right)\mapsto\mu_{\otimes}'\left(\mathbf{e}_{i},\mathbf{f}_{j}\right)$
maps $\left\{ \mathbf{e}_{i}\otimes\mathbf{f}_{j}\right\} $, a basis
for $W$, injectively onto $\left\{ \mathbf{e}_{i}\otimes'\mathbf{f}_{j}\right\} $,
a basis for $W'$, so this map gives an isomorphism $\iota:W\rightarrow W'$.
Hence, any two tensor products of $U$ and $V$ are isomorphic. Note
that $\mu_{\otimes}'=\iota\circ\mu_{\otimes}$, and the main result
in Subsection 4.2(1) implies that $\iota$ is unique since $\mu_{\otimes}'$
is bilinear. Thus, $\iota$ does not depend on the choice of the bases
$\left\{ \mathbf{e}_{i}\right\} $ and $\left\{ \mathbf{f}_{j}\right\} $
but only on $\mu_{\otimes}$ and $\mu_{\otimes}'$.
\begin{rem*}
In modern expositions of algebraic tensor theory, a preferred tensor
product map is usually not defined. Consistent with this approach,
a tensor product space is defined only up to an isomorphism. This
means that we cannot have $X\otimes Y=Z$, only $X\otimes Y\cong Z$.
Hence, we cannot have $K\otimes V=V\otimes K=V$ but only $K\otimes V\cong V\otimes K\cong V$,
and we cannot have $\left(U\otimes V\right)\otimes W=U\otimes\left(V\otimes W\right)$
but only $\left(U\otimes V\right)\otimes W\cong U\otimes\left(V\otimes W\right)$.
Also note that if $i:X\rightarrow Y$ is an isomorphism and $\alpha:X\rightarrow X$
an automorphism then $i\circ\alpha:X\rightarrow Y$ is an isomorphism
as well, so we do not, for example, have a unique isomorphism $\left(U\otimes V\right)\otimes W\rightarrow U\otimes\left(V\otimes W\right)$. 

One way to make the tensor product of $U$ and $V$ 'more unique'
is to define it not as a particular vector space $U\otimes V$ but
as the \emph{pair} $\left(U\otimes V,\mu_{\otimes}\right)$, where
$\mu_{\otimes}$ is a tensor product map $U\times V\rightarrow U\otimes V$.
Then there is a unique correspondence between $\left(U\otimes V,\mu_{\otimes}\right)$
and $\left(U\otimes'V,\mu_{\otimes}'\right)$ in the sense that there
is a unique isomorphism \linebreak{}
$\iota:U\otimes V\rightarrow U\otimes'V$ such that $\mu_{\otimes}'=\iota\circ\mu_{\otimes}$.
If we do not\emph{ }fix the tensor map $\mu_{\otimes}$ by setting
$\mu_{\otimes}=\underline{\mu}\phantom{}_{\otimes}$, the vector space
$U\otimes V$ is still defined only up to isomorphism, however, and
$\iota$ is not the only available isomorphism $U\otimes V\rightarrow U\otimes'V$.
(Singling out $\iota$ merely allows us to regard a heap of isomorphisms
$\left\{ i,i^{-1}\mid i:U\otimes V\rightarrow U\otimes'V\right\} $
as a group of automorphisms $\left\{ \iota^{-1}\circ i\mid i:U\otimes V\rightarrow U\otimes'V\right\} $.)

It is easy to show (by induction) that if $U\otimes\left(V\otimes W\right)=\left(U\otimes V\right)\otimes W$
then we can write $V_{1}\otimes\ldots\otimes V_{n}$ without ambiguity;
any two tensor products of $V_{1},\ldots,V_{n}$ are equal, regardless
of the arrangement of parentheses. This is a highly desirable property,
but isomorphisms are not quite as well-behaved in this respect as
equalities. The point is that two vector spaces can be equal in one
way only, but they can be isomorphic in more than one way. For example,
an isomorphism $I:A\otimes\left(B\otimes\left(C\otimes D\right)\right)\rightarrow\left(\left(A\otimes B\right)\otimes C\right)\otimes D$
can be constructed from a given isomorphism $i:U\otimes\left(V\otimes W\right)\rightarrow\left(U\otimes V\right)\otimes W$
through either one of two chains of isomorphisms, 
\begin{gather*}
A\otimes\left(B\otimes\left(C\otimes D\right)\right)\rightarrow\left(A\otimes B\right)\otimes\left(C\otimes D\right)\rightarrow\left(\left(A\otimes B\right)\otimes C\right)\otimes D\quad\mathrm{or}\\
A\!\otimes\!\left(B\!\otimes\!\left(C\!\otimes\! D\right)\right)\rightarrow A\!\otimes\!\left(\left(B\!\otimes\! C\right)\!\otimes\! D\right)\rightarrow\left(A\!\otimes\!\left(B\!\otimes\! C\right)\right)\!\otimes\! D\rightarrow\left(\left(A\!\otimes\! B\right)\!\otimes\! C\right)\!\otimes\! D,
\end{gather*}
and $I$ depends not only on $i$ but also on the chain used to derive
$I$ from $i$. To ensure the uniqueness of isomorphisms such as $I$,
one can require that certain \emph{coherence axioms} hold. Then isomorphisms
behave as equalities and can be interpreted as equalities. Specifically,
one may introduce coherence axioms that allow the isomorphisms $K\otimes V\cong V$,
$V\otimes K\cong V$ and $U\otimes\left(V\otimes W\right)\cong\left(U\otimes V\right)\otimes W$
(but not $U\otimes V\cong V\otimes U$) to be interpreted as equalities.

In this article, \emph{the} tensor product of $U$ and $V$ is defined
to be a particular vector space $U\otimes V$, given by a unique distinguished
tensor product map as described earlier. Recall that $\underline{\mu}\phantom{}_{\otimes}\left(\mathbf{u},\mathbf{v}\right)=\mathbf{u}\mathbf{v}$,
so the uniqueness of the tensor product follows from the uniqueness
of the corresponding free unital associative algebra on $U$ and $V$
(or, in the general case, $V_{1},\ldots,V_{n}$). 

In this approach, we have $\left(U\otimes V\right)\otimes W=U\otimes\left(V\otimes W\right)$
simply because the algebra $\mathcal{A}$ from which we get the tensor
product is associative. Similarly,\linebreak{}
 $K\otimes V=V\otimes K=V$ because $\mathcal{A}$ is unital. On the
other hand, we do not have $U\otimes V=V\otimes U$ although $U\otimes V\cong V\otimes U$,
because $\mathcal{A}$ is not commutative. In other words, tensor
products of vectors and vector spaces are naturally associative but
not commutative, as assumed in applications; there are also 'native'
unit elements for tensor products.

To finally put Sections 3 and 4 into perspective, the real news is
that the tensor product is defined in terms of the tensor algebra
rather than vice versa, for reasons partly given in this remark.
\end{rem*}

\section{Tensor products of $n$ vectors or vector spaces}

\subsection{Tensor products, tuple products and Cartesian products}

$ $\\
\emph{(1)}. Because of the associative multiplication property (\ref{eq:gentensmap-1}),
we can write tensor products of the form $\mathbf{v}_{1}\otimes\mathbf{v}{}_{2}\otimes\ldots\otimes\mathbf{v}{}_{n}$
and $V_{1}\otimes V_{2}\otimes\ldots\otimes V_{n}$ without ambiguity.
We extend the definition of a tensor from an element of $U\otimes V$
to an element of $V_{1}\otimes V_{2}\otimes\ldots\otimes V_{n}$.
Using ($\otimes)$ and (\ref{eq:gentensmap-1}), it can be shown that
($\otimes^{n}$) if $\left\{ \mathbf{e}_{ij}\right\} $ is a basis
for $V_{i}$ then $\left\{ \mathbf{e}_{1j_{1}}\otimes\ldots\otimes\mathbf{e}_{nj_{n}}\mid\mathbf{e}_{ij_{k}}\in\left\{ \mathbf{e}_{ij}\right\} \right\} $
is a basis for $V_{1}\otimes\ldots\otimes V_{n}$. 

\emph{(2)}. We can also denote the tensor product of $\mathbf{v}_{1},\ldots,\mathbf{v}_{n}$
for $n\geq0$ by 
\[
\left.\bigotimes\right._{i=1}^{n}\mathbf{v}_{i},\qquad\mathrm{setting}\;\,\left.\bigotimes\right._{i=1}^{0}\mathbf{v}_{i}=\xi,
\]
where $\xi\in K$ is a scalar variable rather that a fixed scalar.
Similarly, we may denote the tensor product of $V_{1},\ldots,V_{n}$
for $n\geq0$ by 
\[
\left.\bigotimes\right._{i=1}^{n}V_{i},\qquad\mathrm{where}\;\,\left.\bigotimes\right._{i=1}^{0}V_{i}=K.
\]

We can generalize the notion of simple tensors in an obvious way;
these are vectors in $\bigotimes_{i=1}^{n}V_{i}$ of the form $\bigotimes_{i=1}^{n}\mathbf{v}_{i}$,
where $\mathbf{v}_{i}\in V_{i}$, for $n\geq0$. The set of all such
simple tensors will be denoted by $\left.V_{1}\otimes\ldots\otimes V_{n}\right|_{\mathcal{S}}$
or $\left.\bigotimes_{i=1}^{n}V_{i}\right|_{\mathcal{S}}$. Recall
that $V_{1}\otimes\ldots\otimes V_{n}$ has a basis of simple tensors;
this can be extended to the case $n=0$.

Tensor products of $n$ copies of the same vector space $V$ are of
particular interest. We define the \emph{tensor power} $V^{\otimes n}$
of $V$ for $n\geq0$ by
\[
V^{\otimes n}=\left.\bigotimes\right._{i=1}^{n}V_{i}\qquad(V_{i}=V),\qquad\mathrm{where}\quad V^{\otimes0}=K.
\]

\emph{(3)}.~Let $\mathbf{v}_{i}\in V_{i}$, where $V_{i}$ is a vector
space over $K$. The \emph{tuple product} of $\mathbf{v}_{1},\ldots,\mathbf{v}_{n}$
for $n\geq0$ is defined as 
\[
\left.\prod\right._{i=1}^{n}\mathbf{v}{}_{i}=\left(\mathbf{v}_{1},\ldots,\mathbf{v}_{n}\right)\quad(n\geq1),\qquad\left.\prod\right._{i=1}^{0}\mathbf{v}_{i}=\xi,
\]
where $\xi\in K$ is a scalar variable.

The \emph{(associative) Cartesian product} of $V_{1},\ldots,V_{n}$
for $n\geq0$ is the tuple set
\[
\left.\prod\right._{i=1}^{n}V{}_{i}=V_{1}\times\ldots\times V_{n}\quad(n\geq1),\qquad\left.\prod\right._{i=1}^{0}V{}_{i}=K,
\]
where $V_{1}\times\ldots\times V_{n}$ generalizes the definition
of the (associative) Cartesian product of two tuple sets in Subsection
2.2.

We can also define a \emph{Cartesian power} of a vector space $V$
for $n\geq0$ by setting

\begin{gather*}
V^{n}=\left.\prod\right._{i=1}^{n}V_{i}\quad(V_{i}=V),\qquad\mathrm{where}\quad V^{0}=K.
\end{gather*}

\subsection{Equivalence of linear and separately linear maps}

\emph{$ $}\\
\emph{(1)}. A \emph{separately linear map }(or $n$-\emph{linear}
\emph{map)} is a function 
\[
\boldsymbol{\phi}:\left.\prod\right._{i=1}^{n}V_{i}\rightarrow T\quad\left(n\geq0\right),
\]
which is linear in each argument separately for fixed values of all
other arguments (if any). We let 
\[
\mathscr{L}\left[V_{1}\times\ldots\times V_{n},T\right]\quad\mathrm{or}\quad\mathscr{L}\left[\left.\prod\right._{i=1}^{n}V_{i},T\right]\quad\left(n\geq0\right)
\]
denote the set of all $n$-linear maps of the form displayed. For
example, $\mathscr{L}\left[K,T\right]$ is a set of 0-linear maps,
and $\mathscr{L}\left[V,T\right]$ and $\mathscr{L}\left[\bigotimes_{i=1}^{n}V_{i},T\right]$
are sets of 1-linear maps. Note that the difference between linear
and separately linear maps is sometimes a difference between points
of view, since the separately linear maps $K\rightarrow T$ and $V\rightarrow T$
are linear maps as well.

\emph{(2)}. As usual, we denote $\mathscr{L}\left[V,K\right]$ by
$V^{*}$; for convenience, we write $\left(V^{*}\right)^{n}$ as $V^{*n}$
and $\left(V^{*}\right)^{\otimes n}$ as $V^{*\otimes n}$.

\emph{(3)}. Setting $\left(\eta\boldsymbol{\phi}\right)\!\left(\mathbf{s}\right)=\eta\left(\boldsymbol{\phi}\!\left(\mathbf{s}\right)\right)$
and $\left(\boldsymbol{\phi}_{1}+\boldsymbol{\phi}_{2}\right)\!\left(\mathbf{s}\right)=\boldsymbol{\phi}_{1}\!\left(\mathbf{s}\right)+\boldsymbol{\phi}_{2}\!\left(\mathbf{s}\right)$
clearly makes $\mathscr{L}\left[\bigotimes_{i=1}^{n}V_{i},T\right]$
into a vector space. $\mathscr{L}\left[\prod_{i=1}^{n}V_{i},T\right]$
can be made into a vector space in the same way, setting $\left(\eta\boldsymbol{\phi}\right)\!\left(\mathbf{v}_{1},\ldots,\mathbf{v}_{n}\right)=\eta\left(\boldsymbol{\phi\!}\left(\mathbf{v}_{1},\ldots,\mathbf{v}_{n}\right)\right)$
and $\left(\boldsymbol{\phi}_{1}+\boldsymbol{\phi}_{2}\right)\!\left(\mathbf{v}_{1},\ldots,\mathbf{v}_{n}\right)=\boldsymbol{\phi}_{1}\!\left(\mathbf{v}_{1},\ldots,\mathbf{v}_{n}\right)+\boldsymbol{\phi}_{2}\!\left(\mathbf{v}_{1},\ldots,\mathbf{v}_{n}\right)$.

\emph{(4)}. The equivalence between linear and bilinear maps introduced
in Subsection 4.2 can be generalized from bilinear maps to all separately
linear maps. It can be proved in essentially the same way that for
any separately linear map 
\[
\mu:\left.\prod\right._{i=1}^{n}V_{i}\rightarrow T,
\]
there is an equivalent linear map
\[
\lambda_{\mu}:\left.\bigotimes\right._{i=1}^{n}V_{i}\rightarrow T
\]
such that $\lambda_{\mu}\left(\bigotimes_{i=1}^{n}\mathbf{v}{}_{i}\right)=\mu\left(\prod_{i=1}^{n}\mathbf{v}{}_{i}\right)$,
and conversely for linear maps.

$\lambda_{\mu}$ can also be defined through linear extension of an
internally linear map 
\[
\widetilde{\lambda}_{\mu}:\left.\left.\bigotimes\right._{i=1}^{n}V_{i}\right|_{\mathcal{S}}\rightarrow T
\]
 such that $\widetilde{\lambda}_{\mu}\left(\bigotimes_{i=1}^{n}\mathbf{v}{}_{i}\right)=\mu\left(\prod_{i=1}^{n}\mathbf{v}{}_{i}\right)=\lambda_{\mu}\left(\bigotimes_{i=1}^{n}\mathbf{v}{}_{i}\right)$.

Thus we have bijections of the form 
\begin{equation}
\Lambda_{\mathcal{M}}:\mathscr{L}\left[\left.\prod\right._{i=1}^{n}V_{i},T\right]\rightarrow\mathscr{L}\left[\left.\bigotimes\right._{i=1}^{n}V_{i},T\right],\qquad\mu\mapsto\lambda_{\mu}\label{eq:seplinTenslim}
\end{equation}
for all $n\geq0$. If $\mathscr{L}\left[\prod_{i=1}^{n}V_{i},T\right]$
and $\mathscr{L}\left[\bigotimes_{i=1}^{n}V_{i},T\right]$ are equipped
with vector space structures as just described, the mapping $\Lambda_{\mathcal{M}}$
is clearly linear, and thus a canonical isomorphism. This implies,
in particular, that $\mathscr{L}\left[V^{n},K\right]\cong\mathscr{L}\left[V^{\otimes n},K\right]$,
$\mathscr{L}\left[V^{n},V\right]\cong\mathscr{L}\left[V^{\otimes n},V\right]$
and $\mathscr{L}\left[V^{n},V^{\otimes p}\right]\cong\mathscr{L}\left[V^{\otimes n},V^{\otimes p}\right]$
for any $p>1$.

\specialsection*{\textbf{\textsl{B. Interpretations of classical tensors}}}

\section{Classical tensors as bilateral tensors}

\subsection{Classical tensors and related notions}

\emph{$ $$ $}\\
\emph{(1)}. We previously defined a\emph{ }tensor simply as an\emph{
}element of a tensor product of vector spaces. (As $K=K\otimes K$
and $V=K\otimes V$, this includes scalars and vectors in any vector
space as well.) What is needed here, however, is a definition of \emph{a
(potentially) double-indexed tensor} -- a classical tensor.

In one common interpretation of classical tensors \cite{key-7}, they
are regarded as elements of tensor product spaces of the form
\begin{equation}
V^{\otimes n}\otimes V^{*\otimes m},\qquad(m,n\geq0).\label{eq:eqtensprodusual}
\end{equation}
Any $\boldsymbol{\mathsf{t}}\in V^{\otimes n}\otimes V^{*\otimes m}$
is said to be a tensor of \emph{valence} $\binom{n}{m}$; equivalently,
$\boldsymbol{\mathsf{\mathsf{t}}}$ is said to be \emph{covariant}
of degree $m$ and \emph{contravariant} of degree $n$.
\begin{rem*}
Classical tensors are sometimes defined as\emph{ multilinear forms}
\[
V^{*n}\times V^{m}\rightarrow K.
\]
We know that 
\[
\mathscr{L}\left[V^{*n}\times V^{m},K\right]\cong\mathscr{L}\left[V^{*\otimes n}\otimes V^{\otimes m},K\right],
\]
and it can also be shown that for finite-dimensional vector spaces
\[
\mathscr{L}\left[V^{*\otimes n}\otimes V^{\otimes m},K\right]\cong V^{\otimes n}\otimes V^{*\otimes m},
\]
so the definition of classical tensors as multilinear forms is closely
related to the definition in terms of tensor products of the form
(\ref{eq:eqtensprodusual}).

Classical tensors have also been defined as elements of tensor products
of $n$ copies of $V$ and $m$ copies of $V^{*}$ in any order \cite{key-1,key-3,key-8,key-10}.
Such tensors have been called 'affine tensors' \cite{key-3,key-10},
'homogeneous tensors' \cite{key-8} or simply 'tensors'. 

It should perhaps also be noted that a tensor \emph{field}, i.e.,
an assignment of a classical tensor to each point in some space, is
often informally referred to as a tensor. Although tensors were historically
more or less synonymous with tensor fields, we are only concerned
with the algebraic notion of classical tensors here.
\end{rem*}
\emph{(2)}. As tensor multiplication of vector spaces is associative
and also commutative up to a canonical isomorphism, one may as well
interpret a classical tensor of valence $\binom{n}{m}$ as an element
of 
\begin{equation}
V^{*\otimes m}\otimes V^{\otimes n}\qquad(m,n\geq0),\label{eq:tensprod}
\end{equation}
and this convention turns out to be the most convenient one in the
present context. Elements of tensor products of vector spaces of the
form (\ref{eq:eqtensprodusual}) or (\ref{eq:tensprod}) may be referred
to as\emph{ bilateral tensors}; the form (\ref{eq:tensprod}) will
be used in the definitions below.\emph{ }

\emph{(3)}. Addition and scalar multiplication of bilateral tensors
are well-defined operations simply because a bilateral tensor is a
vector in a tensor product space.

\emph{(4)}. Recall that $V^{*\otimes0}=V^{\otimes0}=K$ so that $V^{*\otimes0}\otimes V^{\otimes0}=K$,
$V^{*\otimes0}\otimes V^{\otimes1}=V$ and $V^{*\otimes1}\otimes V^{\otimes0}=V^{*}$.
It is clear that every $\sigma\in K$ is an element $\phi\otimes\xi=\phi\xi$
of $V^{*\otimes0}\otimes V^{\otimes0}$, every vector $\mathbf{v}\in V$
is an element $\phi\otimes\mathbf{u}=\phi\mathbf{u}$ of $ $$V^{*\otimes0}\otimes V^{\otimes1}$,
and every linear form $\boldsymbol{f}\in V^{*}$ is an element $\boldsymbol{g}\otimes\xi=\xi\boldsymbol{g}$
of $V^{*\otimes1}\otimes V^{\otimes0}$, although these representations
of $\sigma$, $\mathbf{v}$ and $\boldsymbol{f}$ as bilateral tensors
are of course not unique.

\subsection{Notation for bilateral tensors}

$ $\\
There are two ways to refer to bilateral tensors. In \emph{index-free
notation}, a tensor in $V^{*\otimes m}\otimes V^{\otimes n}$ is referred
to by a symbol such as $\boldsymbol{\mathsf{t}}$, in bold sans-serif
typeface. In \emph{double-index notation}, $\boldsymbol{\mathsf{t}}$
is written in the form 
\[
\mathsf{t}_{\mathrm{a}_{1}\ldots\mathrm{a}_{m}}^{\mathrm{b}_{1}\ldots\mathrm{b}_{n}}
\]
with the stem letter in sans-serif typeface and $m+n$ distinct indices.
Roman font is used instead of italics for indices to emphasize that
they are not variables with values such as $1,2,3,\ldots$, but labels
that identify different copies of $V^{*}$ or $V$ in a tensor product
of the form (\ref{eq:tensprod}). Specifically, a subscript points
to a copy of $V^{*}$, while a superscript points to a copy of $V$.
For a simple bilateral tensor $\widetilde{\boldsymbol{\mathsf{t}}}$,
\[
\mathsf{\widetilde{t}}_{\mathrm{a}_{1}\ldots\mathrm{a}_{m}}^{\mathrm{b}_{1}\ldots\mathrm{b}_{n}}=\left(\left.\bigotimes\right._{i=1}^{m}\boldsymbol{f}_{\mathrm{a}_{i}}\right)\otimes\left(\left.\bigotimes\right._{i=1}^{n}\mathbf{v}{}_{\mathrm{b}_{i}}\right)\qquad(\boldsymbol{f}_{\mathrm{a}_{i}}\in V^{*},\mathbf{v}{}_{\mathrm{b}_{i}}\in V).
\]
For a general bilateral tensor we have an expansion of the form
\[
\mathsf{t}_{\mathrm{a}_{1}\ldots\mathrm{a}_{m}}^{\mathrm{b}_{1}\ldots\mathrm{b}_{n}}=\left.\sum\right._{k}t^{k}\left(\left.\bigotimes\right._{i=1}^{m}\boldsymbol{f}_{\mathrm{a}_{i}k}\right)\otimes\left(\left.\bigotimes\right._{i=1}^{n}\mathbf{v}{}_{\mathrm{b}_{i}k}\right),
\]
since $V^{*\otimes m}\otimes V^{\otimes n}$ has a basis of simple
bilateral tensors of the form
\[
\left(\left.\bigotimes\right._{i=1}^{m}\boldsymbol{f}_{i}\right)\otimes\left(\left.\bigotimes\right._{i=1}^{n}\mathbf{v}{}_{i}\right)\qquad(\boldsymbol{f}_{i}\in V^{*},\mathbf{v}{}_{i}\in V).
\]

We have thus already introduced an abstract index notation for bilateral
tensors, where non-numerical indices embellishing double-indexed tensor
symbols characterize the tensor independently of its representation
by scalars. (A formally similar abstract index notation corresponding
to the interpretation of classical tensors as tensor maps will be
introduced later.)

It is sometimes convenient to use an extended notation where $\mathsf{\widetilde{t}}_{\mathrm{a}_{1}\ldots\mathrm{a}_{m}}^{1}$
denotes $\left(\left.\bigotimes\right._{i=1}^{m}\boldsymbol{f}_{\mathrm{a}_{i}}\right)\otimes1$
while $\mathsf{\widetilde{t}}_{1}^{\mathrm{b}_{1}\ldots\mathrm{b}_{n}}$
denotes $1\otimes\left(\left.\bigotimes\right._{i=1}^{n}\mathbf{v}{}_{\mathrm{b}_{i}}\right)$,
and so forth.

\subsection{Permutations of indices}

$ $\\
\emph{(1)}. Tensors can obviously be added, subtracted and compared
only if they belong to the same vector space. This fact leads to certain
\emph{consistency requirements for indices}. For example, $\mathsf{f}_{\mathrm{a}}$
and $\mathsf{v}^{\mathrm{a}}$ clearly belong to different kinds of
vector spaces, so expressions such as $\mathsf{f}_{\mathrm{a}}+\mathsf{v}^{\mathrm{a}}$
or $\mathsf{f}_{\mathrm{a}}=\mathsf{v}^{\mathrm{a}}$ are not meaningful.
In expressions where tensors are combined or compared as vectors,
different tensors are required to have the same subscripts and the
same superscripts. It is not required, however, that subscripts and
superscripts appear in the same order in all tensors. For example,
expressions like $\mathsf{t}^{\mathrm{ab}}=\mathsf{t}^{\mathrm{ba}}$
and $\mathsf{t}_{\mathrm{ab}}-\mathsf{t}_{\mathrm{ba}}$ are meaningful. 

Let us look more closely at this. Recall that $\Lambda_{\mathcal{C}}:U\otimes V\rightarrow V\otimes U,$
given by $\mathbf{u}\otimes\mathbf{v}\mapsto\mathbf{v}\otimes\mathbf{u}$
is a canonical isomorphism, so if $U=V$ then $\Lambda_{\mathcal{C}}:V^{\otimes2}\rightarrow V^{\otimes2}$
is a canonical automorphism. Similarly, 
\[
\mathsf{t}^{\mathrm{ab}}=\sum_{k}t^{k}\,\phi_{k}\otimes\mathbf{v}_{\mathrm{a}k}\otimes\mathbf{v}_{\mathrm{b}k}\;\mapsto\;\mathsf{t}^{\mathrm{ba}}=\sum_{k}t^{k}\,\phi_{k}\otimes\mathbf{v}_{\mathrm{b}k}\otimes\mathbf{v}_{\mathrm{a}k}
\]
gives an automorphism $K\otimes V^{\otimes2}\rightarrow K\otimes V^{\otimes2}$,
so $\mathsf{t}^{\mathrm{ab}}$ and $\mathsf{t}^{\mathrm{ba}}$ belong
to the same vector space. In terms of the standard basis $\left\{ 1\right\} $
for $K$ and a basis $\left\{ \mathbf{e}_{i}\right\} $ for $V$,
we have 
\[
\mathsf{t}^{\mathrm{ab}}=\sum_{i,j}t^{ij}\mathbf{e}{}_{i}\otimes\mathbf{e}_{j}\;\mapsto\;\mathsf{t}^{\mathrm{ba}}=\sum_{i,j}t^{ij}\mathbf{e}_{j}\otimes\mathbf{e}{}_{i}=\sum_{i,j}t^{ji}\mathbf{e}{}_{i}\otimes\mathbf{e}_{j}.
\]

In the general case, let $\mathrm{a}_{i}\mapsto\mathrm{a}_{i}'$ and
$\mathrm{b}_{i}\mapsto\mathrm{b}_{i}'$ be permutations of subscripts
and superscripts, respectively. Then there is a canonical automorphism
\begin{gather*}
V^{*\otimes m}\otimes V^{\otimes n}\rightarrow V^{*\otimes m}\otimes V^{\otimes n},\qquad\mathsf{t}_{\mathrm{a}_{1}\ldots\mathrm{a}_{m}}^{\mathrm{b}_{1}\ldots\mathrm{b}_{n}}\mapsto\mathsf{t}_{\mathrm{a}_{1}'\ldots\mathrm{a}_{m}'}^{\mathrm{b}_{1}'\ldots\mathrm{b}_{n}'}
\end{gather*}
given by
\[
\boldsymbol{f}_{\mathrm{a}_{1}}\otimes\ldots\otimes\boldsymbol{f}{}_{\mathrm{a}_{m}}\otimes\mathbf{v}_{\mathrm{b}_{1}}\otimes\ldots\otimes\mathbf{v}_{\mathrm{b}_{n}}\mapsto\boldsymbol{f}_{\mathrm{a}_{1}'}\otimes\ldots\otimes\boldsymbol{f}{}_{\mathrm{a}_{m}'}\otimes\mathbf{v}_{\mathrm{b}_{1}'}\otimes\ldots\otimes\mathbf{v}_{\mathrm{b}_{n}'}\,,
\]
so $t_{\mathrm{a}_{1}\ldots\mathrm{a}_{m}}^{\mathrm{b}_{1}\ldots\mathrm{b}_{n}}$
and $t_{\mathrm{a}_{1}'\ldots\mathrm{a}_{m}'}^{\mathrm{b}_{1}'\ldots\mathrm{b}_{n}'}$
belong to the same vector space, and it holds in general that bilateral
tensors can be combined and compared as long as they have the same
subscripts and superscripts, regardless of the order in which these
indices are arranged.

\emph{(2)}. There is a subtle problem with the notation used here,
because the indices in symbols such as $\mathsf{t}_{\mathrm{\mathrm{ab}}}$
have no intrinsic meaning; the meaning depends on the context. Hence,
we can write $\mathsf{t}_{\mathrm{ab}}$ as $\mathsf{t}_{\mathrm{xy}}$
or even $\mathsf{t}_{\mathrm{ba}}$. Similarly, we can write $\mathsf{t}_{\mathrm{ab}}-\mathsf{t}_{\mathrm{ba}}$
as $\mathsf{t}_{\mathrm{ba}}-\mathsf{t}_{\mathrm{ab}}$, but then
we run into difficulties, since we can also interpret $\mathsf{t}_{\mathrm{ba}}-\mathsf{t}_{\mathrm{ab}}$
as $-\left(\mathsf{t}_{\mathrm{ab}}-\mathsf{t}_{\mathrm{ba}}\right)$. 

This points to the fact that there is a difference between permutations
of indices due to \emph{formal substitutions} of indices and permutations
of indices due to \emph{automorphisms} $V^{*\otimes m}\otimes V^{\otimes n}\rightarrow V^{*\otimes m}\otimes V^{\otimes n}$.
If necessary, we can eliminate ambiguity due to this overloading phenomenon
by writing, for example, $\mathsf{t}_{\mathrm{b}^{\!2}\mathrm{a}^{\!1}}$
instead of $\mathsf{t}_{\mathrm{ba}}$, creating an implicit context
$\mathsf{t}_{\mathrm{ab}}\mapsto\mathsf{t}_{\mathrm{b'a'}}=\mathsf{t}_{\mathrm{b}^{\!2}\mathrm{a}^{\!1}}$,
where we have written $\mathsf{t}_{\mathrm{ba}}$ as $\mathsf{t}_{\mathrm{b}^{\!2}\mathrm{a}^{\!1}}$
because $\mathsf{t}_{\mathrm{ba}}$ is related to a \emph{reference
tensor} $\mathsf{t}_{\mathrm{ab}}$. It is clear how to extend this
convention to all tensors of valence $\binom{n}{m}$, where $m\geq2$
and/or $n\geq2$.

\subsection{Multiplication of bilateral tensors}

$ $\emph{}\\
\emph{(1)}. Below, $m,n,p,q,\varrho$ will denote non-negative integers
such that $\varrho\leq m,q$. Consider a function
\begin{gather}
\widetilde{\beta}_{\bullet}^{\varrho}:\left.\left(V^{*\otimes m}\otimes V^{\otimes n}\right)\right|_{\mathcal{S}}\times\left.\left(V^{*\otimes p}\otimes V^{\otimes q}\right)\right|_{S}\;\rightarrow\; V^{*\otimes\left(m+p-\varrho\right)}\otimes V^{\left(\otimes n+q-\varrho\right)},\nonumber \\
\left(\bigotimes_{i=1}^{m}\boldsymbol{f}^{i}\otimes\bigotimes_{i=1}^{n}\mathbf{v}_{i}\,,\,\bigotimes_{i=1}^{p}\boldsymbol{\bar{f}}\vphantom{}^{i}\otimes\bigotimes_{i=1}^{q}\bar{\mathbf{v}}_{i}\right)\mapsto\label{eq:bilmult}\\
\prod_{i=1}^{\varrho}\boldsymbol{f}^{i}\!\left(\bar{\mathbf{v}}_{i}\right)\bigotimes_{i=\varrho+1}^{m}\boldsymbol{f}^{i}\otimes\bigotimes_{i=1}^{p}\boldsymbol{\bar{f}}\vphantom{}^{i}\otimes\bigotimes_{i=1}^{n}\mathbf{v}_{i}\otimes\bigotimes_{i=\varrho+1}^{q}\bar{\mathbf{v}}_{i}.\nonumber 
\end{gather}
To take care of special cases, set $\prod_{i=1}^{0}\boldsymbol{f}^{i}\!\left(\bar{\mathbf{v}}_{i}\right)\!=\!1$,
$\bigotimes_{i=\varrho+1}^{\varrho>0}\boldsymbol{f}^{i}\!=\!1$ and
$\bigotimes_{i=\varrho+1}^{\varrho>0}\bar{\mathbf{v}}_{i}\!=\!1$.
(Recall that $\bigotimes_{i=1}^{0}\boldsymbol{f}^{i}=\phi$, $\bigotimes_{i=1}^{0}\boldsymbol{\bar{f}}\vphantom{}^{i}=\bar{\phi}$,
$\bigotimes_{i=1}^{0}\mathbf{v}_{i}=\xi$ and $\bigotimes_{i=1}^{0}\bar{\mathbf{v}}_{i}=\bar{\xi}$,
where $\phi,\bar{\phi},\xi,\bar{\xi}\in K$ are scalar variables.) 
\begin{example*}
Set $\widetilde{\boldsymbol{\mathsf{s}}}=\boldsymbol{f}\otimes\xi=\xi\boldsymbol{f}=\boldsymbol{g}$
and $\widetilde{\boldsymbol{\mathsf{t}}}=\bar{\phi}\otimes\bar{\mathbf{v}}=\bar{\phi}\bar{\mathbf{v}}=\bar{\mathbf{w}}$.
Then we have 
\[
\widetilde{\beta}_{\bullet}^{0}\left(\widetilde{\boldsymbol{\mathsf{s}}},\widetilde{\boldsymbol{\mathsf{t}}}\right)=1\boldsymbol{f}\otimes\bar{\phi}\otimes\xi\otimes\bar{\mathbf{v}}=\boldsymbol{g}\otimes\bar{\mathbf{w}},\qquad\widetilde{\beta}_{\bullet}^{1}\left(\widetilde{\boldsymbol{\mathsf{s}}},\widetilde{\boldsymbol{\mathsf{t}}}\right)=\boldsymbol{f}\!\left(\bar{\mathbf{v}}\right)1\otimes\bar{\phi}\otimes\xi\otimes1=\boldsymbol{g}\!\left(\bar{\mathbf{w}}\right).\quad\square
\]

\end{example*}
Inspection of (\ref{eq:bilmult}) reveals that the function
\[
\mu^{\varrho}:V^{*m}\times V^{n}\times V^{*p}\times V^{q}\;\rightarrow\; V^{*\otimes\left(m+p-\varrho\right)}\otimes V^{\left(\otimes n+q-\varrho\right)}
\]
defined by 
\begin{gather*}
\mu^{\varrho}\left(\prod_{i=1}^{m}\boldsymbol{f}^{i},\prod_{i=1}^{n}\mathbf{v}_{i}\,,\,\prod_{i=1}^{P}\bar{\boldsymbol{f}}\vphantom{}^{i},\prod_{I=1}^{q}\bar{\mathbf{v}}_{i}\right)=\widetilde{\beta}_{\bullet}^{\varrho}\left(\bigotimes_{i=1}^{m}\boldsymbol{f}^{i}\otimes\bigotimes_{i=1}^{n}\mathbf{v}_{i}\,,\,\bigotimes_{i=1}^{p}\boldsymbol{\bar{f}}\vphantom{}^{i}\otimes\bigotimes_{i=1}^{q}\bar{\mathbf{v}}_{i}\right),
\end{gather*}
is separately linear, so corresponding to $\widetilde{\beta}_{\bullet}^{\varrho}$
there is a linear map 
\[
\lambda^{\varrho}:V^{*m}\otimes V^{n}\otimes V^{*p}\otimes V^{q}\;\rightarrow\; V^{*\otimes\left(m+p-\varrho\right)}\otimes V^{\left(\otimes n+q-\varrho\right)},
\]
and corresponding to $\lambda^{\varrho}$ there is a bilinear map
\begin{gather*}
\beta_{\bullet}^{\varrho}:\left(V^{*\otimes m}\otimes V^{\otimes n}\right)\times\left(V^{*\otimes p}\otimes V^{\otimes q}\right)\;\rightarrow\; V^{*\otimes\left(m+p-\varrho\right)}\otimes V^{\otimes\left(n+q-\varrho\right)},\\
\left(\mathsf{s}_{\mathrm{a}_{1}\ldots\mathrm{a}_{m}}^{\mathrm{b}_{1}\ldots\mathrm{b}_{n}},\mathsf{t}_{\mathrm{c}_{1}\ldots\mathrm{d}_{p}}^{\mathrm{d}_{1}\ldots\mathrm{d}_{q}}\right)\mapsto\bullet{}_{\mathrm{a}_{1}\ldots\mathrm{a}_{\varrho}}^{\mathrm{d}_{1}\ldots\mathrm{d}_{\varrho}}\left(\mathsf{s}_{\mathrm{a}_{1}\ldots\mathrm{a}_{m}}^{\mathrm{b}_{1}\ldots\mathrm{b}_{n}},\mathsf{t}_{\mathrm{c}_{1}\ldots\mathrm{d}_{p}}^{\mathrm{d}_{1}\ldots\mathrm{d}_{q}}\right).
\end{gather*}
We have thus defined multiplication of $\boldsymbol{\mathsf{s}}$
and $\boldsymbol{\mathsf{t}}$ connecting the first $\varrho$ subscripts
of $\boldsymbol{\mathsf{s}}$ and the first $\varrho$ superscripts
of $\boldsymbol{\mathsf{t}}$.

As an alternative to this 'basis-free' argument, we could have noted
that \linebreak{}
$V^{*\otimes m}\otimes V^{\otimes n}$ has a basis of the form $\left\{ \bigotimes_{i=1}^{m}\boldsymbol{f}^{i}\otimes\bigotimes_{i=1}^{n}\mathbf{v}_{i}\right\} $,
$V^{*\otimes p}\otimes V^{\otimes q}$ has a basis of the form $\left\{ \bigotimes_{i=1}^{p}\boldsymbol{\bar{f}}\vphantom{}^{i}\otimes\bigotimes_{i=1}^{q}\bar{\mathbf{v}}_{i}\right\} $,
and $V^{*\otimes\left(m+p-\varrho\right)}\otimes V^{\otimes\left(n+q-\varrho\right)}$
has a basis of the form $\left\{ \widetilde{\beta}_{\bullet}^{\varrho}\left(\bigotimes_{i=1}^{m}\boldsymbol{f}^{i}\otimes\bigotimes_{i=1}^{n}\mathbf{v}_{i}\,,\,\bigotimes_{i=1}^{p}\boldsymbol{\bar{f}}\vphantom{}^{i}\otimes\bigotimes_{i=1}^{q}\bar{\mathbf{v}}_{i}\right)\right\} $.
Thus, $\beta_{\bullet}^{\varrho}$ can be obtained by bilinear extension
of $\widetilde{\beta}_{\bullet}^{\varrho}$, as illustrated in Subsection
4.2(1). This argument also makes it clear that $\beta_{\bullet}^{\varrho}$
is a surjective map.

\emph{(2)}. It is customary and convenient to denote $\bullet{}_{\mathrm{a}_{1}\ldots\mathrm{a}_{\varrho}}^{\mathrm{d}_{1}\ldots\mathrm{d}_{\varrho}}\left(\mathsf{\mathsf{s}}{}_{\mathrm{a}_{1}\ldots\mathrm{a}_{m}}^{\mathrm{b}_{1}\ldots\mathrm{b}_{n}},\mathsf{t}{}_{\mathrm{c}_{1}\ldots\mathrm{d}_{p}}^{\mathrm{d}_{1}\ldots\mathrm{d}_{q}}\right)$
by 
\[
\mathsf{s}{}_{\mathrm{a}_{1}\ldots\mathrm{a}_{m}}^{\mathrm{b}_{1}\ldots\mathrm{b}_{n}}\mathsf{t}{}_{\mathrm{c}_{1}\ldots\mathrm{d}_{p}}^{\mathrm{a}_{1}\ldots\mathrm{a}_{\varrho}\mathrm{d}_{\varrho+1}\ldots\mathrm{d}_{q}}\:,
\]
using \emph{matching indices} to indicate how copies of $V^{*}$ in
$V^{*m}\otimes V^{n}$ are paired with copies of $V$ in $V^{*p}\otimes V^{q}$.
Eliminating matching indices, the resulting bilateral tensor can be
written as 
\[
\mathsf{p}{}_{\mathrm{a}_{\varrho+1}\ldots\mathrm{a}_{m}\mathrm{c}_{1}\ldots\mathrm{d}_{p}}^{\mathrm{b}_{1}\ldots\mathrm{b}_{n}\mathrm{d}_{\varrho+1}\ldots\mathrm{d}_{q}}.
\]
For example, 
\[
\bullet\left(\mathsf{s}{}_{\mathrm{a}}^{\mathrm{b}},\mathsf{t}{}_{\mathrm{c}}^{\mathrm{d}}\right)=\mathsf{s}{}_{\mathrm{a}}^{\mathrm{b}}\mathsf{t}{}_{\mathrm{c}}^{\mathrm{d}}=\mathsf{p}{}_{\mathrm{ac}}^{\mathrm{bd}},\quad\bullet_{\mathrm{a}}^{\mathrm{d}}\left(\mathsf{s}{}_{\mathrm{a}}^{\mathrm{b}},\mathsf{t}{}_{\mathrm{c}}^{\mathrm{d}}\right)=\mathsf{s}{}_{\mathrm{a}}^{\mathrm{b}}\mathsf{t}{}_{\mathrm{c}}^{\mathrm{a}}=\mathsf{q}{}_{\mathrm{c}}^{\mathrm{b}},\quad\bullet{}_{\mathrm{ab}}^{\mathrm{de}}\left(\mathsf{s}{}_{\mathrm{a}\mathrm{b}}^{\mathrm{c}},\mathsf{t}{}^{\mathrm{d\mathrm{e}}}\right)=\mathsf{s}{}_{\mathrm{a}\mathrm{b}}^{\mathrm{c}}\mathsf{t}{}^{\mathrm{\mathrm{ab}}}=\mathfrak{\mathsf{r}{}^{\mathrm{c}}}.
\]

\emph{(3)}. It remains to define multiplication of $\boldsymbol{\mathsf{s}}$
and $\boldsymbol{\mathsf{t}}$ connecting any $\varrho$ pairs of
indices formed by one subscript of $\boldsymbol{\mathsf{s}}$ and
one superscript of $\boldsymbol{\mathsf{t}}$. We do this by replacing
the mapping $\widetilde{\beta}_{\bullet}^{\varrho}$ in (\ref{eq:bilmult})
by the more general mapping
\begin{gather}
\left(\bigotimes_{i=1}^{m}\boldsymbol{f}^{i}\otimes\bigotimes_{i=1}^{n}\mathbf{v}_{i}\,,\,\bigotimes_{i=1}^{p}\boldsymbol{\bar{f}}\vphantom{}^{i}\otimes\bigotimes_{i=1}^{q}\bar{\mathbf{v}}_{i}\right)\mapsto\label{eq:multgen}\\
\prod_{j=1}^{\varrho}\boldsymbol{f}^{i_{j}}\!\left(\bar{\mathbf{v}}_{i_{j}}\right)\bigotimes_{i\neq i_{j}}\boldsymbol{f}^{i}\otimes\bigotimes_{i=1}^{p}\boldsymbol{\bar{f}}\vphantom{}^{i}\otimes\bigotimes_{i=1}^{n}\mathbf{v}_{i}\otimes\bigotimes_{i\neq i_{j}}\bar{\mathbf{v}}_{i},\nonumber 
\end{gather}
where $\prod_{j=1}^{0}\boldsymbol{f}^{i_{j}}\!\left(\bar{\mathbf{v}}_{i_{j}}\right)=1$,
$\bigotimes_{i\neq i_{j}}\boldsymbol{f}^{i}=1$ if $\varrho\!=\! m\!>\!0$,
and $\bigotimes_{i\neq i_{j}}\bar{\mathbf{v}}_{i}=1$ if $\varrho\!=\! q\!>\!0$.
\begin{example*}
Set $\mathsf{\widetilde{s}}_{\mathrm{abc}}=\boldsymbol{f}{}_{\mathrm{a}}\otimes\boldsymbol{f}{}_{\mathrm{b}}\otimes\boldsymbol{f}{}_{\mathrm{c}}\otimes\xi$
and $\widetilde{\mathsf{t}}^{\mathrm{def}}=\bar{\phi}\otimes\bar{\mathbf{v}}_{\mathrm{d}}\otimes\bar{\mathbf{v}}_{\mathrm{e}}\otimes\bar{\mathbf{v}}_{\mathrm{f}}$.
Then we have 
\[
\bullet_{\mathrm{ac}}^{\mathrm{fe}}\left(\mathsf{\widetilde{s}}_{\mathrm{abc}},\widetilde{\mathsf{t}}^{\mathrm{def}}\right)=\boldsymbol{f}_{\mathrm{a}}\!\left(\bar{\mathbf{v}}_{\mathrm{f}}\right)\boldsymbol{f}_{\mathrm{c}}\!\left(\bar{\mathbf{v}}_{\mathrm{e}}\right)\boldsymbol{f}{}_{\mathrm{b}}\otimes\bar{\phi}\otimes\xi\otimes\bar{\mathbf{v}}_{\mathrm{d}}=\boldsymbol{f}_{\mathrm{a}}\!\left(\bar{\mathbf{v}}_{\mathrm{f}}\right)\boldsymbol{f}_{\mathrm{c}}\!\left(\bar{\mathbf{v}}_{\mathrm{e}}\right)\xi\boldsymbol{f}{}_{\mathrm{b}}\otimes\bar{\phi}\bar{\mathbf{v}}_{\mathrm{d}}.
\]
 Thus, $\mathsf{\widetilde{s}}_{\mathrm{abc}}^{1}=\boldsymbol{f}{}_{\mathrm{a}}\otimes\boldsymbol{f}{}_{\mathrm{b}}\otimes\boldsymbol{f}{}_{\mathrm{c}}\otimes1$
and $\widetilde{\mathsf{t}}_{1}^{\mathrm{def}}=1\otimes\bar{\mathbf{v}}_{\mathrm{d}}\otimes\bar{\mathbf{v}}_{\mathrm{e}}\otimes\bar{\mathbf{v}}_{\mathrm{f}}$
gives
\[
\bullet_{\mathrm{ac}}^{\mathrm{fe}}\left(\mathsf{\widetilde{s}}_{\mathrm{abc}}^{1},\widetilde{\mathsf{t}}_{1}^{\mathrm{def}}\right)=\boldsymbol{f}_{\mathrm{a}}\!\left(\bar{\mathbf{v}}_{\mathrm{f}}\right)\boldsymbol{f}_{\mathrm{c}}\!\left(\bar{\mathbf{v}}_{\mathrm{e}}\right)\boldsymbol{f}{}_{\mathrm{b}}\otimes1\otimes1\otimes\bar{\mathbf{v}}_{\mathrm{d}}=\boldsymbol{f}_{\mathrm{a}}\!\left(\bar{\mathbf{v}}_{\mathrm{f}}\right)\boldsymbol{f}_{\mathrm{c}}\!\left(\bar{\mathbf{v}}_{\mathrm{e}}\right)\boldsymbol{f}{}_{\mathrm{b}}\otimes\bar{\mathbf{v}}_{\mathrm{d}}.\quad\square
\]

\end{example*}
It is clear how to define a general bilinear function
\[
\left(\mathsf{s}_{\mathrm{a}_{1}\ldots\mathrm{a}_{m}}^{\mathrm{b}_{1}\ldots\mathrm{b}_{n}},\mathsf{t}_{\mathrm{c}_{1}\ldots\mathrm{d}_{p}}^{\mathrm{d}_{1}\ldots\mathrm{d}_{q}}\right)\mapsto\bullet{}_{\mathrm{a}_{i_{1}}\ldots\mathrm{a}_{i_{\varrho}}}^{\mathrm{d}_{j_{1}}\ldots\mathrm{d}_{j_{\varrho}}}\left(\mathsf{s}_{\mathrm{a}_{1}\ldots\mathrm{a}_{m}}^{\mathrm{b}_{1}\ldots\mathrm{b}_{n}},\mathsf{t}_{\mathrm{c}_{1}\ldots\mathrm{d}_{p}}^{\mathrm{d}_{1}\ldots\mathrm{d}_{q}}\right).
\]
We have thus introduced\emph{ multiplication of bilateral tensors
in the general case.}

The notation with matching indices turns out to be convenient in the
general case, too; we can define complicated products of bilateral
tensors such as 
\[
\bullet_{\mathrm{ac}}^{\mathrm{fe}}\left(\mathsf{s}_{\mathrm{abc}},\mathsf{t}^{\mathrm{def}}\right)=\mathsf{s}_{\mathrm{abc}}\mathsf{t}^{\mathrm{dca}}=\mathsf{p}_{\mathrm{b}}^{\mathrm{d}},\qquad\bullet_{\mathrm{ab}}^{\mathrm{ec}}\left(\mathsf{s}{}_{\mathrm{ab}},\mathsf{t}{}_{\mathrm{f}}^{\mathrm{cde}}\right)=\mathsf{s}{}_{\mathrm{ab}}\mathsf{t}{}_{\mathrm{f}}^{\mathrm{bda}}=\mathsf{q}{}_{\mathrm{f}}^{\mathrm{d}}.
\]

\emph{(4)}. We call
\[
\bullet\left(\mathsf{s}{}_{\mathrm{a}_{1}\ldots\mathrm{a}_{m}}^{\mathrm{b}_{1}\ldots\mathrm{b}_{n}},\mathsf{t}{}_{\mathrm{c}_{1}\ldots\mathrm{d}_{p}}^{\mathrm{d}_{1}\ldots\mathrm{d}_{q}}\right)=\mathsf{s}{}_{\mathrm{a}_{1}\ldots\mathrm{a}_{m}}^{\mathrm{b}_{1}\ldots\mathrm{b}_{n}}\mathsf{t}{}_{\mathrm{c}_{1}\ldots\mathrm{d}_{p}}^{\mathrm{d}_{1}\ldots\mathrm{d}_{q}}=\mathsf{p}{}_{\mathrm{a}_{1}\ldots\mathrm{a}_{m}\mathrm{c}_{1}\ldots\mathrm{d}_{p}}^{\mathrm{b}_{1}\ldots\mathrm{b}_{n}\mathrm{d}_{1}\ldots\mathrm{d}_{q}}\,,
\]
in index-free form denoted by $\mathsf{\boldsymbol{\mathsf{s}}}\boldsymbol{\mathsf{t}}$,
the \emph{outer product} of $\boldsymbol{\mathsf{s}}$ and $\boldsymbol{\mathsf{t}}$.
\begin{rem*}
The 'outer product' is of course the usual tensor product. Note that
the outer product $\mathsf{\boldsymbol{\mathsf{s}}}\boldsymbol{\mathsf{t}}$
of $\mathsf{\boldsymbol{\mathsf{s}}},\boldsymbol{\mathsf{t}}\in V^{*}\otimes V$
belongs to $V^{*}\otimes V^{*}\otimes V\otimes V$ rather than $\left(V^{*}\otimes V\right)\otimes\left(V^{*}\otimes V\right)$,
because $\mathsf{\boldsymbol{\mathsf{s}}}\boldsymbol{\mathsf{t}}$
is also a bilateral tensor. Thus, we do not in general have $\mathsf{\boldsymbol{\mathsf{s}}}\boldsymbol{\mathsf{t}}=\boldsymbol{\mathsf{s}}\otimes\boldsymbol{\mathsf{t}}$,
only $\mathsf{\boldsymbol{\mathsf{s}}}\boldsymbol{\mathsf{t}}\cong\boldsymbol{\mathsf{s}}\otimes\boldsymbol{\mathsf{t}}$.
\end{rem*}
\emph{(5)}. Set $\boldsymbol{\mathsf{r}}=\boldsymbol{f}\otimes\mathbf{u}$,
$\boldsymbol{\mathsf{s}}=\boldsymbol{g}\otimes\mathbf{v}$ and $\boldsymbol{\mathsf{t}}=\boldsymbol{h}\otimes\mathbf{w}$.
Then 
\begin{gather*}
\left(\boldsymbol{\mathsf{r}}\boldsymbol{\mathsf{s}}\right)\boldsymbol{\mathsf{t}}=\left(\boldsymbol{f}\otimes\boldsymbol{g}\otimes\mathbf{u}\otimes\mathbf{v}\right)\left(\boldsymbol{h}\otimes\mathbf{w}\right)=\left(\left(\boldsymbol{f}\otimes\boldsymbol{g}\right)\otimes\boldsymbol{h}\right)\otimes\left(\left(\mathbf{u}\otimes\mathbf{v}\right)\otimes\mathbf{w}\right)=\\
\left(\boldsymbol{f}\otimes\left(\boldsymbol{g}\otimes\boldsymbol{h}\right)\right)\otimes\left(\mathbf{u}\otimes\left(\mathbf{v}\otimes\mathbf{w}\right)\right)=\left(\boldsymbol{f}\otimes\mathbf{u}\right)\left(\boldsymbol{g}\otimes\boldsymbol{h}\otimes\mathbf{v}\otimes\mathbf{w}\right)=\boldsymbol{\mathsf{r}}\left(\boldsymbol{\mathsf{s}}\boldsymbol{\mathsf{t}}\right),
\end{gather*}
and it is shown in the same way that $\left(\boldsymbol{\mathsf{r}}\boldsymbol{\mathsf{s}}\right)\boldsymbol{\mathsf{t}}=\boldsymbol{\mathsf{r}}\left(\boldsymbol{\mathsf{s}}\boldsymbol{\mathsf{t}}\right)$
if $\boldsymbol{\mathsf{r}}=\left(\bigotimes_{i}\boldsymbol{f}_{i}\right)\otimes\left(\bigotimes_{i}\mathbf{u}_{i}\right)$
etc. This can be generalized to general multiplication of simple bilateral
tensors and finally to general multiplication of general bilateral
tensors. In other words, since the tensor product is associative,
multiplication of bilateral tensors is associative as well.

\subsection{Contraction of bilateral tensors}

$ $\\
\emph{(1)}. Let $\varrho,m,n$ be non-negative integers such that
$\varrho\leq m,n$. The mapping

\begin{equation}
\widetilde{\lambda}:\left(\bigotimes_{i=1}^{m}\boldsymbol{f}^{i}\right)\otimes\left(\bigotimes_{i=1}^{n}\mathbf{v}_{i}\right)\;\mapsto\;\prod_{i=1}^{\varrho}\boldsymbol{f}^{i}\!\left(\mathbf{v}{}_{i}\right)\left(\bigotimes_{i=\varrho+1}^{m}\!\boldsymbol{f}^{i}\right)\otimes\left(\bigotimes_{i=\varrho+1}^{n}\!\mathsf{\mathbf{v}}_{i}\right),\label{eq:contrbil}
\end{equation}
where $\prod_{i=1}^{0}\boldsymbol{f}^{i}\!\left(\mathbf{v}{}_{i}\right)=1$,
$\bigotimes_{i=\varrho+1}^{\varrho>0}\boldsymbol{f}^{i}=1$ and $\bigotimes_{i=\varrho+1}^{\varrho>0}\mathbf{v}{}_{i}=1$,
has a corresponding separately linear map
\[
\mu:V^{*m}\times V^{n}\rightarrow V^{*\otimes\left(m-\varrho\right)}\otimes V^{\otimes\left(n-\varrho\right)}
\]
such that 
\[
\mu\left(\left.\prod\right._{i=1}^{m}\boldsymbol{f}^{i}\,,\,\left.\prod\right._{i=1}^{n}\mathbf{v}_{i}\right)=\widetilde{\lambda}\left(\left(\left.\bigotimes\right._{i=1}^{m}\boldsymbol{f}^{i}\right)\otimes\left(\left.\bigotimes\right._{i=1}^{n}\mathbf{v}_{i}\right)\right),
\]
so $\widetilde{\lambda}$ gives a linear map
\begin{gather*}
\mathfrak{\mathbf{\blacklozenge}{}_{\mathrm{a}_{1}\ldots\mathrm{a}_{\varrho}}^{\mathrm{b}_{1}\ldots\mathrm{b}_{\varrho}}}:V^{*\otimes m}\otimes V^{\otimes n}\rightarrow V^{*\otimes\left(m-\varrho\right)}\otimes V^{\otimes\left(n-\varrho\right)},\\
\mathsf{t}_{\mathrm{a}_{1}\ldots\mathrm{a}_{m}}^{\mathrm{b}_{1}\ldots\mathrm{b}_{n}}\mapsto\mathbf{\blacklozenge}{}_{\mathrm{a}_{1}\ldots\mathrm{a}_{\varrho}}^{\mathrm{b}_{1}\ldots\mathrm{b}_{\varrho}}\left(\mathsf{t}_{\mathrm{a}_{1}\ldots\mathrm{a}_{m}}^{\mathrm{b}_{1}\ldots\mathrm{b}_{n}}\right).
\end{gather*}
$\mathbf{\blacklozenge}{}_{\mathrm{a}_{1}\ldots\mathrm{a}_{\varrho}}^{\mathrm{b}_{1}\ldots\mathrm{b}_{\varrho}}\left(\mathsf{t}_{\mathrm{a}_{1}\ldots\mathrm{a}_{m}}^{\mathrm{b}_{1}\ldots\mathrm{b}_{n}}\right)$
is said to be the \emph{contraction} of $\mathsf{t}_{\mathrm{a}_{1}\ldots\mathrm{a}_{m}}^{\mathrm{b}_{1}\ldots\mathrm{b}_{n}}$
with respect to $\mathrm{\mathsf{a}}_{1}\ldots\mathrm{\mathsf{a}}_{\varrho}$
and $\mathrm{\mathsf{b}}_{1}\ldots\mathrm{\mathsf{b}}_{\varrho}$.
(One may choose to reserve the term 'contraction' for the case where
$\varrho>0$.) 

It is clear that we can introduce a general contraction mapping
\[
\mathsf{t}_{\mathrm{a}_{1}\ldots\mathrm{a}_{m}}^{\mathrm{b}_{1}\ldots\mathrm{b}_{n}}\mapsto\mathbf{\blacklozenge}{}_{\mathrm{a}_{i_{1}}\ldots\mathrm{a}_{i_{\varrho}}}^{\mathrm{b}_{j_{1}}\ldots\mathrm{b}_{j_{\varrho}}}\left(\mathsf{t}_{\mathrm{a}_{1}\ldots\mathrm{a}_{m}}^{\mathrm{b}_{1}\ldots\mathrm{b}_{n}}\right)
\]
defined by
\begin{equation}
\widetilde{\lambda}:\left(\bigotimes_{i=1}^{m}\boldsymbol{f}^{i}\right)\otimes\left(\bigotimes_{i=1}^{n}\mathbf{v}_{i}\right)\;\mapsto\;\prod_{j=1}^{\varrho}\boldsymbol{f}^{i_{j}}\!\left(\mathbf{v}_{i_{j}}\right)\left(\bigotimes_{i\neq i_{j}}\boldsymbol{f}^{i}\right)\otimes\left(\bigotimes_{i\neq i_{j}}\mathbf{v}_{i}\right),\label{eq:contract2}
\end{equation}
where $\prod_{j=1}^{0}\boldsymbol{f}^{i_{j}}\!\left(\bar{\mathbf{v}}_{i_{j}}\right)=1$,
$\bigotimes_{i\neq i_{j}}\boldsymbol{f}^{i}=1$ if $\varrho\!=\! m\!>\!0$,
and $\bigotimes_{i\neq i_{j}}\bar{\mathbf{v}}_{i}=1$ if $\varrho\!=\! q\!>\!0$.

It is customary to use matching indices to specify contractions. For
example, 
\[
\blacklozenge{}_{\mathrm{a}}^{\mathrm{b}}\left(\mathsf{r}_{\mathrm{a}}^{\mathrm{b}}\right)=\mathsf{r}_{\mathrm{a}}^{\mathrm{a}}=\widehat{\mathsf{r}},\quad\blacklozenge{}_{\mathrm{ab}}^{\mathrm{dc}}\left(\mathsf{s}_{\mathrm{ab}}^{\mathrm{cd}}\right)=\mathsf{s}_{\mathrm{ab}}^{\mathrm{ba}}=\widehat{\mathsf{s}},\quad\blacklozenge{}_{\mathrm{b}}^{\mathrm{c}}\left(\mathsf{t}_{\mathrm{ab}}^{\mathrm{cd}}\right)=\mathsf{t}_{\mathrm{ab}}^{\mathrm{bd}}=\widehat{\mathsf{t}}_{\mathrm{a}}^{\mathrm{d}}.
\]
\emph{(2)}. Let us consider another way of defining contractions.
Set 
\[
\widetilde{\boldsymbol{\mathsf{t}}}=\left(\left.\bigotimes\right._{i=1}^{m}\boldsymbol{f}^{i}\right)\otimes\left(\left.\bigotimes\right._{i=1}^{n}\mathbf{v}{}_{i}\right),\quad\widetilde{\boldsymbol{\mathsf{p}}}=\left(\left.\bigotimes\right._{i=1}^{m}\boldsymbol{f}^{i}\right)\otimes1,\quad\widetilde{\boldsymbol{\mathsf{q}}}=1\otimes\left(\left.\bigotimes\right._{i=1}^{n}\mathbf{v}{}_{i}\right).
\]
Then $\widetilde{\boldsymbol{\mathsf{t}}}=\widetilde{\boldsymbol{\mathsf{p}}}\widetilde{\boldsymbol{\mathsf{q}}}$,
and it is clear that any simple bilateral tensor $\widetilde{\boldsymbol{\mathsf{t}}}$
can be written in a unique way as the outer product $\widetilde{\boldsymbol{\mathsf{p}}}\widetilde{\boldsymbol{\mathsf{q}}}$
of two simple bilateral tensors of the forms shown; in extended double-index
notation we have 
\[
\widetilde{\mathsf{t}}_{\mathrm{a}_{1}\ldots\mathrm{a}_{m}}^{\mathrm{b}_{1}\ldots\mathrm{b}_{n}}=\mathbf{\bullet}\left(\widetilde{\mathsf{p}}_{\mathrm{a}_{1}\ldots\mathrm{a}_{m}}^{1},\widetilde{\mathsf{q}}_{1}^{\mathrm{b}_{1}\ldots\mathrm{b}_{n}}\right)=\widetilde{\mathsf{p}}_{\mathrm{a}_{1}\ldots\mathrm{a}_{m}}^{1}\widetilde{\mathsf{q}}_{1}^{\mathrm{b}_{1}\ldots\mathrm{b}_{n}}.
\]
It follows from (\ref{eq:multgen}) that the mapping
\[
\widetilde{\lambda'}:\widetilde{\mathsf{t}}_{\mathrm{a}_{1}\ldots\mathrm{a}_{m}}^{\mathrm{b}_{1}\ldots\mathrm{b}_{n}}=\mathbf{\bullet}\left(\widetilde{\mathsf{p}}_{\mathrm{a}_{1}\ldots\mathrm{a}_{m}}^{1},\widetilde{\mathsf{q}}_{1}^{\mathrm{b}_{1}\ldots\mathrm{b}_{n}}\right)\mapsto\mathbf{\bullet}{}_{\mathrm{a}_{i_{1}}\ldots\mathrm{a}_{i_{\varrho}}}^{\mathrm{b}_{j_{1}}\ldots\mathrm{b}_{j_{\varrho}}}\left(\widetilde{\mathsf{p}}_{\mathrm{a}_{1}\ldots\mathrm{a}_{m}}^{1},\widetilde{\mathsf{q}}_{1}^{\mathrm{b}_{1}\ldots\mathrm{b}_{n}}\right)
\]
is actually the mapping $\widetilde{\lambda}$ defined in (\ref{eq:contract2}).
Thus, the linear map given by $\widetilde{\lambda'}$ is precisely
the contraction map $\mathbf{\blacklozenge}{}_{\mathrm{a}_{i_{1}}\ldots\mathrm{a}_{i_{\varrho}}}^{\mathrm{b}_{j_{1}}\ldots\mathrm{b}_{j_{\varrho}}}$.
\begin{rem*}
General multiplication of bilateral tensors is usually defined in
terms of outer multiplication and contraction, but as we have seen
it is also possible to define outer multiplication and contraction
in terms of general multiplication. 
\end{rem*}

\section{Classical tensors as tensor maps}
\begin{quote}
Differential geometers are very fond of the isomorphism $L\left(E,E\right)\rightarrow E^{\vee}\otimes E,$
and often use $E^{\vee}\otimes E$ when they think geometrically of
$L\left(E,E\right)$, thereby emphasizing an unnecessary dualization,
and an irrelevant formalism, when it is easier to deal directly with
$L\left(E,E\right)$.\\
\emph{$\hphantom{xxxxxxxxxxxxxxxxxxxxxxxxxxxxx}$Serge Lang} \cite{key-2}
(p. 628).
\end{quote}

\subsection{Linear tensor maps}

$ $\\
There is a well-known isomorphism $V^{*}\otimes V\rightarrow\mathscr{L}\left[V,V\right]$,
and as pointed out by Lang, maps in $\mathscr{L}\left[V,V\right]$
seem to be more natural objects than elements of $V^{*}\otimes V$.
However, the same can be said about maps in $\mathscr{L}\left[V^{\otimes m},V^{\otimes n}\right]$
versus elements of $V^{*\otimes m}\otimes V^{\otimes n}$. We shall
now make a conceptual leap to a position where classical tensors are
actually not seen as (bilateral) tensors but as linear maps between
certain tensor product spaces.

\emph{(1)}. A\emph{ linear tensor map} $\boldsymbol{t}$ on a vector
space $V$ over $K$ is a linear function
\begin{equation}
\boldsymbol{t}:V^{\otimes m}\rightarrow V^{\otimes n}\quad(m,n\geq0).\label{eq:tensmap}
\end{equation}
Instead of this \emph{index-free notation}, we can use \emph{double-index
notation for linear tensor maps}, writing 
\begin{equation}
t_{\mathrm{a}_{1}\ldots\mathrm{a}_{m}}^{\mathrm{b}_{1}\ldots\mathrm{b}_{n}}:V_{\mathrm{a}_{1}}\otimes\ldots\otimes V_{\mathrm{a}_{m}}\rightarrow V_{\mathrm{b}_{1}}\otimes\ldots\otimes V_{\mathrm{b}_{n}}\quad(V_{\mathrm{a}_{i}},V_{\mathrm{b}_{j}}=V),\label{eq:tensmapdind}
\end{equation}
where no index occurs twice. The indices thus identify copies of $V$;
subscripts identify copies of $V$ associated with 'inputs' to $\boldsymbol{t}$,
while superscripts identify copies of $V$ associated with 'outputs'
from $\boldsymbol{t}$. As in the case of bilateral tensors, a tensor
map with $m$ subscripts and $n$ superscripts is said to have valence
$\binom{n}{m}$.

We use Roman letters rather than italics to suggest that indices are
labels identifying 'slots' for arguments (inputs) or values (outputs)
rather than indices identifying scalars in systems of scalars.

Italics are used to identify particular scalars, vectors, tensors
or tensor maps in collections of such objects, as when we write $s_{i},\mathbf{v}{}_{i},\boldsymbol{\mathsf{t}}{}_{i},v_{i}^{\mathrm{a}},\boldsymbol{t}_{i},t_{\mathrm{a}i}^{\mathrm{b}}$
etc. 

As usual, a sequence of zero or more indices can be replaced by a
multi-index such as $\mathrm{A}=\left(\mathsf{a}_{1},\ldots,\mathsf{a}_{n}\right)$
or $I=\left(i,j,k\right)$.

\emph{(2)}. A linear tensor map $s:K\rightarrow K$ is said to be
\emph{scalar-like}. Since $s$ is linear, $s\!\left(\eta\right)=\eta s\!\left(1\right)=\eta\sigma$,
so we can identify $s$ with the scalar $\sigma$. Similarly, a linear
tensor map $v^{\mathrm{a}}:K\rightarrow V$ is said to be\emph{ vector-like}.
In this case, $v^{\mathrm{a}}\!\left(\eta\right)=\eta v^{\mathrm{a}}\!\left(1\right)=\eta\mathbf{v}$,
so it is natural to identify $v^{\mathrm{a}}$ with the vector $\mathbf{v}$.
Finally, a linear tensor map $f_{\mathrm{a}}:V\rightarrow K$ is a
\emph{linear form} $\mathbf{v}\mapsto\boldsymbol{f}\!\left(\mathbf{v}\right)$.

\subsection{Linear and separately linear tensor maps}

$ $\\
As mentioned in Subsection 5.2, the sets of linear maps $\mathscr{L}\left[V^{\otimes m},V^{\otimes n}\right]$
and separately linear maps $\mathscr{L}\left[V^{m},V^{\otimes n}\right]$
can be regarded as vector spaces. We also know that there is a canonical
isomorphism
\[
\Lambda_{\mathcal{M}}^{-1}:\mathscr{L}\left[V^{\otimes m},V^{\otimes n}\right]\rightarrow\mathscr{L}\left[V^{m},V^{\otimes n}\right].
\]
The separately linear map $\Lambda_{\mathcal{M}}^{-1}\left(t_{\mathrm{a}_{1}\ldots\mathrm{a}_{m}}^{\mathrm{b}_{1}\ldots\mathrm{b}_{n}}\right)$
will be denoted by $t_{\mathrm{a}_{1},\ldots,\mathrm{a}_{m}}^{\mathrm{b}_{1}\ldots\mathrm{b}_{n}}$
(double-index notation) or $\boldsymbol{t}$ (index-free notation).
For example, $g_{\mathrm{ab}}$ is a linear scalar-valued tensor map,
while $g_{\mathrm{a},\mathrm{b}}$ is a separately linear (bilinear)
scalar-valued tensor map. (This notation does not lead to ambiguity,
since a separately linear tensor map with $0$ or $1$ subscripts
is in effect a linear tensor map.)

Since a separately linear tensor map can be identified with the corresponding
linear tensor map, we can use the term \emph{tensor map} to refer
to both without distinction.

\subsection{Bases for spaces of tensor maps}

$ $\\
Let $U$ and $V$ be $N$-dimensional vector spaces with bases $\left\{ \mathbf{e}_{k}\right\} $
and $\left\{ \mathbf{f}_{j}\right\} $, respectively. Any $\boldsymbol{t}\in\mathscr{L}\left[U,V\right]$
is given by expansions of the form $\boldsymbol{t}\!\left(\mathbf{e}_{k}\right)=\sum_{j}\alpha^{kj}\mathbf{f}_{j}$
for all $\mathbf{e}_{k}$. Let $\boldsymbol{T}_{ij}\in\mathscr{L}\left[U,V\right]$
be given by $\boldsymbol{T}_{ij}\!\left(\mathbf{e}_{k}\right)=\mathbf{f}_{j}$
if $i=k$ and $\boldsymbol{T}_{ij}\!\left(\mathbf{e}_{k}\right)=0$
if $i\neq k$ for all $\mathbf{e}_{k}$. Then $\alpha^{kj}\mathbf{f}_{j}=\sum_{i}\alpha^{ij}\boldsymbol{T}_{ij}\!\left(\mathbf{e}_{k}\right)$,
so $\boldsymbol{t}\!\left(\mathbf{e}_{k}\right)=\sum_{i,j}\alpha^{ij}\boldsymbol{T}_{ij}\!\left(\mathbf{e}_{k}\right)$
for all $\mathbf{e}_{k}$, so $\boldsymbol{t}\!\left(\mathbf{u}\right)=\sum_{i,j}\alpha^{ij}\boldsymbol{T}_{ij}\!\left(\mathbf{u}\right)$
for all $\boldsymbol{t}\in\mathscr{L}\left[U,V\right]$ and $\mathbf{u}\in U$,
and as the coefficients $\alpha^{ij}$ are uniquely determined this
implies that $\left\{ \boldsymbol{T}_{ij}\right\} $ is a basis for
$\mathscr{L}\left[U,V\right]$.

The fact that $\boldsymbol{T}_{ij}\!\left(\mathbf{e}_{i}\right)=\mathbf{f}_{j}$
but $\boldsymbol{T}_{ij}\!\left(\mathbf{e}_{k}\right)=0$ if $i\neq k$
suggests that we write $\boldsymbol{T}_{ij}$ as $\left\langle \mathbf{e}_{i}\mapsto\mathbf{f}_{j}\right\rangle $.
With this notation, $\mathscr{L}\left[U,V\right]$ has a basis with
$N^{2}$ elements of the form
\[
\left\{ \left\langle \mathbf{e}_{i}\mapsto\mathbf{f}_{j}\right\rangle \right\} ,
\]
and as $\left\{ \mathbf{e}_{i_{1}}\otimes\ldots\otimes\mathbf{e}_{i_{m}}\right\} $
is a basis for $V^{\otimes m}$ and $\left\{ \mathbf{e}_{i_{1}}\otimes\ldots\otimes\mathbf{e}_{i_{n}}\right\} $
a basis for $V^{\otimes n}$ we conclude that
\[
\left\{ \left\langle \mathbf{e}_{i_{1}}\otimes\ldots\otimes\mathbf{e}_{i_{m}}\mapsto\mathbf{e}_{j_{1}}\otimes\ldots\otimes\mathbf{e}_{j_{n}}\right\rangle \right\} 
\]
is a basis with $N^{m+n}$ elements for $\mathscr{L}\left[V^{\otimes m},V^{\otimes n}\right]$. 

Since $\left\{ 1\right\} $ is basis for $K$, $\left\{ \left\langle 1\mapsto1\right\rangle \right\} $
is a basis for $\mathscr{L}\left[K,K\right]$. By definition,\linebreak{}
 $\left\langle 1\mapsto1\right\rangle \left(1\right)=1$, so by linearity
$\left\langle 1\mapsto1\right\rangle \left(\eta\right)=\eta$, so
$\left\langle 1\mapsto1\right\rangle =\left(\eta\mapsto\eta\right)$.
Similarly, if $\left\{ \mathbf{e}_{i}\right\} $ is a basis for $V$
then $\left\{ \left\langle 1\mapsto\mathbf{e}_{i}\right\rangle \right\} =\left\{ \left(\eta\mapsto\eta\mathbf{e}_{i}\right)\right\} $
is a basis for $\mathscr{L}\left[K,V\right]$.

\subsection{Bilateral tensors and corresponding tensor maps}

$ $\\
\emph{(1)}. Consider the mapping
\[
\widetilde{\Lambda}_{\mathcal{B}}:\left.\left(V^{*\otimes m}\otimes V^{\otimes n}\right)\right|_{\mathcal{S}}\rightarrow\mathscr{L}\left[V^{\otimes m},V^{\otimes n}\right],\quad\widetilde{\boldsymbol{\mathsf{t}}}\mapsto\widetilde{\boldsymbol{t}},
\]
 where 
\[
\widetilde{\boldsymbol{\mathsf{t}}}=\left(\left.\bigotimes\right._{i=1}^{m}\boldsymbol{f}_{i}\right)\otimes\left(\left.\bigotimes\right._{i=1}^{n}\mathbf{v}_{i}\right),\quad\widetilde{\boldsymbol{t}}=\left(\left.\bigotimes\right._{i=1}^{m}\mathbf{u}_{i}\mapsto\left.\prod\right._{i=1}^{m}\boldsymbol{f}_{i}\!\left(\mathbf{u}_{i}\right)\left.\bigotimes\right._{j=1}^{n}\mathbf{v}_{j}\right),
\]
and $\prod_{i=1}^{0}\boldsymbol{f}_{i}\!\left(\mathbf{u}_{i}\right)=\phi\eta$,
where $\phi,\eta\in K$ are scalar variables.

Inspection of $\widetilde{\boldsymbol{\mathsf{t}}}$ and $\widetilde{\boldsymbol{t}}$
reveals that the mapping $\mu:V^{*m}\times V^{n}\rightarrow\mathscr{L}\left[V^{\otimes m},V^{\otimes n}\right]$
defined by setting 
\[
\mu\left(\left.\prod\right._{i=1}^{m}\boldsymbol{f}_{i}\,,\,\left.\prod\right._{i=1}^{n}\mathbf{v}_{i}\right)=\widetilde{\Lambda}_{\mathcal{B}}\left(\left(\left.\bigotimes\right._{i=1}^{m}\boldsymbol{f}_{i}\right)\otimes\left(\left.\bigotimes\right._{i=1}^{n}\mathbf{v}_{i}\right)\right)
\]
is separately linear, so there is a corresponding linear map
\[
\Lambda_{\mathcal{B}}:V^{*m}\otimes V^{n}\rightarrow\mathscr{L}\left[V^{\otimes m},V^{\otimes n}\right],\qquad\boldsymbol{\mathsf{t}}\mapsto\Lambda_{\mathcal{B}}\!\left(\boldsymbol{\mathsf{t}}\right)=\boldsymbol{t}.
\]
This 'basis-free' argument does not prove that $\Lambda_{\mathcal{B}}$
is an isomorphism, however, so let us introduce a basis $\left\{ \mathbf{e}_{1},\ldots,\mathbf{e}_{N}\right\} $
for $V$ and a dual basis $\left\{ \boldsymbol{e}^{1},\ldots,\boldsymbol{e}^{N}\right\} $
for $V^{*}$. Then 
\[
\left\{ \boldsymbol{e}^{i_{1}}\otimes\ldots\otimes\boldsymbol{e}^{i_{m}}\otimes\mathbf{e}{}_{j_{1}}\otimes\ldots\otimes\mathbf{e}{}_{j_{n}}\mid\boldsymbol{e}^{i_{k}}\in\left\{ \boldsymbol{e}^{1},\ldots,\boldsymbol{e}^{N}\right\} ,\mathbf{e}{}_{j_{\ell}}\in\left\{ \mathbf{e}_{1},\ldots,\mathbf{e}_{N}\right\} \right\} 
\]
 is a basis with $N^{m+n}$ elements for $V^{*\otimes m}\otimes V^{\otimes n}$.
By definition,
\[
\widetilde{\Lambda}_{\mathcal{B}}\!\left(\boldsymbol{e}^{i_{1}}\otimes\ldots\otimes\boldsymbol{e}^{i_{m}}\otimes\mathbf{e}{}_{j_{1}}\otimes\ldots\otimes\mathbf{e}{}_{j_{n}}\right)\!\left(\mathbf{e}{}_{k_{1}}\otimes\ldots\otimes\mathbf{e}{}_{k_{m}}\right)=\prod_{\ell=1}^{m}\boldsymbol{e}^{_{i_{\ell}}}\!\left(\mathbf{e}{}_{k_{\ell}}\right)\mathbf{e}{}_{j_{1}}\otimes\ldots\otimes\mathbf{e}{}_{j_{n}},
\]
and as $\boldsymbol{e}^{_{i_{\ell}}}\!\left(\mathbf{e}{}_{k_{\ell}}\right)=1$
if $i_{\ell}=k_{\ell}$ but $\boldsymbol{e}^{_{i_{\ell}}}\!\left(\mathbf{e}{}_{k_{\ell}}\right)=0$
if $i_{\ell}\neq k_{\ell}$, this means that
\[
\widetilde{\Lambda}_{\mathcal{B}}\left(\boldsymbol{e}^{i_{1}}\otimes\ldots\otimes\boldsymbol{e}^{i_{m}}\otimes\mathbf{e}{}_{j_{1}}\otimes\ldots\otimes\mathbf{e}{}_{j_{n}}\right)=\left\langle \mathbf{e}{}_{i_{1}}\otimes\ldots\otimes\mathbf{e}{}_{i_{m}}\mapsto\mathbf{e}{}_{j_{1}}\otimes\ldots\otimes\mathbf{e}{}_{j_{n}}\right\rangle .
\]
$\widetilde{\Lambda}_{\mathcal{B}}$ thus maps a basis for $V^{*\otimes m}\otimes V^{\otimes n}$
bijectively to a basis for $\mathscr{L}\left[V^{\otimes m},V^{\otimes n}\right]$,
so $\widetilde{\Lambda}_{\mathcal{B}}$ can be extended by linearity
to an isomorphism $\Lambda_{\mathcal{B}}:V^{*\otimes m}\otimes V^{\otimes n}\rightarrow\mathscr{L}\left[V^{\otimes m},V^{\otimes n}\right]$.
\begin{rem*}
The proof that $V^{*\otimes m}\otimes V^{\otimes n}\cong\mathscr{L}\left[V^{\otimes m},V^{\otimes n}\right]$
does not work if $V$ is infinite-dimensional, since the dual set
$\left\{ \boldsymbol{e}^{i}\right\} $ is not a basis for $V^{*}$
in that case. No results specifically concerned with the infinite-dimensional
case will be presented in this article.
\end{rem*}
\emph{(2)}. As an illustration, we have the following mappings 
\begin{gather*}
K\otimes K\ni\phi\otimes\xi\overset{\widetilde{\Lambda}_{\mathcal{B}}}{\longmapsto}\left(\eta\mapsto\left(\phi\eta\right)\xi\right)=\left(\eta\mapsto\eta\left(\phi\xi\right)\right)\in\mathscr{L}\left[K,K\right],\\
K\otimes V\ni\phi\otimes\mathbf{v}\overset{\widetilde{\Lambda}_{\mathcal{B}}}{\longmapsto}\left(\eta\mapsto\left(\phi\eta\right)\mathbf{v}\right)=\left(\eta\mapsto\eta\left(\phi\mathbf{v}\right)\right)\in\mathscr{L}\left[K,V\right],\\
V^{*}\otimes K\ni\boldsymbol{f}\otimes\xi\overset{\widetilde{\Lambda}_{\mathcal{B}}}{\longmapsto}\left(\mathbf{u}\mapsto\boldsymbol{f}\!\left(\mathbf{u}\right)\xi\right)=\left(\mathbf{u}\mapsto\left(\xi\boldsymbol{f}\right)\!\left(\mathbf{u}\right)\right)\in\mathscr{L}\left[V,K\right],\\
V^{*}\otimes V\ni\boldsymbol{f}\otimes\mathbf{v}\overset{\widetilde{\Lambda}_{\mathcal{B}}}{\longmapsto}\left(\mathbf{u}\mapsto\boldsymbol{f}\!\left(\mathbf{u}\right)\!\mathbf{v}\right)\in\mathscr{L}\left[V,V\right].
\end{gather*}
In the first three cases, $\Lambda_{\mathcal{B}}=\widetilde{\Lambda}_{\mathcal{B}}$,
and the tensor maps considered are elementary tensor maps corresponding
to simple bilateral tensors in $K\otimes K$, $K\otimes V$ and $V^{*}\otimes K$,
respectively. In the fourth case, $\Lambda_{\mathcal{B}}:V^{*}\otimes V\rightarrow\mathscr{L}\left[V,V\right]$
is obtained from $\widetilde{\Lambda}_{\mathcal{B}}$ by linear extension.
Accordingly, $K\otimes K\cong\mathscr{L}\left[K,K\right]$, $K\otimes V\cong\mathscr{L}\left[K,V\right]$,
$V^{*}\otimes K\cong\mathscr{L}\left[V,K\right]$ and $V^{*}\otimes V\cong\mathscr{L}\left[V,V\right]$.
Note that the scalar $\sigma=\phi\xi=\phi\otimes\xi$ is represented
by the scalar-like tensor map $\eta\mapsto\eta\sigma$, the vector
$\mathbf{w}=\phi\mathbf{v}=\phi\otimes\mathbf{v}$ is represented
by the vector-like tensor map $\eta\mapsto\eta\mathbf{w}$, and the
linear form $\boldsymbol{g}=\xi\boldsymbol{f}=\boldsymbol{f}\otimes\xi$
is the tensor map $\mathbf{u}\mapsto\boldsymbol{g}\!\left(\mathbf{u}\right)$.

\section{Composition and contraction of tensor maps}

\subsection{Composition of tensor maps}

$ $\\
\emph{(1)}. Let $m,n,p,q,\varrho$ be non-negative integers such that
$\varrho\leq m,q$, and consider the mapping 
\begin{gather}
\widetilde{\beta}_{\circ}^{\varrho}:\Lambda_{\mathcal{B}}\left(\left.V^{*\otimes m}\otimes V^{\otimes n}\right|_{\mathcal{S}}\right)\times\Lambda_{\mathcal{B}}\left(\left.V^{*\otimes p}\otimes V^{\otimes q}\right|_{\mathcal{S}}\right)\longrightarrow\mathscr{L}\left[V^{\otimes m+p-\varrho},V^{\otimes n+q-\varrho}\right],\nonumber \\
\left(\left(\bigotimes_{i=1}^{m}\mathbf{u}_{i}\mapsto\prod_{i=1}^{m}\boldsymbol{f}^{i}\!\left(\mathbf{u}_{i}\right)\bigotimes_{i=1}^{n}\mathbf{v}_{i}\right),\left(\bigotimes_{i=1}^{p}\bar{\mathbf{u}}_{i}\mapsto\prod_{i=1}^{p}\bar{\boldsymbol{f}}\vphantom{}^{i}\!\left(\bar{\mathbf{u}}_{i}\right)\bigotimes_{i=1}^{q}\bar{\mathbf{v}}_{i}\right)\right)\longmapsto\label{eq:multiso}\\
\left(\bigotimes_{i=\varrho+1}^{m}\!\!\mathbf{u}_{i}\otimes\,\bigotimes_{i=1}^{p}\bar{\mathbf{u}}_{i}\mapsto\prod_{i=1}^{\varrho}\boldsymbol{f}^{i}\!\left(\bar{\mathbf{v}}_{i}\right)\!\prod_{i=\varrho+1}^{m}\!\!\boldsymbol{f}^{i}\!\left(\mathbf{u}_{i}\right)\prod_{i=1}^{p}\bar{\boldsymbol{f}}\vphantom{}^{i}\!\left(\bar{\mathbf{u}}_{i}\right)\bigotimes_{i=1}^{n}\mathbf{v}_{i}\otimes\bigotimes_{i=\varrho+1}^{q}\!\!\bar{\mathbf{v}}_{i}\right).\nonumber 
\end{gather}
where $\prod_{i=1}^{0}\boldsymbol{f}^{i}\!\left(\bar{\mathbf{v}}_{i}\right)=1$,
$\prod_{i=1}^{0}\boldsymbol{f}^{i}\!\left(\mathbf{u}_{i}\right)=\phi\eta$,
$\prod_{i=\varrho+1}^{\varrho>0}\boldsymbol{f^{i}}\!\left(\mathbf{u}_{i}\right)=1$,
$\prod_{i=1}^{0}\bar{\boldsymbol{f}}\vphantom{}^{i}\!\left(\bar{\mathbf{u}}_{i}\right)=\bar{\phi}\bar{\eta}$,
$\bigotimes_{i=\varrho+1}^{\varrho>0}\mathbf{u}_{i}=1$ and $\bigotimes_{i=\varrho+1}^{\varrho>0}\overline{\mathbf{v}}_{i}=1,$
and where $\phi,\eta,\bar{\phi},\bar{\eta}\in K$ are scalar variables.
(Recall that $\bigotimes_{i=1}^{0}\mathbf{u}_{i}=\eta$, $\bigotimes_{i=1}^{0}\bar{\mathbf{u}}_{i}=\bar{\eta},$
$\bigotimes_{i=1}^{0}\mathbf{v}_{i}=\xi$ and $\bigotimes_{i=1}^{0}\overline{\mathbf{v}}_{i}=\overline{\xi}$.)
\begin{example*}
Set $\widetilde{\boldsymbol{s}}=\left(\mathbf{u}\mapsto\boldsymbol{f}\!\left(\mathbf{u}\right)\xi=\boldsymbol{g}\!\left(\mathbf{u}\right)\right)$
and $\widetilde{\boldsymbol{t}}=\left(\bar{\eta}\mapsto\bar{\phi}\bar{\eta}\,\bar{\mathbf{v}}=\bar{\eta}\,\bar{\mathbf{w}}\right)$.
Then 
\begin{gather*}
\widetilde{\beta}_{\circ}^{0}\left(\widetilde{\boldsymbol{s}},\widetilde{\boldsymbol{t}}\right)=\left(\mathbf{u}\otimes\bar{\eta}\mapsto1\boldsymbol{\boldsymbol{f}}\!\left(\mathbf{u}\right)\bar{\phi}\bar{\eta}\,\xi\otimes\bar{\mathbf{v}}\right)=\left(\bar{\eta}\mathbf{u}\mapsto\bar{\eta}\xi\boldsymbol{\boldsymbol{f}}\!\left(\mathbf{u}\right)\bar{\phi}\bar{\mathbf{v}}\right)=\left(\mathbf{u}\mapsto\boldsymbol{g}\!\left(\mathbf{u}\right)\bar{\mathbf{w}}\right),\\
\widetilde{\beta}_{\circ}^{1}\left(\widetilde{\boldsymbol{s}},\widetilde{\boldsymbol{t}}\right)=\left(1\otimes\bar{\eta}\mapsto\boldsymbol{\boldsymbol{f}}\!\left(\bar{\mathbf{v}}\right)1\bar{\phi}\bar{\eta}\,\xi\otimes1\right)=\left(\bar{\eta}\mapsto\bar{\eta}\xi\boldsymbol{\boldsymbol{f}}\!\left(\bar{\phi}\bar{\mathbf{v}}\right)\right)=\left(\left(\bar{\eta}\mapsto\bar{\eta}\boldsymbol{\boldsymbol{g}}\!\left(\bar{\mathbf{w}}\right)\right)\right).\quad\square
\end{gather*}

\end{example*}
We can use $\widetilde{\beta}_{\circ}^{\varrho}$ and the isomorphism
$\Lambda_{\mathcal{B}}:V^{*m}\otimes V^{n}\rightarrow\mathscr{L}\left[V^{\otimes m},V^{\otimes n}\right]$
to define a bilinear map
\[
\beta_{\circ}^{\varrho}:\mathscr{L}\left[V^{\otimes m},V^{\otimes n}\right]\times\mathscr{L}\left[V^{\otimes p},V^{\otimes q}\right]\;\rightarrow\;\mathscr{L}\left[V^{\otimes m+p-\varrho},V^{\otimes n+q-\varrho}\right]
\]
by a 'basis-free' argument similar to that used to define multiplication
of bilateral tensors in Section 6.4. Alternatively, it suffices to
note that since (in the finite-dimensional case) we can choose a basis
$\left\{ \mathbf{e}_{i}\right\} $ for $V$ and a dual basis $\left\{ \boldsymbol{e}^{i}\right\} $
for $V^{*}$, $\mathscr{L}\left[V^{\otimes m},V^{\otimes n}\right]$
has a basis of the form $\left\{ \bigotimes_{i=1}^{m}\mathbf{u}_{i}\mapsto\prod_{i=1}^{m}\boldsymbol{f}^{i}\left(\mathbf{u}_{i}\right)\bigotimes_{i=1}^{n}\mathbf{v}_{i}\right\} $,
$\mathscr{L}\left[V^{\otimes p},V^{\otimes q}\right]$ has a basis
of the form $\left\{ \bigotimes_{i=1}^{p}\bar{\mathbf{u}}_{i}\mapsto\prod_{i=1}^{p}\bar{\boldsymbol{f}}\vphantom{}^{i}\left(\bar{\mathbf{u}}_{i}\right)\bigotimes_{i=1}^{q}\bar{\mathbf{v}}_{i}\right\} $,
and $\mathscr{L}\left[V^{\otimes m+p-\varrho},V^{\otimes n+q-\varrho}\right]$
has a basis of the form 
\[
\left\{ \widetilde{\beta}_{\circ}^{\varrho}\left(\left(\bigotimes_{i=1}^{m}\mathbf{u}_{i}\mapsto\prod_{i=1}^{m}\boldsymbol{f}^{i}\left(\mathbf{u}_{i}\right)\bigotimes_{i=1}^{n}\mathbf{v}_{i}\right),\left(\bigotimes_{i=1}^{p}\bar{\mathbf{u}}_{i}\mapsto\prod_{i=1}^{p}\bar{\boldsymbol{f}}\vphantom{}^{i}\left(\bar{\mathbf{u}}_{i}\right)\bigotimes_{i=1}^{q}\bar{\mathbf{v}}_{i}\right)\right)\right\} .
\]
Thus, $\beta_{\circ}^{\varrho}$ can be obtained by bilinear extension
of $\widetilde{\beta}_{\circ}^{\varrho}$, and $\beta_{\circ}^{\varrho}$
is a surjective map.

\emph{(2)}. In double-index notation, we denote $\beta_{\circ}^{\varrho}\left(\boldsymbol{s},\boldsymbol{t}\right)$
by 
\[
\circ{}_{\mathrm{a}_{1}\ldots\mathrm{a}_{\varrho}}^{\mathrm{d}_{1}\ldots\mathrm{d}_{\varrho}}\left(s{}_{\mathrm{a}_{1}\ldots\mathrm{a}_{m}}^{\mathrm{b}_{1}\ldots\mathrm{b}_{n}},t{}_{\mathrm{c}_{1}\ldots\mathrm{c}_{p}}^{\mathrm{d}_{1}\ldots\mathrm{d}_{q}}\right)\quad\mathrm{or}\quad s{}_{\mathrm{a}_{1}\ldots\mathrm{a}_{m}}^{\mathrm{b}_{1}\ldots\mathrm{b}_{n}}\circ t{}_{\mathrm{c}_{1}\ldots\mathrm{c}_{p}}^{\mathrm{a}_{1}\ldots\mathrm{a}_{\varrho}\mathrm{d}_{\varrho+1}\ldots\mathrm{d}_{q}}\:,
\]
where we have used matching indices\emph{ }in the last case. It is
clear that we can introduce a general bilinear function
\[
\left(s_{\mathrm{a}_{1}\ldots\mathrm{a}_{m}}^{\mathrm{b}_{1}\ldots\mathrm{b}_{n}},t_{\mathrm{c}_{1}\ldots\mathrm{c}_{p}}^{\mathrm{d}_{1}\ldots\mathrm{d}_{q}}\right)\mapsto\circ{}_{\mathrm{a}_{i_{1}}\ldots\mathrm{a}_{i_{\varrho}}}^{\mathrm{d}_{j_{1}}\ldots\mathrm{d}_{j_{\varrho}}}\left(s_{\mathrm{a}_{1}\ldots\mathrm{a}_{m}}^{\mathrm{b}_{1}\ldots\mathrm{b}_{n}},t_{\mathrm{c}_{1}\ldots\mathrm{c}_{p}}^{\mathrm{d}_{1}\ldots\mathrm{d}_{q}}\right)
\]
in the same way as we did for multiplication of bilateral tensors.
We have thus defined \emph{general composition of linear tensor maps}. 

Using matching indices to specify the composition, we have, for example,
\[
\circ\!\left(s_{\mathrm{a}}^{\mathrm{b}},t_{\mathrm{c}}^{\mathrm{d}}\right)=s_{\mathrm{a}}^{\mathrm{b}}\circ t_{\mathrm{c}}^{\mathrm{d}}=p_{\mathrm{ac}}^{\mathrm{bd}},\quad\circ_{\mathrm{a}}^{\mathrm{d}}\!\left(s_{\mathrm{a}}^{\mathrm{b}},t_{\mathrm{c}}^{\mathrm{d}}\right)=s_{\mathrm{a}}^{\mathrm{b}}\circ t_{\mathrm{c}}^{\mathrm{a}}=q_{\mathrm{c}}^{\mathrm{b}},\quad\circ{}_{\mathrm{b}}^{\mathrm{d}}\!\left(s_{\mathrm{a}\mathrm{b}}^{\mathrm{c}},t^{\mathrm{d\mathrm{e}}}\right)=s_{\mathrm{a}\mathrm{b}}^{\mathrm{c}}\circ t^{\mathrm{b\mathrm{e}}}=r_{\mathrm{a}}^{\mathrm{ce}}.
\]

\subsection{Outer and inner composition of tensor maps}

$ $\\
\emph{(1)}. When $\varrho=0$ we have composition without matching
indices, where
\[
\circ\left(s_{\mathrm{a}_{1}\ldots\mathrm{a}_{m}}^{\mathrm{b}_{1}\ldots\mathrm{b}_{n}},t_{\mathrm{c}_{1}\ldots\mathrm{c}_{p}}^{\mathrm{d}_{1}\ldots\mathrm{d}_{q}}\right)=s_{\mathrm{a}_{1}\ldots\mathrm{a}_{m}}^{\mathrm{b}_{1}\ldots\mathrm{b}_{n}}\circ t_{\mathrm{c}_{1}\ldots\mathrm{c}_{p}}^{\mathrm{d}_{1}\ldots\mathrm{d}_{q}}=p_{\mathrm{a}_{1}\ldots\mathrm{a}_{m}\mathrm{c}_{1}\ldots\mathrm{c}_{p}}^{\mathrm{b}_{1}\ldots\mathrm{b}_{n}\mathrm{d}_{1}\ldots\mathrm{d}_{q}}
\]
is defined by 
\begin{gather}
\left(\left(\bigotimes_{i=1}^{m}\mathbf{u}_{i}\mapsto\prod_{i=1}^{m}\boldsymbol{f}^{i}\left(\mathbf{u}_{i}\right)\bigotimes_{i=1}^{n}\mathbf{v}_{i}\right),\left(\bigotimes_{i=1}^{p}\bar{\mathbf{u}}_{i}\mapsto\prod_{i=1}^{p}\bar{\boldsymbol{f}}\vphantom{}^{i}\left(\bar{\mathbf{u}}_{i}\right)\bigotimes_{i=1}^{q}\bar{\mathbf{v}}_{i}\right)\right)\longmapsto\nonumber \\
\bigotimes_{i=1}^{m}\mathbf{u}_{i}\otimes\bigotimes_{i=1}^{p}\bar{\mathbf{u}}{}_{i}\;\mapsto\;\prod_{i=1}^{m}\boldsymbol{f}^{i}\left(\mathbf{u}_{i}\right)\prod_{i=1}^{p}\bar{\boldsymbol{f}}\vphantom{}^{i}\left(\bar{\mathbf{u}}{}_{i}\right)\bigotimes_{i=1}^{n}\mathbf{v}_{i}\otimes\bigotimes_{i=1}^{q}\bar{\mathbf{v}}_{i},\label{eq:outercomp}
\end{gather}
where $ $$\prod_{i=1}^{0}\boldsymbol{f}^{i}\left(\mathbf{u}_{i}\right)=\phi\eta$
and $\prod_{i=1}^{0}\bar{\boldsymbol{f}}\vphantom{}^{i}\left(\bar{\mathbf{u}}{}_{i}\right)=\bar{\phi}\bar{\eta}$,
$ $according to (\ref{eq:multiso}). We call this tensor map the
\emph{outer composition} of $\boldsymbol{s}$ and $\boldsymbol{t}$
and denote it in index-free form by 
\[
\boldsymbol{s}\otimes\boldsymbol{t}.
\]

\emph{(2)}. If, on the other hand, $\varrho=m=q$, then 
\[
\circ{}_{\mathrm{a}_{1}\ldots\mathrm{a}_{\varrho}}^{\mathrm{d}_{1}\ldots\mathrm{d}_{\varrho}}\left(s_{\mathrm{a}_{1}\ldots\mathrm{a}_{\varrho}}^{\mathrm{b}_{1}\ldots\mathrm{b}_{n}},t_{\mathrm{c}_{1}\ldots\mathrm{c}_{p}}^{\mathrm{d}_{1}\ldots\mathrm{d}_{\varrho}}\right)=s_{\mathrm{a}_{1}\ldots\mathrm{a}_{\varrho}}^{\mathrm{b}_{1}\ldots\mathrm{b}_{n}}\circ t_{\mathrm{c}_{1}\ldots\mathrm{c}_{p}}^{\mathrm{a}_{1}\ldots\mathrm{a}_{\varrho}}=p_{\mathrm{c}_{1}\ldots\mathrm{c}_{p}}^{\mathrm{b}_{1}\ldots\mathrm{b}_{n}}
\]
 is defined by
\begin{gather}
\left(\left(\bigotimes_{i=1}^{m}\mathbf{u}_{i}\mapsto\prod_{i=1}^{m}\boldsymbol{f}^{i}\left(\mathbf{u}_{i}\right)\bigotimes_{i=1}^{n}\mathbf{v}_{i}\right),\left(\bigotimes_{i=1}^{p}\overline{\mathbf{u}}_{i}\mapsto\prod_{i=1}^{p}\bar{\boldsymbol{f}}\vphantom{}^{i}\left(\bar{\mathbf{u}}_{i}\right)\bigotimes_{i=1}^{q}\bar{\mathbf{v}}_{i}\right)\right)\longmapsto\nonumber \\
\bigotimes_{i=1}^{p}\bar{\mathbf{u}}{}_{i}\;\mapsto\;\prod_{i=1}^{\varrho}\boldsymbol{f}^{i}\left(\bar{\mathbf{v}}_{i}\right)\prod_{i=1}^{p}\bar{\boldsymbol{f}}\vphantom{}^{i}\left(\bar{\mathbf{u}}{}_{i}\right)\bigotimes_{i=1}^{n}\mathbf{v}_{i},\label{eq:incomp}
\end{gather}
where $\prod_{i=1}^{0}\boldsymbol{f}^{i}\left(\bar{\mathbf{v}}_{i}\right)=\phi\bar{\xi}$
and $\prod_{i=1}^{0}\bar{\boldsymbol{f}}\vphantom{}^{i}\left(\bar{\mathbf{u}}{}_{i}\right)=\bar{\phi}\bar{\eta}$.
For $\varrho>0$ this follows immediately from (\ref{eq:multiso}),
and for $\varrho=0$ we can use (\ref{eq:multiso}) or (\ref{eq:outercomp})
to obtain
\begin{gather*}
\left(\left(\bigotimes_{i=1}^{m}\mathbf{u}_{i}\mapsto\prod_{i=1}^{m}\boldsymbol{f}^{i}\left(\mathbf{u}_{i}\right)\bigotimes_{i=1}^{n}\mathbf{v}_{i}\right),\left(\bigotimes_{i=1}^{p}\bar{\mathbf{u}}_{i}\mapsto\prod_{i=1}^{p}\bar{\boldsymbol{f}}\vphantom{}^{i}\left(\bar{\mathbf{u}}_{i}\right)\bigotimes_{i=1}^{q}\bar{\mathbf{v}}_{i}\right)\right)\longmapsto\\
\eta\otimes\bigotimes_{i=1}^{p}\bar{\mathbf{u}}_{i}\;\mapsto\;\phi\eta\prod_{i=1}^{p}\bar{\boldsymbol{f}}\vphantom{}^{i}\left(\bar{\mathbf{u}}{}_{i}\right)\left(\left.\bigotimes\right._{i=1}^{n}\mathbf{v}_{i}\right)\otimes\bar{\xi},
\end{gather*}
which simplifies to (\ref{eq:incomp}). We call the tensor map in
the case $\varrho=m=q$ the \emph{inner composition} of $\boldsymbol{s}$
and $\boldsymbol{t}$ and denote it in index-free form by 
\[
\boldsymbol{s}\circ\boldsymbol{t}.
\]

Note that if $\varrho=m=q=0$ then the composition of $\boldsymbol{s}$
and $\boldsymbol{t}$ according to (\ref{eq:multiso}) is both an
outer composition and an inner composition, so $\boldsymbol{s}\otimes\boldsymbol{t}=\boldsymbol{s}\circ\boldsymbol{t}$.

\subsection{Contraction of tensor maps}

$ $\\
Set $\widetilde{\boldsymbol{t}}=\left(\bigotimes_{i=1}^{m}\mathbf{u}_{i}\mapsto\prod_{i=1}^{m}\boldsymbol{f}^{i}\!\left(\mathbf{u}_{i}\right)\bigotimes_{j=1}^{n}\mathbf{v}_{j}\right)$,
$\widetilde{\boldsymbol{p}}=\left(\bigotimes_{i=1}^{m}\mathbf{u}_{i}\mapsto\prod_{i=1}^{m}\boldsymbol{f}^{i}\!\left(\mathbf{u}_{i}\right)\right)$
and $\widetilde{\boldsymbol{q}}=\left(\eta\mapsto\eta\bigotimes_{j=1}^{n}\mathbf{v}_{j}\right)$.
Then $\widetilde{\boldsymbol{t}}=\widetilde{\boldsymbol{p}}\otimes\widetilde{\boldsymbol{q}}$,
and it is clear that any tensor map $\widetilde{\boldsymbol{t}}$
of the form shown can be written in a unique way as the outer product
$\widetilde{\boldsymbol{p}}\otimes\widetilde{\boldsymbol{q}}$ of
two tensor maps of the forms shown. In double-index notation, we have
\[
\widetilde{t}_{\mathrm{a}_{1}\ldots\mathrm{a}_{m}}^{\mathrm{b}_{1}\ldots\mathrm{b}_{n}}=\circ\left(\widetilde{p}_{\mathrm{a}_{1}\ldots\mathrm{a}_{m}},\widetilde{q}^{\mathrm{b}_{1}\ldots\mathrm{b}_{n}}\right)=\widetilde{p}_{\mathrm{a}_{1}\ldots\mathrm{a}_{m}}\circ\widetilde{q}^{\mathrm{b}_{1}\ldots\mathrm{b}_{n}}.
\]

We can use a mapping 
\[
\widetilde{t}_{\mathrm{a}_{1}\ldots\mathrm{a}_{m}}^{\mathrm{b}_{1}\ldots\mathrm{b}_{n}}=\circ\left(\widetilde{p}_{\mathrm{a}_{1}\ldots\mathrm{a}_{m}},\widetilde{q}^{\mathrm{b}_{1}\ldots\mathrm{b}_{n}}\right)\mapsto\circ{}_{\mathrm{a}_{i_{1}}\ldots\mathrm{a}_{i_{\varrho}}}^{\mathrm{b}_{j_{1}}\ldots\mathrm{b}_{j_{\varrho}}}\left(\widetilde{p}_{\mathrm{a}_{1}\ldots\mathrm{a}_{m}},\widetilde{q}^{\mathrm{b}_{1}\ldots\mathrm{b}_{n}}\right),
\]
to define a linear map
\begin{gather*}
\Diamond{}_{\mathrm{a}_{i_{1}}\ldots\mathrm{a}_{i_{\varrho}}}^{\mathrm{b}_{j_{1}}\ldots\mathrm{b}_{j_{\varrho}}}:\mathscr{L}\left[V^{\otimes m},V^{\otimes n}\right]\rightarrow\mathscr{L}\left[V^{\otimes\left(m-\varrho\right)},V^{\otimes\left(n-\varrho\right)}\right],\\
t_{\mathrm{a}_{1}\ldots\mathrm{a}_{m}}^{\mathrm{b}_{1}\ldots\mathrm{b}_{n}}\mapsto\Diamond{}_{\mathrm{a}_{i_{1}}\ldots\mathrm{a}_{i_{\varrho}}}^{\mathrm{b}_{j_{1}}\ldots\mathrm{b}_{j_{\varrho}}}\left(t_{\mathrm{a}_{1}\ldots\mathrm{a}_{m}}^{\mathrm{b}_{1}\ldots\mathrm{b}_{n}}\right)=\circ{}_{\mathrm{a}_{i_{1}}\ldots\mathrm{a}_{i_{\varrho}}}^{\mathrm{b}_{j_{1}}\ldots\mathrm{b}_{j_{\varrho}}}\left(p_{\mathrm{a}_{1}\ldots\mathrm{a}_{m}},q^{\mathrm{b}_{1}\ldots\mathrm{b}_{n}}\right)
\end{gather*}
in the same way as we did for bilateral tensors. $\Diamond{}_{\mathrm{a}_{i_{1}}\ldots\mathrm{a}_{i_{\varrho}}}^{\mathrm{b}_{j_{1}}\ldots\mathrm{b}_{j_{\varrho}}}\left(t_{\mathrm{a}_{1}\ldots\mathrm{a}_{m}}^{\mathrm{b}_{1}\ldots\mathrm{b}_{n}}\right)$
is the \emph{contraction} of $t_{\mathrm{a}_{1}\ldots\mathrm{a}_{m}}^{\mathrm{b}_{1}\ldots\mathrm{b}_{n}}$
with respect to $\mathrm{\mathsf{a}}_{j_{1}}\ldots\mathrm{\mathsf{a}}_{j_{\varrho}}$
and $\mathrm{\mathsf{b}}_{j_{1}}\ldots\mathrm{\mathsf{b}}_{j_{\varrho}}$. 

We can again use matching indices to specify contractions. For example,
with $x_{\mathrm{a}}^{\mathrm{b}}=p_{\mathrm{a}}\circ q^{\mathrm{b}}$
and $y_{\mathrm{ab}}^{\mathrm{cd}}=r{}_{\mathrm{ab}}\circ s{}^{\mathrm{cd}}$
we can form contractions 
\begin{gather*}
\lozenge_{\mathrm{a}}^{\mathrm{b}}\left(x_{\mathrm{a}}^{\mathrm{b}}\right)=\circ{}_{\mathrm{a}}^{\mathrm{b}}\left(p{}_{\mathrm{a}},q{}^{\mathrm{b}}\right)=p{}_{\mathrm{a}}\circ q{}^{\mathrm{a}}=x{}_{\mathrm{a}}^{\mathrm{a}}=\hat{x},\\
\Diamond_{\mathrm{ab}}^{\mathrm{cd}}\left(y_{\mathrm{ab}}^{\mathrm{cd}}\right)=\circ{}_{\mathrm{ab}}^{\mathrm{cd}}\left(r{}_{\mathrm{ab}},s{}^{\mathrm{cd}}\right)=r{}_{\mathrm{ab}}\circ s{}^{\mathrm{ab}}=y{}_{\mathrm{ab}}^{\mathrm{ab}}=\hat{y},\\
\Diamond_{\mathrm{a}}^{\mathrm{d}}\left(y_{\mathrm{ab}}^{\mathrm{cd}}\right)=\circ_{\mathrm{a}}^{\mathrm{d}}\left(r{}_{\mathrm{ab}},s{}^{\mathrm{cd}}\right)=r{}_{\mathrm{ab}}\circ s{}^{\mathrm{ca}}=y{}_{\mathrm{ab}}^{\mathrm{ca}}=\hat{y}_{\mathrm{b}}^{\mathrm{c}}.
\end{gather*}

\subsection{Multiplication of bilateral tensors and composition of tensor maps}

Recall the isomorphism $\Lambda_{\mathcal{B}}:V^{*m}\otimes V^{n}\rightarrow\mathscr{L}\left[V^{\otimes m},V^{\otimes n}\right]$.
It can be verified by straightforward calculation that
\[
\beta_{\circ}^{\varrho}\left(\Lambda_{\mathcal{B}}\left(\boldsymbol{\mathsf{s}}\right),\Lambda_{\mathcal{B}}\left(\boldsymbol{\mathsf{t}}\right)\right)=\Lambda_{\mathcal{B}}\left(\beta_{\bullet}^{\varrho}\left(\boldsymbol{\mathsf{s}},\boldsymbol{\mathsf{t}}\right)\right).
\]
Using matching indices $\mathsf{a}_{1},\ldots,\mathsf{a}_{\varrho}$,
we can write
\[
s_{\mathrm{a}_{1}\ldots\mathrm{a}_{m}}^{\mathrm{b}_{1}\ldots\mathrm{b}_{n}}\circ t{}_{\mathrm{c}_{1}\ldots\mathrm{c}_{p}}^{\mathrm{a}_{1}\ldots\mathrm{a}_{\varrho}\mathrm{d}_{\varrho+1}\ldots\mathrm{d}_{q}}=\Lambda_{\mathcal{B}}\left(\mathsf{s}_{\mathrm{a}_{1}\ldots\mathrm{a}_{m}}^{\mathrm{b}_{1}\ldots\mathrm{b}_{n}}\mathsf{t}{}_{\mathrm{c}_{1}\ldots\mathrm{c}_{p}}^{\mathrm{a}_{1}\ldots\mathrm{a}_{\varrho}\mathrm{d}_{\varrho+1}\ldots\mathrm{d}_{q}}\right).
\]
In the general case, where we make no assumptions about where matching
indices occur, we have
\begin{equation}
\circ{}_{\mathrm{a}_{i_{1}}\ldots\mathrm{a}_{i_{\varrho}}}^{\mathrm{d}_{j_{1}}\ldots\mathrm{d}_{j_{\varrho}}}\left(\Lambda_{\mathcal{B}}\left(\mathsf{s}_{\mathrm{a}_{1}\ldots\mathrm{a}_{m}}^{\mathrm{b}_{1}\ldots\mathrm{b}_{n}}\right),\Lambda_{\mathcal{B}}\left(\mathsf{t}{}_{\mathrm{c}_{1}\ldots\mathrm{c}_{p}}^{\mathrm{a}_{1}\ldots\mathrm{d}_{q}}\right)\right)=\Lambda_{\mathcal{B}}\left(\bullet{}_{\mathrm{a}_{i_{1}}\ldots\mathrm{a}_{i_{\varrho}}}^{\mathrm{d}_{j_{1}}\ldots\mathrm{d}_{j_{\varrho}}}\left(\mathsf{s}_{\mathrm{a}_{1}\ldots\mathrm{a}_{m}}^{\mathrm{b}_{1}\ldots\mathrm{b}_{n}},\mathsf{t}_{\mathrm{c}_{1}\ldots\mathrm{d}_{p}}^{\mathrm{d}_{1}\ldots\mathrm{d}_{q}}\right)\right).\label{eq:isobilmapmult}
\end{equation}

As $V^{*m}\otimes V^{n}$ and $\mathscr{L}\left[V^{\otimes m},V^{\otimes n}\right]$
are canonically isomorphic as vector spaces, and multiplication of
tensors in $V^{*m}\otimes V^{n}$ is compatible with composition of
tensor maps in $\mathscr{L}\left[V^{\otimes m},V^{\otimes n}\right]$
according to (\ref{eq:isobilmapmult}), we conclude that for finite-dimensional
vector spaces the interpretation of classical tensors as tensor maps
is equivalent to the traditional interpretation of classical tensors
as bilateral tensors.

Note that the fact that the multiplication operators $\bullet$ and
$\circ$ are compatible implies that the contraction operator $\blacklozenge$
for bilateral tensors is compatible with the contraction operator
$\Diamond$ for tensor maps,

Finally, since general multiplication of bilateral tensors is associative,
general composition of tensor maps is also associative.

\section{Interpretations of tensor map composition}

\subsection{Composition of tensor maps as generalized function composition}

Composition of tensor maps turns out to be related not only to multiplication
of bilateral tensors but also to ordinary function composition.

Let us first look at scalar-like tensor maps. If $s=\left(\eta\mapsto\phi\eta\xi\right)=\left(\eta\mapsto\eta\sigma\right)$,
$t=\left(\bar{\eta}\mapsto\bar{\phi}\bar{\eta}\bar{\xi}\right)=\left(\bar{\eta}\mapsto\bar{\eta}\tau\right)$
and $\circ$ denotes ordinary function composition then\linebreak{}
$s\circ t\!\left(\bar{\eta}\right)=s\!\left(t\!\left(\bar{\eta}\right)\right)=s\!\left(\bar{\eta}\tau\right)=\bar{\eta}\tau\sigma$,
so $s\circ t=\left(\bar{\eta}\rightarrow\bar{\eta}\tau\sigma\right)$.
On the other hand, by (\ref{eq:incomp}) inner composition of tensor
maps yields 
\[
s\circ t=\left(\bar{\eta}\mapsto\phi\bar{\xi}\,\bar{\phi}\bar{\eta}\xi\right)=\left(\bar{\eta}\mapsto\bar{\eta}\bar{\phi}\,\bar{\xi}\phi\xi\right)=\left(\bar{\eta}\mapsto\bar{\eta}\tau\sigma\right),
\]
so ordinary composition of functions and inner composition of tensor
maps coincide.

For another example, consider the tensor maps $\boldsymbol{f}=\left(\mathbf{u}\mapsto\boldsymbol{g}\!\left(\mathbf{u}\right)\xi\right)=\left(\mathbf{u}\mapsto\boldsymbol{f}\!\left(\mathbf{u}\right)\right)$
and $\boldsymbol{v}=\left(\bar{\eta}\mapsto\bar{\phi}\bar{\eta}\,\bar{\mathbf{u}}\right)=\left(\bar{\eta}\mapsto\bar{\eta}\,\bar{\mathbf{v}}\right)$.
In terms of ordinary function composition, we have $\boldsymbol{f}\circ\boldsymbol{v}\!\left(\bar{\eta}\right)=\boldsymbol{f}\!\left(\boldsymbol{v}\!\left(\bar{\eta}\right)\right)=\boldsymbol{f}\!\left(\bar{\eta}\,\bar{\mathbf{v}}\right)=\bar{\eta}\boldsymbol{f}\!\left(\bar{\mathbf{v}}\right)$,
so $\boldsymbol{f}\circ\boldsymbol{v}=\left(\bar{\eta}\mapsto\bar{\eta}\boldsymbol{f}\!\left(\bar{\mathbf{v}}\right)\right)$.
Using (\ref{eq:incomp}), we obtain
\[
\boldsymbol{f}\circ\boldsymbol{v}=f_{\mathrm{a}}\circ v^{\mathrm{a}}=\left(\bar{\eta}\mapsto\boldsymbol{g}\!\left(\bar{\mathbf{u}}\right)\bar{\phi}\bar{\eta}\xi\right)=\left(\bar{\eta}\mapsto\bar{\eta}\xi\,\boldsymbol{g}\!\left(\bar{\phi}\bar{\mathbf{u}}\right)\right)=\left(\bar{\eta}\mapsto\bar{\eta}\boldsymbol{f}\!\left(\bar{\mathbf{v}}\right)\right),
\]
so inner composition of tensor maps coincides with ordinary function
composition again.

We now turn to the general case. Let $\left\{ \mathbf{e}_{i}\right\} $
be a basis for $V$ and $\left\{ \boldsymbol{e}^{i}\right\} $ the
dual basis for $V^{*}$. Then $\left(\left.\bigotimes\right._{\ell=1}^{p}\mathbf{u}_{\ell}\mapsto\left.\prod\right._{\ell=1}^{p}\boldsymbol{e}^{i_{\ell}}\!\left(\mathbf{u}_{\ell}\right)\left.\bigotimes\right._{\ell=1}^{n}\mathbf{e}_{j_{\ell}}\right)=\left\langle \left.\bigotimes\right._{\ell=1}^{p}\mathbf{e}_{i_{\ell}}\mapsto\left.\bigotimes\right._{\ell=1}^{n}\mathbf{e}_{j_{\ell}}\right\rangle .$

With $\circ$ denoting the inner product of tensor maps, (\ref{eq:incomp})
thus gives

\begin{center}
$\left\langle \bigotimes_{\ell=1}^{m}\mathbf{e}_{i_{\ell}}\mapsto\bigotimes_{\ell=1}^{n}\mathbf{e}_{j_{\ell}}\right\rangle \circ\left\langle \bigotimes_{\ell=1}^{p}\mathbf{e}_{k_{\ell}}\mapsto\bigotimes_{\ell=1}^{m}\mathbf{e}_{i_{\ell}}\right\rangle $
\par\end{center}

\begin{center}
$ $$=\left(\bigotimes_{\ell=1}^{m}\mathbf{u}_{\ell}\mapsto\prod_{\ell=1}^{m}\boldsymbol{e}^{i_{\ell}}\!\left(\mathbf{u}_{\ell}\right)\bigotimes_{\ell=1}^{n}\mathbf{e}_{j_{\ell}}\right)\circ\left(\bigotimes_{\ell=1}^{p}\bar{\mathbf{u}}_{\ell}\mapsto\prod_{\ell=1}^{p}\boldsymbol{e}^{k_{\ell}}\!\left(\bar{\mathbf{u}}_{\ell}\right)\bigotimes_{\ell=1}^{m}\mathbf{e}_{i_{\ell}}\right)$
\par\end{center}

\begin{center}
$=\left(\bigotimes_{\ell=1}^{p}\bar{\mathbf{u}}_{\ell}\mapsto\prod_{\ell=1}^{m}\boldsymbol{e}^{i_{\ell}}\!\left(\mathbf{e}_{i_{\ell}}\right)\prod_{\ell=1}^{p}\boldsymbol{e}^{k_{\ell}}\!\left(\bar{\mathbf{u}}_{\ell}\right)\bigotimes_{\ell=1}^{n}\mathbf{e}_{j_{\ell}}\right)$
\par\end{center}

\begin{center}
$=\left(\bigotimes_{\ell=1}^{p}\bar{\mathbf{u}}_{\ell}\mapsto\prod_{\ell=1}^{p}\boldsymbol{e}^{k_{\ell}}\!\left(\bar{\mathbf{u}}_{\ell}\right)\bigotimes_{\ell=1}^{n}\mathbf{e}_{j_{\ell}}\right)=\left\langle \bigotimes_{\ell=1}^{p}\mathbf{e}_{k_{\ell}}\mapsto\bigotimes_{\ell=1}^{n}\mathbf{e}_{j_{\ell}}\right\rangle .$
\par\end{center}

On the other hand, with $\circ$ denoting ordinary function composition
we have

\begin{center}
$\left\langle \bigotimes_{\ell=1}^{m}\mathbf{e}_{i_{\ell}}\mapsto\bigotimes_{\ell=1}^{n}\mathbf{e}_{j_{\ell}}\right\rangle \circ\left\langle \bigotimes_{\ell=1}^{p}\mathbf{e}_{k_{\ell}}\mapsto\bigotimes_{\ell=1}^{m}\mathbf{e}_{i_{\ell}}\right\rangle =\left\langle \bigotimes_{\ell=1}^{p}\mathbf{e}_{k_{\ell}}\mapsto\bigotimes_{\ell=1}^{n}\mathbf{e}_{j_{\ell}}\right\rangle .$
\par\end{center}

Since $\left\{ \left\langle \bigotimes_{\ell=1}^{r}\mathbf{e}_{i_{\ell}}\mapsto\bigotimes_{\ell=1}^{s}\mathbf{e}_{j_{\ell}}\right\rangle \right\} $
is a basis for $\mathscr{L}\left[V^{\otimes r},V^{\otimes s}\right]$,
this means that inner composition of tensor maps coincides with ordinary
function composition; if $\boldsymbol{s}$ and $\boldsymbol{t}$ are
linear tensor maps then $\boldsymbol{s}\circ\boldsymbol{t}$ is the
same map regardless of whether $\circ$ is interpreted as inner composition
of tensor maps or as ordinary function composition.

\subsection{Composition of tensor maps as generalized function application}

$ $\\
Set $s=\left(\eta\mapsto\eta\sigma\right)$, $t=\left(\eta\mapsto\eta\tau\right)$,
$u^{\mathrm{a}}=\left(\eta\mapsto\eta\mathbf{u}\right)$, $v^{\mathrm{a}}=\left(\eta\mapsto\eta\mathbf{v}\right)$
and let \linebreak{}
$\mathrm{B}=\left(\mathrm{b}_{1},\ldots,\mathrm{b}_{n}\right)$ be
a multi-index. Since tensor map composition is bilinear,
\begin{gather*}
r^{\mathrm{B}}\circ\xi s=\xi r^{\mathrm{B}}\circ s,\qquad r^{\mathrm{B}}\circ\left(s+t\right)=r^{\mathrm{B}}\circ s+r^{\mathrm{B}}\circ t,\\
r_{\mathrm{a}}^{\mathrm{B}}\circ\xi u^{\mathrm{a}}=\xi r_{\mathrm{a}}^{\mathrm{B}}\circ u^{\mathrm{a}},\qquad r_{\mathrm{a}}^{\mathrm{B}}\circ\left(u^{\mathrm{a}}+v{}^{\mathrm{a}}\right)=r_{\mathrm{a}}^{\mathrm{B}}\circ u^{\mathrm{a}}+r_{\mathrm{a}}^{\mathrm{B}}\circ v{}^{\mathrm{a}},
\end{gather*}
but on the other hand,
\begin{gather*}
r^{\mathrm{B}}\!\left(\xi\sigma\right)=\xi r^{\mathrm{B}}\!\left(\sigma\right),\qquad r^{\mathrm{B}}\!\left(\sigma+\tau\right)=r^{\mathrm{B}}\!\left(\sigma\right)+r^{\mathrm{B}}\!\left(\tau\right),\\
r_{\mathrm{a}}^{\mathrm{B}}\!\left(\xi\mathbf{u}\right)=\xi r_{\mathrm{a}}^{\mathrm{B}}\!\left(\mathbf{u}\right),\qquad r_{\mathrm{a}}^{\mathrm{B}}\!\left(\mathbf{u}+\mathbf{v}\right)=r_{\mathrm{a}}^{\mathrm{B}}\!\left(\mathbf{u}\right)+r_{\mathrm{a}}^{\mathrm{B}}\!\left(\mathbf{v}\right),
\end{gather*}
since tensor maps are linear. Hence, there is an analogy between tensor
map composition and application of a linear map to a scalar argument
or a vector argument. Let us look more closely at this, using the
fact that inner composition of tensor maps is the same as usual function
composition. 

If $s=\left(\eta\mapsto\eta\sigma\right)$ so that $s$ represents
the scalar $\sigma$, then $r^{\mathrm{B}}\circ s\!\left(\eta\right)=r^{\mathrm{B}}\!\left(s\!\left(\eta\right)\right)=r^{\mathrm{B}}\!\left(\eta\sigma\right)=\eta r^{\mathrm{B}}\!\left(\sigma\right)$,
so $r^{\mathrm{B}}\circ s=\left(\eta\mapsto\eta r^{\mathrm{B}}\!\left(\sigma\right)\right)$,
so $r^{\mathrm{B}}\circ s$ represents $r^{\mathrm{B}}\!\left(\sigma\right)$. 

Similarly, if $v^{\mathrm{a}}=\left(\eta\mapsto\eta\mathbf{v}\right)$
so that $v^{\mathrm{a}}$ represents $\mathbf{v}$, then $f_{\mathrm{a}}^{\mathrm{B}}\circ v^{\mathrm{a}}\!\left(\eta\right)=f_{\mathrm{a}}^{\mathrm{B}}\!\left(v^{\mathrm{a}}\!\left(\eta\right)\right)=f_{\mathrm{a}}^{\mathrm{B}}\!\left(\eta\mathbf{v}\right)=\eta f_{\mathrm{a}}^{\mathrm{B}}\!\left(\mathbf{v}\right)$,
so $f_{\mathrm{a}}^{\mathrm{B}}\circ v^{\mathrm{a}}=\left(\eta\mapsto\eta f_{\mathrm{a}}^{\mathrm{B}}\!\left(\mathbf{v}\right)\right)$,
so $f_{\mathrm{a}}^{\mathrm{B}}\circ v^{\mathrm{a}}$ represents $f_{\mathrm{a}}^{\mathrm{B}}\!\left(\mathbf{v}\right)$. 

More generally, set $t^{\mathrm{a}_{1}\ldots\mathrm{a}_{n}}=\left(\eta\mapsto\eta\mathbf{v}_{1}\otimes\ldots\mathbf{v}_{n}\right)$
so that $t^{\mathrm{a}_{1}\ldots\mathrm{a}_{n}}$ represents \linebreak{}
$\mathbf{v}_{1}\otimes\ldots\otimes\mathbf{v}_{n}$. Then 
\[
f_{\mathrm{a}_{1}\ldots\mathrm{a}_{n}}^{\mathrm{B}}\circ t^{\mathrm{a}_{1}\ldots\mathrm{a}_{n}}\!\left(\eta\right)=f_{\mathrm{a}_{1}\ldots\mathrm{a}_{n}}^{\mathrm{B}}\!\left(\eta\mathbf{v}_{1}\otimes\ldots\otimes\mathbf{v}_{n}\right)=\eta f_{\mathrm{a}_{1}\ldots\mathrm{a}_{n}}^{\mathrm{B}}\!\left(\mathbf{v}_{1}\otimes\ldots\otimes\mathbf{v}_{n}\right),
\]
so $f_{\mathrm{a}_{1}\ldots\mathrm{a}_{n}}^{\mathrm{B}}\circ t^{\mathrm{a}_{1}\ldots\mathrm{a}_{n}}$
represents $f_{\mathrm{a}_{1}\ldots\mathrm{a}_{n}}^{\mathrm{B}}\!\left(\mathbf{v}_{1}\otimes\ldots\otimes\mathbf{v}_{n}\right)$
or $\phi{}_{\mathrm{a}_{1},\ldots,\mathrm{a}_{n}}^{\mathrm{B}}\!\left(\mathbf{v}_{1},\ldots\mathbf{,v}_{n}\right)$,
where $\phi{}_{\mathrm{a}_{1},\ldots,\mathrm{a}_{n}}^{\mathrm{B}}$
is the separately linear map equivalent to $f_{\mathrm{a}_{1}\ldots\mathrm{a}_{n}}^{\mathrm{B}}$. 

Furthermore, with $v_{i}^{\mathrm{a}_{i}}=\boldsymbol{v}_{i}=\left(\eta\mapsto\eta\mathbf{v}_{i}\right)$
we have 
\[
v_{1}^{\mathrm{a}_{1}}\circ\ldots\circ v_{n}^{\mathrm{a}_{n}}=\boldsymbol{v}_{1}\otimes\ldots\otimes\boldsymbol{v}{}_{n}=\left(\eta^{n}\mapsto\eta^{n}\mathbf{v}_{1}\otimes\ldots\otimes\mathbf{v}_{n}\right)=t^{\mathrm{a}_{1}\ldots\mathrm{a}_{n}},
\]
so we conclude that $f_{\mathrm{a}_{1}\ldots\mathrm{a}_{n}}^{\mathrm{B}}\circ v^{\mathrm{a}_{1}}\circ\ldots\circ v^{\mathrm{a}_{n}}$
represents $f_{\mathrm{a}_{1}\ldots\mathrm{a}_{n}}^{\mathrm{B}}\!\left(\mathbf{v}_{1}\otimes\ldots\otimes\mathbf{v}_{n}\right)$
or\linebreak{}
 $\phi{}_{\mathrm{a}_{1},\ldots,\mathrm{a}_{n}}^{\mathrm{B}}\!\left(\mathbf{v}_{1},\ldots\mathbf{,v}_{n}\right)$.
In index-free notation, $\boldsymbol{f}\circ\boldsymbol{v}_{1}\otimes\ldots\otimes\boldsymbol{v}_{n}$
represents\linebreak{}
 $\boldsymbol{f}\!\left(\mathbf{v}_{1}\otimes\ldots\otimes\mathbf{v}_{n}\right)$
or $\boldsymbol{\phi}\!\left(\mathbf{v}_{1},\ldots\mathbf{,v}_{n}\right)$.

We can thus interpret $r\circ s$ as \emph{$r\left(\sigma\right)$,
}a linear map applied to a scalar argument, and we can interpret $\boldsymbol{f}\circ\boldsymbol{v}_{1}\otimes\ldots\otimes\boldsymbol{v}_{n}$
as $\boldsymbol{f}\!\left(\mathbf{v}_{1}\otimes\ldots\otimes\mathbf{v}_{n}\right)$,
a linear map applied to a simple tensor, or as $\boldsymbol{\phi}\!\left(\mathbf{v}_{1},\ldots,\mathbf{v}_{n}\right)$,\emph{
}an $n$-linear map applied to $n$ vector arguments.

\section{Isomorphisms involving spaces of tensor maps: 'index gymnastics'}

\subsection{Automorphisms on spaces of tensor maps: moving indices around}

$ $\\
Recall that any two bilateral tensors in the same vector space have
the same subscripts and the same superscripts, although subscripts
and superscripts may be rearranged by permutations. This is true also
for the indices of tensor maps in the same vector space. 

For example, let $p_{\mathrm{ab}}^{\mathrm{ba}}:V^{\otimes2}\rightarrow V^{\otimes2}$
be the automorphism given by $\mathbf{u}\otimes\mathbf{v}\mapsto\mathbf{v}\otimes\mathbf{u}$.
For any tensor map $t_{\mathrm{ab}}:V^{\otimes2}\rightarrow K$ we
can define a tensor map $t_{\mathrm{ba}}:V^{\otimes2}\rightarrow K$
by setting $t_{\mathrm{ba}}=t_{\mathrm{ab}}\circ p_{\mathrm{ba}}^{\mathrm{ab}}$,
and for any $t^{\mathrm{ab}}:K\rightarrow V^{\otimes2}$ we can define
$t^{\mathrm{ba}}:K\rightarrow V^{\otimes2}$ by setting $t^{\mathrm{ba}}=p_{\mathrm{ab}}^{\mathrm{ba}}\circ t^{\mathrm{ab}}$.

In the general case, let $\mathsf{a}_{i}\mapsto\mathsf{a}_{i}'$ be
a permutation of $\mathsf{\mathrm{a}}_{1},\ldots,\mathrm{\mathsf{a}}_{m}$,
let $\mathsf{b}_{i}\mapsto\mathsf{b}_{i}'$ be a permutation of $\mathrm{\mathsf{b}}_{1},\ldots,\mathrm{\mathsf{b}}_{n}$,
and let $p_{\mathrm{a}_{1}'\ldots\mathrm{a}_{m}'}^{\mathrm{a}_{1}\ldots\mathrm{a}_{m}}:V^{\otimes m}\rightarrow V^{\otimes m}$
and $q_{\mathrm{b}_{1}\ldots\mathrm{b}_{n}}^{\mathrm{b}_{1}'\ldots\mathrm{b}_{n}'}:V^{\otimes n}\rightarrow V^{\otimes n}$
be automorphisms defined by 
\begin{gather*}
\mathbf{v}_{\mathrm{a}_{1}'}\otimes\ldots\otimes\mathbf{v}_{\mathrm{a}_{m}'}\mapsto\mathbf{v}_{\mathrm{a}_{1}}\otimes\ldots\otimes\mathbf{v}_{\mathrm{a}_{m}},\qquad\mathbf{v}_{\mathrm{b}_{1}}\otimes\ldots\otimes\mathbf{v}_{\mathrm{b}_{n}}\mapsto\mathbf{v}_{\mathrm{b}_{1}'}\otimes\ldots\otimes\mathbf{v}_{\mathrm{b}_{n}'},
\end{gather*}
respectively. Given $t_{\mathrm{a}_{1}\ldots\mathrm{a}_{m}}^{\mathrm{b}_{1}\ldots\mathrm{b}_{n}}$,
we can define a tensor map with the same subscripts and superscripts
in the same vector space by setting
\[
t_{\mathrm{a_{1}'\ldots\mathrm{a}_{m}'}}^{\mathrm{b}_{1}'\ldots\mathrm{b}_{n}'}=q_{\mathrm{b}_{1}\ldots\mathrm{b}_{n}}^{\mathrm{b}_{1}'\ldots\mathrm{b}_{n}'}\circ t_{\mathrm{a}_{1}\ldots\mathrm{a}_{m}}^{\mathrm{b}_{1}\ldots\mathrm{b}_{n}}\circ p_{\mathrm{a}_{1}'\ldots\mathrm{a}_{m}'}^{\mathrm{a}_{1}\ldots\mathrm{a}_{m}}\,.
\]

\begin{rem*}
Automorphism-based permutations of indices can be used to define symmetric
and alternating tensor maps in the usual way. Since these constructions
are well-known, they will not be discussed here.
\end{rem*}

\subsection{Dual tensor maps: shifting of indices}

$ $\\
\emph{(1)}. Let finite-dimensional vector spaces $V,W$ be given and
consider a mapping
\[
\mathfrak{D}:\mathscr{L}\left[V,W\right]\rightarrow\mathscr{L}\left[W^{*},V^{*}\right],
\]
 given by 
\[
\mathfrak{D}\!\left(\boldsymbol{t}\right)\!\left(\boldsymbol{f}\right)\!\left(\mathbf{v}\right)=\boldsymbol{f}\!\left(\boldsymbol{t}\!\left(\mathbf{v}\right)\right)\in K\quad\forall\boldsymbol{t}\in\mathscr{L}\left[V,W\right],\:\forall\boldsymbol{f}\in W^{*},\:\forall\mathbf{v}\in V.
\]
$\mathfrak{D}$ can be shown to be a canonical isomorphism, so for
any $V$ and $m,n\geq0$ we have a canonical isomorphism 
\[
\mathscr{L}\left[V^{\otimes m},V^{\otimes n}\right]\rightarrow\mathscr{L}\left[\left(V^{\otimes n}\right)^{*},\left(V^{\otimes m}\right)^{*}\right],
\]
and in view of the canonical isomorphism $i:\left(V^{*}\right)^{\otimes N}\rightarrow\left(V^{\otimes N}\right)^{*}$
given by\linebreak{}
 $i\left(\boldsymbol{\boldsymbol{f}_{1}\otimes\ldots\otimes\boldsymbol{f}_{n}}\right)\left(\mathbf{u}_{1}\otimes\ldots\otimes\mathbf{u}_{n}\right)=\boldsymbol{f}_{1}\!\left(\mathbf{u}_{1}\right)\ldots\boldsymbol{f}_{n}\!\left(\mathbf{u}_{n}\right)$
there are also canonical isomorphisms of the form
\[
\mathscr{L}\left[V^{\otimes m},V^{\otimes n}\right]\rightarrow\mathscr{L}\left[V^{*\otimes n},V^{*\otimes m}\right],\qquad t_{\mathrm{a}_{1}\ldots\mathrm{a}_{m}}^{\mathrm{b_{1}\ldots\mathrm{b}_{n}}}\mapsto t_{\mathrm{b_{1}^{*}\ldots\mathrm{b}_{n}^{*}}}^{\mathrm{a}_{1}^{*}\ldots\mathrm{a}_{m}^{*}}.
\]
In particular, there are isomorphisms
\begin{gather*}
\mathfrak{D}^{\!\vee}:\mathscr{L}\left[K,V\right]\rightarrow\mathscr{L}\left[V^{*},K\right],\quad\left(v^{\mathrm{a}}\mapsto v^{\mathrm{a}^{\!\vee}}=v_{\mathrm{a^{*}}}\right)=\left(\left(\eta\mapsto\eta\mathbf{v}\right)\mapsto\left(\boldsymbol{f}\mapsto\boldsymbol{f}\!\left(\mathbf{v}\right)\right)\right),\\
\mathfrak{D}_{\!\vee}:\mathscr{L}\left[V,K\right]\rightarrow\mathscr{L}\left[K,V^{*}\right],\quad\left(f_{\mathrm{a}}\mapsto f_{\mathrm{a}^{\!\vee}}=f^{\mathrm{a}^{*}}\right)=\left(\left(\mathbf{v}\mapsto\boldsymbol{f}\!\left(\mathbf{v}\right)\right)\mapsto\left(\eta\mapsto\eta\boldsymbol{f}\right)\right).
\end{gather*}
We call $v^{\mathrm{a}^{\!\vee}}=v_{\mathrm{a^{*}}}$ and $f_{\mathrm{a}^{\!\vee}}=f^{\mathrm{a^{*}}}$
the \emph{dual tensor maps} (or \emph{outer transposes)} corresponding
to $v^{\mathrm{a}}$ and $f_{\mathrm{a}}$, respectively. 

Note that while $v^{\mathrm{a}}$ and $f_{\mathrm{a}}$ are tensor
maps on $V$, $v_{\mathrm{a^{*}}}$ and $f^{\mathrm{a^{*}}}$ are
tensor maps on $V^{*}$; indices are embellished with asterisks to
emphasize that we are dealing with $V^{*}$ rather than $V$. Specifically,
the dual of a vector-like tensor map $v^{\mathrm{a}}$ on $V$ is
a linear form $v_{\mathrm{a^{*}}}$ on $V^{*}$, while the dual of
a linear form $f_{\mathrm{a}}$ on $V$ is a vector-like tensor map
$f^{\mathrm{a}^{*}}$ on $V^{*}$. The connection between tensor maps
and their duals is illustrated by the fact that
\begin{gather*}
f_{\mathrm{a}}\circ v^{\mathrm{a}}=\left(\mathbf{v}\mapsto\boldsymbol{f}\!\left(\mathbf{v}\right)\right)\circ\left(\eta\mapsto\eta\mathbf{v}\right)=\left(\eta\mapsto\eta\boldsymbol{f}\!\left(\mathbf{v}\right)\right)\\
=\left(\boldsymbol{f}\mapsto\boldsymbol{f}\!\left(\mathbf{v}\right)\right)\circ\left(\eta\mapsto\eta\boldsymbol{f}\right)=v_{\mathrm{a^{*}}}\circ f^{\mathrm{a^{*}}}\!.
\end{gather*}

\emph{(2)}. The notions introduced above can be generalized. For finite-dimensional
vector spaces we have canonical isomorphisms $i':V\rightarrow V^{**}$
given by $i'\!\left(\mathbf{v}\right)\!\left(\boldsymbol{f}\right)=\boldsymbol{f}\!\left(\mathbf{v}\right)$
and $i'':\mathscr{L}\left[U\times V,W\right]\rightarrow\mathscr{L}\left[U,\mathscr{L}\left[V,W\right]\right]$
given by $i''\!\left(\mu\right)\!\left(\mathbf{u}\right)\!\left(\mathbf{v}\right)=\mu\!\left(\mathbf{u},\mathbf{v}\right)$.
Hence, there are canonical isomorphisms
\[
\mathscr{L}\left[V,V\right]\rightarrow\mathscr{L}\left[V,V^{**}\right]=\mathscr{L}\left[V,\mathscr{L}\left[V^{*},K\right]\right]\rightarrow\mathscr{L}\left[V\times V^{*},K\right]\rightarrow\mathscr{L}\left[V\otimes V^{*},K\right],
\]
so there are mappings 
\[
t_{\mathrm{a}}^{\mathrm{b}}\mapsto t_{\mathrm{a}\mathrm{b}^{*}}\mapsto t^{\left(\mathrm{ab^{*}}\right)^{*}}\mapsto t^{\mathrm{a^{*}b^{**}}}\mapsto t^{\mathrm{a^{*}b}},
\]
where $t^{\left(\mathrm{ab^{*}}\right)^{*}}$ is a map $K\rightarrow\left(V\otimes V^{*}\right)^{*}$
and $t^{\mathrm{a^{*}b^{**}}}$ a map $K\rightarrow V^{*}\otimes V^{**}$.
Specifically, given $t_{\mathrm{a}}^{\mathrm{b}}:V\rightarrow V$,
we can let $t_{\mathrm{a}}^{\mathrm{b\!}^{\vee}}=t_{\mathrm{a}\mathrm{b}^{*}}$
be a corresponding map $V\otimes V^{*}\rightarrow K$ and $t_{\mathrm{a}^{\vee}}^{\mathrm{b}}=t^{\mathrm{a}^{*}\mathrm{b}}$
a corresponding map $K\rightarrow V^{*}\otimes V$ . Similarly, given
$t_{\mathrm{a}}^{\mathrm{bc}}:V\rightarrow V\otimes V$, we can let
$t_{\mathrm{a}}^{\mathrm{b}\mathrm{c}^{\!\vee}}=t_{\mathrm{ac^{*}}}^{\mathrm{b}}$
be a corresponding map $V\otimes V^{*}\rightarrow V$, $t_{\mathrm{a}^{\!\vee}}^{\mathrm{bc}}=t^{\mathrm{a}^{*}\mathrm{b}\mathrm{c}}$
a corresponding map $K\rightarrow V^{*}\otimes V\otimes V$, and so
forth. Note that lowered superscripts are placed after any subscripts,
while raised subscripts are placed before any superscripts.

In general, we can transpose a tensor map with regard to one or more
individual indices; the replacement of a subscript $\mathsf{a}$ with
a superscript $\mathsf{a}^{*}$ or of a superscript $\mathsf{b}$
with a subscript $\mathsf{b}^{*}$ can be called \emph{shifting of
indices}. Note, though, that the resulting mappings are not tensor
maps in the sense defined earlier but generalized, \emph{comprehensive}
\emph{tensor maps}, defined on copies of $V$ and/or $V^{*}$. 
\begin{rem*}
In particular, given $t_{\mathrm{a}_{1}\ldots\mathrm{a_{m}}}^{\mathrm{b}_{1}\ldots\mathrm{b}_{n}}:V^{\otimes m}\rightarrow V^{\otimes n}$
we obtain by shifting of all superscripts a map
\[
t_{\mathrm{a}_{1}\ldots\mathrm{a_{m}\mathrm{b}_{1}^{*}\ldots\mathrm{b}_{n}^{*}}}:V^{\otimes n}\otimes V^{*\otimes m}\rightarrow K,
\]
corresponding to a separately linear map
\[
\tau_{\mathrm{a}_{1},\ldots,\mathrm{a_{m}},\mathrm{b}_{1}^{*},\ldots,\mathrm{b}_{n}^{*}}:V^{n}\times V^{*m}\rightarrow K.
\]
Conversely, tensor maps can be defined in terms of separately linear
maps of this form. (That bilateral tensors can be defined similarly
was noted earlier.)
\end{rem*}

\subsection{Adjoint tensor maps: raising and lowering indices}

$ $\\
\emph{(1)}. Let $\boldsymbol{g}:V\times V\rightarrow K$ be a bilinear
tensor map, and assume that $\boldsymbol{g}\left(\mathbf{u},\mathbf{v}\right)=\boldsymbol{g}\left(\mathbf{v},\mathbf{u}\right)$
and that $\boldsymbol{g}\left(\mathbf{u},\mathbf{v}\right)=0$ for
all $\mathbf{v}\in V$ implies $\mathbf{u}=0$; then $\boldsymbol{g}$
and the corresponding linear tensor map $g_{\mathrm{ab}}:V\otimes V\rightarrow K$
are said to be \emph{symmetric} and \emph{non-degenerate}.

We can define a linear map $I_{\boldsymbol{g}}:V\rightarrow V^{*}$
by setting
\[
I_{\boldsymbol{g}}\!\left(\mathbf{u}\right)\!\left(\mathbf{v}\right)=\boldsymbol{g}\!\left(\mathbf{u},\mathbf{v}\right)\quad\forall\mathbf{u}\in V,\;\forall\mathbf{v}\in V.
\]
If $\boldsymbol{g}$ is symmetric then this function is identical
to that obtained by substituting $\boldsymbol{g}\left(\mathbf{v},\mathbf{u}\right)$
for $\boldsymbol{g}\left(\mathbf{u},\mathbf{v}\right)$. Also, if
$V$ is finite-dimensional and $\boldsymbol{g}$ non-degenerate then
$I_{\boldsymbol{g}}$ is an isomorphism by (a special case of) the
Riesz representation theorem. In the finite-dimensional case, we thus
have a preferred isomorphism $V\rightarrow V^{*}$ after a choice
of a symmetric non-degenerate tensor map $g_{\mathrm{ab}}$. Hence,
we can identify $V$ and $V^{*}$, and any mapping which sends a tensor
map $t:V^{\otimes m}\rightarrow V^{\otimes n}$ to a dual map $\mathfrak{D}\!\left(\boldsymbol{t}\right):V^{*\otimes n}\rightarrow V^{*\otimes m}$
can be reinterpreted as a mapping which sends $t:V^{\otimes m}\rightarrow V^{\otimes n}$
to an \emph{adjoint map} (or \emph{inner transpose}) $\mathfrak{A}\!\left(\boldsymbol{t}\right):V^{\otimes n}\rightarrow V^{\otimes m}$. 

\emph{(2)}. Let us now reinterpret the mappings $\mathfrak{D}^{\!\vee}$
and $\mathfrak{D}_{\!\vee}$ in the light of the observation just
made. Set $u^{\mathrm{a}}=\left(\eta\mapsto\eta\mathbf{u}\right)$
and $v^{\mathrm{b}}=\left(\eta\rightarrow\eta\mathbf{v}\right)$.
Then $g_{\mathrm{ab}}\circ u^{\mathrm{a}}\circ v^{\mathrm{b}}=\left(\eta\mapsto\eta\boldsymbol{g}\!\left(\mathbf{u},\mathbf{v}\right)\right)$
and $\left(g_{\mathrm{ab}}\circ u^{\mathrm{a}}\right)\left(\mathbf{v}\right)=\boldsymbol{g}\!\left(\mathbf{u},\mathbf{v}\right)=I_{\boldsymbol{g}}\!\left(\mathbf{u}\right)\!\left(\mathbf{v}\right)$,
so $u^{\mathrm{a}}\mapsto f_{\mathrm{b}}=g_{\mathrm{ab}}\circ u^{\mathrm{a}}$
is the mapping $\left(K\rightarrow V\right)\rightarrow V^{*}$ corresponding
to $I_{\boldsymbol{g}}$. Also, if $\boldsymbol{g}$ is non-degenerate,
there is a mapping $I_{\boldsymbol{g}}^{-1}:V^{*}\rightarrow V$ such
that $I_{\boldsymbol{g}}^{-1}\circ I_{\boldsymbol{g}}$ is the identity
map on $V$, and there is a tensor map $g^{\mathrm{bc}}:K\rightarrow V\otimes V$
such that $g_{\mathrm{ab}}\circ g^{\mathrm{bc}}:V\rightarrow V$ is
the identity map. This means that corresponding to $g_{\mathrm{ab}}$
there is a tensor map $g^{\mathrm{ab}}$ such that $f_{\mathrm{a}}\mapsto v^{\mathrm{b}}=f_{\mathrm{a}}\circ g^{\mathrm{ab}}$
is the mapping $V^{*}\rightarrow\left(K\rightarrow V\right)$ corresponding
to $I_{\boldsymbol{g}}^{-1}$. We can thus define mappings
\begin{gather*}
\mathfrak{A}^{\wedge}:\mathscr{L}\left[K,V\right]\rightarrow\mathscr{L}\left[V,K\right],\quad v^{\mathrm{a}}\mapsto v^{\mathrm{a}^{\wedge}}=f_{\mathrm{a'}}=g_{\mathrm{a}\mathrm{a'}}\circ v^{\mathrm{a}},\\
\mathfrak{A}_{\wedge}:\mathscr{L}\left[V,K\right]\rightarrow\mathscr{L}\left[K,V\right],\quad f_{\mathrm{a}}\mapsto f_{\mathrm{a}^{\wedge}}=v^{\mathrm{a'}}=f_{\mathrm{a}}\circ g^{\mathrm{aa'}}.
\end{gather*}

In classical tensor terminology, we say that $\mathfrak{A}^{\wedge}$
effectuates the \emph{lowering of indices by }$\boldsymbol{g}$, while
$\mathfrak{A}_{\wedge}$ effectuates the \emph{raising of indices
by} $\boldsymbol{g}$.

\emph{(3)}. We can raise and lower indices on any tensor maps, and
we can use $\boldsymbol{g}$ several times on the same tensor map
to raise or lower indices. For example, 
\begin{gather*}
t_{\mathrm{a}}^{\mathrm{bc}^{\wedge}}=g_{\mathrm{cc'}}\circ t_{\mathrm{a}}^{\mathrm{bc}}=t_{\mathrm{c'a}}^{\mathrm{b}}\,,\qquad t_{\mathrm{a}^{\wedge}\mathrm{b}^{\wedge}}^{\mathrm{c}}=t_{\mathrm{a}\mathrm{b}}^{\mathrm{c}}\circ g^{\mathrm{aa'}}\circ g^{\mathrm{bb'}}=t^{\mathrm{ca'b'}},\\
t_{\mathrm{a}^{\wedge}}^{\mathrm{b}^{\wedge}}=g_{\mathrm{bb'}}\circ t_{\mathrm{a}}^{\mathrm{b}}\circ g^{\mathrm{aa'}}=t_{\mathrm{b'}}^{\mathrm{a'}}\,.
\end{gather*}
In general,
\[
t_{\ldots\mathrm{a}_{i}^{\wedge}\ldots}^{\ldots\mathrm{b}_{j}^{\wedge}\ldots}=\ldots\circ g_{\mathrm{b}_{j}\mathrm{b}_{j}'}\circ\ldots\circ t_{\ldots\mathrm{a}_{i}\ldots}^{\ldots\mathrm{b}_{j}\ldots}\circ\ldots\circ g^{\mathrm{a}_{i}\mathrm{a}_{i}'}\circ\ldots=t_{\ldots\mathrm{b}_{j}'\ldots}^{\ldots\mathrm{a}_{i}'\ldots}\,,
\]
where lowered subscripts are placed before any original subscripts,
while raised subscripts are placed after any original superscripts.
\begin{rem*}
Given $t_{\mathrm{a}_{1}\ldots\mathrm{a_{m}}}^{\mathrm{b}_{1}\ldots\mathrm{b}_{n}}:V^{\otimes m}\rightarrow V^{\otimes n}$,
we obtain by lowering of all superscripts a map
\[
t_{\mathrm{b}_{1}'\ldots\mathrm{b}_{n}'\mathrm{a}_{1}\ldots\mathrm{a_{m}}}:V^{\otimes n}\otimes V^{\otimes m}\rightarrow K,
\]
corresponding to a separately linear map
\[
\tau_{\mathrm{b}_{1}',\ldots,\mathrm{b}_{n}',\mathrm{a}_{1},\ldots,\mathrm{a_{m}}}:V^{n}\times V^{m}\rightarrow K.
\]
Conversely, tensor maps can be defined in terms of separately linear
maps of this form. Thus, when a vector space $V$ is equipped with
a symmetric non-degenerate bilinear form $\boldsymbol{g}$ which makes
it possible to identify $V^{*}$ with $V$, any tensor map on $V$
can be regarded as a linear or separately linear form on $V$. 
\end{rem*}
\emph{(3)}. We can combine raising and lowering of indices with permutations
of indices in various ways. Let us look at two simple cases:

(a) We can lower or raise an index first and then permute subscripts
or superscripts in the resulting tensor map. For example, compare
\[
g_{\mathrm{bb'}}\circ t_{\mathrm{a}}^{\mathrm{b}}=t_{\mathrm{b'a}}\quad\mathrm{and}\quad\left(g_{\mathrm{bb'}}\circ t_{\mathrm{a}}^{\mathrm{b}}\right)\circ p_{\mathrm{a'b''}}^{\mathrm{b'a}}=t_{\mathrm{b'a}}\circ p_{\mathrm{a'b''}}^{\mathrm{b'a}}=t_{\mathrm{\mathrm{a}'b''}}\,.
\]

(b) We can permute subscripts or superscripts in a given tensor map
first and then raise or lower indices. For example, compare
\[
g_{\mathrm{bb'}}\circ t^{\mathrm{ba}}=t_{\mathrm{b}'}^{\mathrm{a}}\qquad\mathrm{and}\qquad g_{\mathrm{b'b''}}\circ\left(p_{\mathrm{ba}}^{\mathrm{a'b'}}\circ t^{\mathrm{ba}}\right)=g_{\mathrm{b'b''}}\circ t^{\mathrm{a'b'}}=t_{\mathrm{b}''}^{\mathrm{a'}}\,.
\]

In case (a), we can write $t_{\mathrm{b'a}}$ as $t_{\mathrm{ba}}$
and $t_{\mathrm{\mathrm{a}'b''}}$ as $t_{\mathrm{ab}}$ or $ $$t_{\mathrm{a}^{\!2}\mathrm{b}^{\!1}}$
to avoid ambiguity, relating $t_{\mathrm{ab}}$ to the reference tensor
$t_{\mathrm{ba}}$ as discussed in Subsection 6.3(2). In case (b),
the ambiguity problem is even worse, since we can write both $t_{\mathrm{b}'}^{\mathrm{a}}$
and $t_{\mathrm{b}''}^{\mathrm{a'}}$ as $t_{\mathrm{b}}^{\mathrm{a}}$,
and the two instances of $t_{\mathrm{b}}^{\mathrm{a}}$ clearly have
different meanings. To eliminate this ambiguity, we can again label
indices with respect to how they are arranged in a reference tensor.
Thus, we can write the first $t_{\mathrm{b}}^{\mathrm{a}}$ as $t_{\mathrm{b}^{\!1}}^{\mathrm{a}^{\!2}}$,
since the reference tensor is $t^{\mathrm{ba}}$, and the second $t_{\mathrm{b}}^{\mathrm{a}}$
as $t_{\mathrm{b}^{\!2}}^{\mathrm{a}^{\!1}}$, since the reference
tensor this time is $t^{\mathrm{ab}}$.
\begin{rem*}
Alternatively, one can use stacked indices, writing $t_{\mathrm{b}^{\!1}}^{\mathrm{a}^{\!2}}$
as $t_{\mathrm{b}}^{\phantom{\mathrm{b}}\mathrm{a}}$ and $t_{\mathrm{b}^{\!2}}^{\mathrm{a}^{\!1}}$
as $t_{\phantom{\mathrm{a}}\mathrm{b}}^{\mathrm{a}}$. The underlying
issue here is really the same as that noted in Subsection 6.3(2),
however, namely the contextual interpretation of indices, so there
is a reason to use the same kind of notation in both cases.\newpage{}
\end{rem*}

\specialsection*{\textbf{\textsl{C. Representations of classical tensors}}}

\section{Arrays and array-like notation}

\subsection{Arrays as generalized matrices}

\emph{$ $}\\
A \emph{plain matrix} is a finite collection of not necessarily distinct
elements called \emph{entries}, arranged in one or more rows and columns
as exemplified below: 
\[
R=\left[\begin{array}{cc}
1 & 2\end{array}\right],\: C=\left[\begin{array}{c}
3\\
4
\end{array}\right],\: X=\left[\begin{array}{cc}
1 & 2\\
3 & 4
\end{array}\right],\: Z=\left[\begin{array}{cccc}
1 & 2 & 3 & 4\\
5 & 6 & 7 & 8\\
9 & 0 & 1 & 2\\
3 & 4 & 5 & 6
\end{array}\right].
\]
Each plain matrix can be converted into an \emph{array} by adding
an \emph{indexation}. As an example, $\left[x_{i}^{j}\right]$ is
an array given by an indexation of $X$ such that $x_{1}^{1}=1$,
$x_{1}^{2}=2$, $x_{2}^{1}=3$ and $x_{2}^{2}=4$. In all such indexations,
superscripts are column indices, identifying entries within a row,
while subscripts are row indices, identifying entries within a column.

Each one of the plain matrices $R$, $C$, and $X$ has only one possible
indexation, so that the distinction between plain and indexed matrices
(arrays) appears to be redundant. In other cases, however, a plain
matrix can be indexed in different ways, so that one plain matrix
corresponds to several arrays. For example, $Z$ can be indexed as
$\left[z{}_{i}^{j}\right]$, where both indices run from 1 to 4, or
as $\left[z{}_{ij}^{k\ell}\right]$, where all indices run from 1
to 2. We shall assume here that all indices in an array $A$ run from
$1$ to some $N_{A}>1$. Indices are ordered lexicographically within
rows and columns so that, for example, 
\[
\left[z_{ij}^{k\ell}\right]=\left[\begin{array}{cccc}
z_{11}^{11} & z_{11}^{12} & z_{11}^{21} & z_{11}^{22}\\
z_{12}^{11} & z_{12}^{12} & z_{12}^{21} & z_{12}^{22}\\
z_{21}^{11} & z_{21}^{12} & z_{21}^{21} & z_{21}^{22}\\
z_{22}^{11} & z_{22}^{12} & z_{22}^{21} & z_{22}^{22}
\end{array}\right].
\]

A general array is denoted 
\[
\left[a_{i_{1}\ldots i_{m}}^{j_{1}\ldots j_{n}}\right]\qquad\left(m,n\geq0\right).
\]
 An array of this form is said to have \emph{valence} $\binom{n}{m}$.
Arrays of the forms $\left[a_{i}^{j}\right]$, $\left[a^{j}\right]$
and $\left[a_{i}\right]$ correspond to square matrices, row matrices
and column matrices, respectively. A singleton array $\left[a\right]$
is by convention identified with the entry $a$.

As usual, a multi-index $I$ can refer to a (possibly empty) tuple
of indices $i_{1}\ldots i_{n}$. By an abuse of notation we can write
$\left[c_{I}^{J}\right]$ where $c_{I}^{J}=\lambda a_{I}^{J}$ as
$\left[\lambda a_{I}^{J}\right]$, $\left[c_{I}^{J}\right]$ where
$c_{I}^{J}=a_{I}^{J}+b_{I}^{J}$ as $\left[a_{I}^{J}+b_{I}^{J}\right]$,
$\left[c_{I}^{J}\right]$ where $c_{I}^{J}=f\!\left(x_{I}^{J}\right)$
as $\left[f\!\left(x_{I}^{J}\right)\right]$ etc.

Arrays of scalars will be used to represent tensor maps and related
objects, but note that the entries of an array need not be scalars.
We will make use of arrays the elements of which are scalars, vectors,
tensors or tensor maps. An ordered basis for a vector space can thus
be regarded as an array $\left[\mathbf{e}_{i}\right]$, where each
entry $\mathbf{e}_{i}$ is a basis vector. (This is necessarily an
array with distinct entries.)

\subsection{Sets of arrays as vector spaces}

$ $\\
If $\eta$ is a scalar and $\left[a{}_{i_{1}\cdots i_{m}}^{j_{1}\cdots j_{n}}\right]$
an array of entries that allow scalar multiplication, we can define
their scalar product by 
\[
\eta\left[a{}_{i_{1}\cdots i_{m}}^{j_{1}\cdots j_{n}}\right]=\left[\eta a{}_{i_{1}\cdots i_{m}}^{j_{1}\cdots j_{n}}\right],
\]
and if $\left[a{}_{i_{1}\cdots i_{m}}^{j_{1}\cdots j_{n}}\right]$
and $\left[b{}_{i_{1}\cdots i_{m}}^{j_{1}\cdots j_{n}}\right]$ are
arrays of entries that can be added, we can set 
\[
\left[a{}_{i_{1}\cdots i_{m}}^{j_{1}\cdots j_{n}}\right]+\left[b{}_{i_{1}\cdots i_{m}}^{j_{1}\cdots j_{n}}\right]=\left[a{}_{i_{1}\cdots i_{m}}^{j_{1}\cdots j_{n}}+a{}_{i_{1}\cdots i_{m}}^{j_{1}\cdots j_{n}}\right].
\]

A set of arrays of the same valence and with the same kind of entries,
allowing addition and multiplication with scalars from $K$, obviously
constitutes a vector space over $K$. The vector space of arrays of
valence $\binom{n}{m}$ with entries in $K$ and indices ranging from
$1$ to $N$ is denoted $\mathscr{A}_{m}^{n}\left[K^{N}\right]$.

\subsection{Multiplication of arrays}

$ $\\
\emph{(1).}~~Consider scalar arrays $\left[a{}_{J}^{i_{1}\ldots i_{k}K}\right],\left[b{}_{i_{1}\ldots i_{k}J'}^{K'}\right]$,
where the multi-indices $J,J',K,K'$ have no individual indices in
common. We define the product of these arrays by
\[
\left[a{}_{J}^{i_{1}\ldots i_{k}K}\right]\left[b{}_{i_{1}\ldots i_{k}J'}^{K'}\right]=\left[\sum_{i_{1},\ldots,i_{k}}a{}_{J}^{i_{1}\ldots i_{k}K}b{}_{i_{1}\ldots i_{k}J'}^{K'}\right]=\left[c{}_{JJ'}^{KK'}\right].
\]
For example, $\left[a^{i}\right]\left[b_{i}\right]=\left[\sum_{i}a^{i}b_{i}\right]=\left[c\right]$,
$\left[a^{ij}\right]\left[b_{i}\right]=\left[\sum_{i}a^{ij}b_{i}\right]=\left[c^{j}\right]$
and $\left[a_{j}^{i}\right]\left[b_{i}^{k}\right]=\left[\begin{array}{c}
\sum_{i}a_{j}^{i}b_{i}^{k}\end{array}\right]=\left[c_{j}^{k}\right]$. As we did when we defined general multiplication of bilateral tensors
and tensor maps, we can generalize this definition so that it applies
to cases where the indices summed over are not necessarily the first
$k$ superscripts and subscripts. Multiplication of arrays which have
no indices in common is also possible. In this case, entries are multiplied
but no summation occurs. For example, 
\[
\left[a_{i}^{j}\right]\left[b_{k}^{\ell}\right]=\begin{array}{c}
\left[a_{i}^{j}b_{k}^{\ell}\right]\end{array}=\left[c_{ik}^{j\ell}\right].
\]

It is easy to verify that multiplication of scalar arrays is associative
but not commutative.

\emph{(2)}. The concept of matrix inverses can be generalized to inverses
of general scalar arrays of valence $\binom{n}{n}$. Set
\[
\delta_{i_{1}\ldots i_{n}}^{j_{1}\ldots j_{n}}=\begin{cases}
1 & \mathrm{if\;}i_{1}=j_{1},\ldots,i_{n}=j_{n},\\
0 & \mathrm{otherwise.}
\end{cases}
\]
 The \emph{inverse} $\left[a_{i_{1}\ldots i_{n}}^{\ell_{1}\ldots\ell_{n}}\right]^{-1}$
of $\left[a_{j_{1}\ldots j_{n}}^{k_{1}\ldots k_{n}}\right]$ is an
array $\left[\alpha{}_{i_{1}\ldots i_{n}}^{\ell_{1}\ldots\ell_{n}}\right]$
such that 
\[
\left[a_{j_{1}\ldots j_{n}}^{k_{1}\ldots k_{n}}\right]\left[\alpha{}_{k_{1}\ldots k_{n}}^{\ell_{1}\ldots\ell_{n}}\right]=\left[\delta_{j_{1}\ldots j_{n}}^{\ell_{1}\ldots\ell_{n}}\right],\qquad\left[\alpha{}_{i_{1}\ldots i_{n}}^{j_{1}\ldots j_{n}}\right]\left[a_{j_{1}\ldots j_{n}}^{k_{1}\ldots k_{n}}\right]=\left[\delta_{i_{1}\ldots i_{n}}^{k_{1}\ldots k_{n}}\right].
\]
 It is clear that this definition agrees with the usual definition
of the inverse of a square matrix $\left[\alpha_{i}^{j}\right]$.

\emph{(3)}. Multiplication of other arrays than scalar arrays is also
possible, provided that array elements can be multiplied in some sense
and the resulting products added. In particular, a scalar array can
be multiplied with an array of vectors, tensors or tensor maps, based
on scalar multiplication of individual entries; arrays of tensors
can be multiplied, based on tensor multiplication of entries; and
arrays of tensor maps can be multiplied, based on composition of entries
as tensor maps. For example,
\[
\left[f_{\mathrm{a}}^{i}\right]\circ\left[e_{i}^{\mathrm{a}}\right]=\left[\sum_{i}f_{\mathrm{a}}^{i}\circ e_{i}^{\mathrm{a}}\right]=\left[f_{\mathrm{a}}^{1}\circ e_{1}^{\mathrm{a}}+...+f_{\mathrm{a}}^{n}\circ e_{n}^{\mathrm{a}}\right]=\left[s_{1}+\ldots+s_{n}\right]=s,
\]
where $s_{i}=\left(\eta\mapsto\eta\sigma_{i}\right)$ and $s=\left(\eta\mapsto\eta\sigma\right)$
are scalar-like tensor maps.

\subsection{Transposes of arrays}

$ $\\
\emph{(1)}. We can define different kinds of \emph{transposes} of
arrays, analogous to the transpose of a matrix. For example, $\left[a_{i}^{j}\right]$
has three possible transposes:
\begin{enumerate}
\item An array $\left[a_{i}^{j^{\mathrm{T}}}\right]=\left[\alpha_{ij}\right]$
, where $\alpha{}_{ij}=a_{i}^{j}$.
\item An array $\left[a_{i^{\mathrm{T}}}^{j}\right]=\left[\alpha^{ij}\right]$
, where $\alpha^{ij}=a_{i}^{j}$.
\item An array $\left[a_{i^{\mathrm{T}}}^{j^{\mathrm{T}}}\right]=\left[\alpha_{j}^{i}\right]$
, where $\alpha_{j}^{i}=a_{i}^{j}$. 
\end{enumerate}
Transposition of general arrays is defined similarly; it involves
lowering of superscripts and raising of subscripts. By convention,
lowered superscripts are placed after any subscripts, and raised subscripts
are placed before any superscripts. 

The array $\left[a_{i^{\mathrm{T}}}^{j^{\mathrm{T}}}\right]$ is the
usual matrix transpose $A^{\mathrm{T}}$ of the matrix-like array
$A=\left[a_{i}^{j}\right]$, and it is convenient to write $\left[a_{i_{1}\ldots i_{m}}^{j_{1}\ldots j_{n}}\right]^{\mathrm{T}}$
instead of $\left[a_{i_{1}^{\mathrm{T}}\ldots i_{m}^{\mathrm{T}}}^{j_{1}^{\mathrm{T}}\ldots j_{n}^{\mathrm{T}}}\right]$
when all indices in an array are transposed.

\emph{(2)}. Because of the way multiplication of arrays is defined
we have, for example
\[
\left(\left[a_{j}^{i}\right]\left[b_{i}^{k}\right]\right)^{\mathrm{T}}=\left[b_{i}^{k}\right]^{\mathrm{T}}\left[a_{j}^{i}\right]^{\mathrm{T}},
\]
which is a well-known identity from matrix algebra, but the definition
of array multiplication does also imply other types of identities
such as 
\[
\left[a^{i}\right]\!\left[b_{ij}\right]=\left[b_{i^{\mathrm{T}}j}\right]\!\left[a^{i^{\mathrm{T}}}\right],\quad\left[a_{j}^{i}\right]\!\left[b_{i}^{k}\right]=\left[b_{i^{\mathrm{T}}}^{k}\right]\!\left[a_{j}^{i^{\mathrm{T}}}\right],\quad\left[u^{i}\right]\!\left[v^{j}\right]\!\left[g_{ij}\right]=\left[u^{i}\right]\!\left[g_{ij^{\mathrm{T}}}\right]\!\left[v^{j^{\mathrm{T}}}\right].
\]
 Note, however, that such identities presuppose that multiplication
of array entries is commutative.

\subsection{Conventional matrix notation and index-free array notation}

\emph{$ $}\\
\emph{(1)}. ~Array notation involves partly arbitrary choices of
indices; for example, $\left[a_{j}^{i}\right]\left[b_{i}^{k}\right]$
and $\left[a_{i}^{k}\right]\left[b_{k}^{\ell}\right]$ have the same
meaning although different indices are used. In other words, array
notation is characterized by some redundancy, and this makes it possible
to eliminate the indices in certain special situations by introducing
suitable conventions. For example, each of the array expressions just
considered can be written in conventional matrix notation as $AB$,
using the 'row by column' convention for matrix multiplication. 

Using fraktur style but otherwise adhering to the conventional notation,
we let $\boldsymbol{\mathfrak{A}},\boldsymbol{\mathfrak{B}},...$
represent square arrays $\left[a_{i}^{j}\right],\left[b_{i}^{j}\right],\ldots\:$,
let $\boldsymbol{\mathfrak{a}},\boldsymbol{\mathfrak{b}},...$ represent
column arrays $\left[a_{i}\right],\left[b_{i}\right],\ldots\:$, and
let $\boldsymbol{\mathfrak{a}}^{\mathrm{T}},\boldsymbol{\mathfrak{b}}^{\mathrm{T}},...$
represent row arrays $\left[a^{i}\right],\left[b^{i}\right],\ldots\;$.
For example, $\boldsymbol{\mathfrak{a}}{}^{\mathrm{T}}\boldsymbol{\mathfrak{B}}$,
$\boldsymbol{\mathfrak{B}}\boldsymbol{\mathfrak{a}}$ and $\boldsymbol{\mathfrak{A}}\boldsymbol{\mathfrak{B}}$
represent essentially unique array expressions such as $\left[a^{i}\right]\left[b_{i}^{j}\right]$,
$\left[b_{i}^{j}\right]\left[a_{j}\right]$ and $\left[a_{j}^{i}\right]\left[b_{i}^{k}\right]$,
respectively. 

\emph{(2)}. ~Going one step further, we can generalize conventional
matrix notation to \emph{index-free array notation,} where we write
$\left[a_{i_{1}...i_{m}}^{j_{1}...j_{n}}\right]$ as $\underset{m}{\overset{n}{\boldsymbol{\mathfrak{A}}}}$;
in particular, $\left[a\right]=\underset{0}{\overset{0}{\boldsymbol{\mathfrak{A}}}}$,
$\left[a_{i}\right]=\underset{1}{\overset{0}{\boldsymbol{\mathfrak{A}}}}$,
$\left[a^{j}\right]=\underset{0}{\overset{1}{\boldsymbol{\mathfrak{A}}}}$,
and $\left[a_{i}^{j}\right]=\underset{1}{\overset{1}{\boldsymbol{\mathfrak{A}}}}$
(although one may use forms such as $\boldsymbol{\mathfrak{a}}$,
$\boldsymbol{\mathfrak{a}}^{\mathrm{T}}$ and $\boldsymbol{\mathfrak{A}}$
together with forms such as $\underset{m}{\overset{n}{\boldsymbol{\mathfrak{A}}}}$
for convenience). 

The \emph{inner product} of $\underset{m}{\overset{n}{\boldsymbol{\mathfrak{A}}}}=\left[a_{i_{1}...i_{m}}^{j_{1}...j_{n}}\right]$
and $\underset{n}{\overset{p}{\boldsymbol{\mathfrak{B}}}}=\left[b_{k_{1}...k_{n}}^{\ell_{1}...\ell_{p}}\right]$
is
\[
\underset{m}{\overset{n}{\boldsymbol{\mathfrak{A}}}}\underset{n}{\overset{p}{\boldsymbol{\mathfrak{B}}}}=\left[a_{i_{1}...i_{m}}^{j_{1}...j_{n}}\right]\left[b_{j_{1}...j_{n}}^{\ell_{1}...\ell_{p}}\right]=\left[\sum_{j_{1},\ldots,j_{n}}a_{i_{1}...i_{m}}^{j_{1}...j_{n}}b_{j_{1}...j_{n}}^{\ell_{1}...\ell_{p}}\right]=\left[c_{i_{1}...i_{m}}^{\ell_{1}...\ell_{p}}\right]=\underset{m}{\overset{p}{\boldsymbol{\mathfrak{C}}}}.
\]
The inner product obviously generalizes ordinary matrix multiplication. 

The \emph{outer product} of $\underset{m}{\overset{n}{\boldsymbol{\mathfrak{A}}}}=\left[a_{i_{1}...i_{m}}^{j_{1}...j_{n}}\right]$
and $\underset{p}{\overset{q}{\boldsymbol{\mathfrak{B}}}}=\left[b_{k_{1}...k_{p}}^{\ell_{1}...\ell_{q}}\right]$
is
\[
\underset{m}{\overset{n}{\boldsymbol{\mathfrak{A}}}}\otimes\underset{p}{\overset{q}{\boldsymbol{\mathfrak{B}}}}=\left[a_{i_{1}...i_{m}}^{j_{1}...j_{n}}\right]\left[b_{k_{1}...k_{p}}^{\ell_{1}...\ell_{q}}\right]=\left[a_{i_{1}...i_{m}}^{j_{1}...j_{n}}b_{k_{1}...k_{p}}^{\ell_{1}...\ell_{q}}\right]=\left[c_{i_{1}...i_{m}k_{1}...k_{p}}^{j_{1}...j_{n}\ell_{1}...\ell_{q}}\right]=\underset{m+p}{\overset{n+q}{\boldsymbol{\mathfrak{C}}}}.
\]
The outer product generalizes the so-called Kronecker product of matrices.

It should be pointed out that an index-free notation for arrays has
the same kind of limitations as an index-free notation for classical
tensors. For example, the simple array product $\left[a^{ij}\right]\left[b_{j}\right]$
cannot be expressed as an ordinary matrix product, nor as an inner
product of arrays, nor as a Kronecker product of matrices, nor as
an outer product of arrays.

\subsection{Array notation and conventional indicial notation}

$ $\\
Conventional indicial notation traditionally used in tensor analysis
is an alternative to array notation. For example, using Einstein's
summation convention we can write
\[
v^{i}t_{i}^{j}
\]
instead of $\left.\sum\right._{i}v^{i}t_{i}^{j}$ or $\left[v^{i}\right]\left[t_{i}^{j}\right]$.
When using array notation instead of indicial notation, we require
that an index summed over occurs first as a superscript in one array,
then as a subscript in a subsequent array, so that we get 'row by
column' multiplication of two adjacent arrays. In indicial notation
with Einstein's summation convention, an index summed over may also
occur first as a subscript and then as a superscript. For example,
as $\left[a_{i}^{k}\right]\left[b_{k}^{\ell}\right]\left[c_{\ell}^{j}\right]=\left[\sum_{k,\ell}a_{i}^{k}b_{k}^{\ell}c_{\ell}^{j}\right]$,
and since the scalars in each term commute, we can write $\left[a_{i}^{k}\right]\left[b_{k}^{\ell}\right]\left[c_{\ell}^{j}\right]$
as $a_{i}^{k}b_{k}^{\ell}c_{\ell}^{j}$, $a_{i}^{k}c_{\ell}^{j}b_{k}^{\ell}$,
$b_{k}^{\ell}a_{i}^{k}c_{\ell}^{j}$, $b_{k}^{\ell}c_{\ell}^{j}a_{i}^{k}$,
$c_{\ell}^{j}a_{i}^{k}b_{k}^{\ell}$ or $c_{\ell}^{j}b_{k}^{\ell}a_{i}^{k}$
in indicial notation.
\begin{rem*}
Indicial notation is versatile and flexible. In particular, one does
not have to manipulate transposes explicitly. For example, instead
of rewriting $\left[u^{i}\right]\!\left[v^{j}\right]\!\left[g_{ij}\right]$
as $\left[u^{i}\right]\!\left[g_{ij^{\mathrm{T}}}\right]\!\left[v^{j^{\mathrm{T}}}\right]$
one simply identifies $u^{i}v^{j}g_{ij}$ with $u^{i}g_{ij}v^{j}$.
However, the fact that the notation conveys the impression that nothing
is happening in such cases may hide what is going on conceptually.
In addition, the fact that expressions in indicial notation have no
'canonical form' means that there may be some more-or-less arbitrary
choices to make when using this notation.

Partly because it is so compact, indicial notation can be ambiguous.
For example, $a^{i}b_{i}$ may mean (i) a scalar $a^{i}b_{i}$ for
a definite $i$, (ii) a row matrix $\left[c^{i}\right]$, where $c^{i}=a^{i}b_{i}$,
or a column matrix $\left[c_{i}\right]$, where $c_{i}=a^{i}b_{i}$,
or (iii) a singleton array $\left[a^{i}\right]\left[b_{i}\right]=\left[\sum_{i}a^{i}b_{i}\right]$.
The occasional note {}``summation not implied'' in texts using indicial
notation serves to distinguish between cases (ii) and (iii).

Last but not least, array notation is closer to conventional matrix
notation than indicial notation is, so array notation helps to preserve
the unity of mathematical notation, and array notation does not blur
the distinction between a tensor and its representation.
\end{rem*}

\section{Array-type representations of vectors, tensors and tensor maps }

\subsection{Array representations of vectors, tensors and tensor maps}

$ $\\
Below, let $\left[\mathsf{\mathbf{e}}_{i}\right]$ be an indexed basis
for $V$. Then $\left[\mathbf{e}{}_{i}\otimes\mathbf{e}{}_{j}\right]$
is an indexed basis for $V\otimes V$, etc. Furthermore, $\left[\boldsymbol{e}_{i}\right]=\left[\eta\mapsto\eta\mathbf{e}{}_{i}\right]$
is an indexed basis for $\mathscr{L}\left[K,V\right]$, $\left[\boldsymbol{e}_{i}\otimes\boldsymbol{e}_{j}\right]=\left[\eta\mapsto\eta\left(\mathbf{e}{}_{i}\otimes\mathbf{e}{}_{j}\right)\right]$
is an indexed basis for $\mathscr{L}\left[K,V\otimes V\right]$, and
so forth. Every vector, tensor and tensor map can be represented by
a scalar array relative to some basis as described in this subsection.

\emph{(1)}. There is a unique expansion of every $\mathbf{v}\in V$
with regard to $\left[\mathbf{e}{}_{i}\right]$, namely 
\[
\mathbf{v}=\left.\sum\right._{i}v^{i}\mathbf{e}_{i}=\left[v^{i}\right]\left[\mathbf{e}{}_{i}\right].
\]
In terms of the corresponding vector-like tensor map $\boldsymbol{v}=v^{\mathrm{a}}:K\rightarrow V$,
this expansion becomes
\[
\boldsymbol{v}=\left.\sum\right._{i}v^{i}\boldsymbol{e}_{i}=\left[v^{i}\right]\left[\boldsymbol{e}_{i}\right].
\]
This means that $\left[v^{i}\right]$ represents $\mathbf{v}$ relative
to $\left[\mathbf{e}{}_{i}\right]$ and $\boldsymbol{v}$ relative
to $\left[\boldsymbol{e}_{i}\right]$. Indirectly, $\left[v^{i}\right]$
also represents $\boldsymbol{v}$ relative to $\left[\mathbf{e}_{i}\right]$.

\emph{(2)}. As $\left[\mathbf{e}{}_{i}\otimes\mathbf{e}{}_{j}\right]$
is an indexed basis for $V\otimes V$, any $\mathbf{t}\in V\otimes V$
has a unique expansion
\[
\mathbf{t}=\left.\sum\right._{i,j}t^{ij}\left(\mathbf{e}{}_{i}\otimes\mathbf{e}{}_{j}\right)=\left[t^{ij}\right]\left[\mathbf{e}{}_{i}\otimes\mathbf{e}{}_{j}\right]
\]
relative to $\left[\mathbf{e}{}_{i}\otimes\mathbf{e}{}_{j}\right]$,
and ultimately with respect to $\left[\mathbf{e}{}_{i}\right]$. In
terms of the tensor map $\boldsymbol{t}=t^{\mathrm{ab}}:K\rightarrow V\otimes V$,
this expansion is
\[
t^{\mathrm{ab}}=\left.\sum\right._{i,j}t^{ij}\left(e_{i}^{\mathrm{a}}\circ e_{j}^{\mathrm{b}}\right)=\left[t^{ij}\right]\left[e_{i}^{\mathrm{a}}\circ e_{j}^{\mathrm{b}}\right],
\]
or, using index-free notation for tensor maps,
\[
\boldsymbol{t}=\left.\sum\right._{i,j}t^{ij}\left(\boldsymbol{e_{i}}\otimes\boldsymbol{e_{j}}\right)=\left[t^{ij}\right]\left[\boldsymbol{e_{i}}\otimes\boldsymbol{e_{j}}\right].
\]
Thus, both $\mathbf{t}$ and $\boldsymbol{t}$ are represented by
$\left[t^{ij}\right]$, ultimately with respect to $\left[\mathbf{e}{}_{i}\right]$.
This result can obviously be generalized to higher tensor powers of
$V$. 

\emph{(3)}. Consider a linear form $\boldsymbol{f}=f_{\mathrm{a}}:V\rightarrow K$.
We have
\[
\boldsymbol{f}\!\left(\mathbf{v}\right)=\boldsymbol{f}\!\left(\left.\sum\right._{i}v^{i}\mathbf{e}_{i}\right)=\left.\sum\right._{i}v^{i}\boldsymbol{f}\!\left(\mathbf{e}_{i}\right),
\]
 or in array notation, 
\[
\boldsymbol{f}\!\left(\mathbf{v}\right)=\boldsymbol{f}\!\left(\left[v^{i}\right]\left[\mathbf{e}{}_{i}\right]\right)=\left[v^{i}\right]\left[f_{i}\right],
\]
where $f_{i}=\boldsymbol{f}\!\left(\mathbf{e}_{i}\right)$. Since
$\left[v^{i}\right]$ represents $\mathbf{v}$ relative to $\left[\mathbf{e}{}_{i}\right]$,
$\left[f_{i}\right]$ represents $\boldsymbol{f}$ relative to $\left[\mathbf{e}{}_{i}\right]$. 

In terms of tensor map composition we have
\[
f_{\mathrm{a}}\circ v^{\mathrm{a}}=f_{\mathrm{a}}\circ\left[v^{i}\right]\left[e_{i}^{\mathrm{a}}\right]=\left[v^{i}\right]\left[f_{\mathrm{a}}\circ e_{i}^{\mathrm{a}}\right]=\left[v^{i}\right]\left[\boldsymbol{f}{}_{i}\right],
\]
or in index-free notation
\[
\boldsymbol{f}\circ\boldsymbol{v}=\boldsymbol{f}\circ\left[v^{i}\right]\left[\boldsymbol{e}_{i}\right]=\left[v^{i}\right]\left[\boldsymbol{f}\circ\boldsymbol{e}_{i}\right]=\left[v^{i}\right]\left[\boldsymbol{f}{}_{i}\right],
\]
where 
\[
\left[\boldsymbol{f}{}_{i}\right]=\left[\eta\mapsto\eta\boldsymbol{f}\!\left(\mathbf{e}_{i}\right)\right]=\left[\eta\mapsto\eta f_{i}\right]=\left[f_{i}\right]\left[\eta\mapsto\eta\right],
\]
so $\left[f_{i}\right]$ represents $\boldsymbol{f}$ relative to
$\left[\boldsymbol{e}{}_{i}\right]$ and thus relative to $\left[\mathbf{e}{}_{i}\right]$.

It follows that $\mathbf{v}\mapsto\boldsymbol{f}\!\left(\mathbf{v}\right)$
and $\boldsymbol{v}\mapsto\boldsymbol{f}\circ\boldsymbol{v}$ can
be represented in array form as 
\[
\left[v^{i}\right]\mapsto\left[v^{i}\right]\left[f_{i}\right],
\]
 and in matrix form as $\boldsymbol{\mathfrak{v}}^{\mathrm{T}}\mapsto\boldsymbol{\mathfrak{v}}^{\mathrm{T}}\boldsymbol{\mathfrak{f}}$
or, since $\left(\boldsymbol{\mathfrak{v}}^{\mathrm{T}}\boldsymbol{\mathfrak{f}}\right)^{\mathrm{T}}=\boldsymbol{\mathfrak{f}}{}^{\mathrm{T}}\boldsymbol{\mathfrak{v}}$,
as $\boldsymbol{\mathfrak{v}}\mapsto\boldsymbol{\mathfrak{f}}{}^{\mathrm{T}}\boldsymbol{\mathfrak{v}}$.

\emph{(4)}. Let $\boldsymbol{\gamma}=\gamma_{\mathrm{a,b}}:V\times V\rightarrow K$
be a separately linear (bilinear) map. Then
\[
\boldsymbol{\gamma}\!\left(\mathbf{u},\mathbf{v}\right)=\boldsymbol{\gamma}\!\left(\left[u^{i}\right]\left[\mathbf{e}_{i}\right],\left[v^{j}\right]\left[\mathbf{e}_{j}\right]\right)=\left[u^{i}\right]\left[v^{j}\right]\boldsymbol{\gamma}\!\left(\left[\mathbf{e}_{i}\right],\left[\mathbf{e}_{j}\right]\right)=\left[u^{i}\right]\left[v^{j}\right]\left[g_{ij}\right],
\]
where $g_{ij}=\boldsymbol{\gamma}\left(\mathbf{e}_{i},\mathbf{e}_{j}\right)$.
If $\boldsymbol{g}:V\otimes V\rightarrow K$ is a corresponding linear
tensor map we have
\[
\boldsymbol{g}\!\left(\mathbf{u}\otimes\mathbf{v}\right)=\boldsymbol{g}\!\left(\left[u^{i}\right]\left[\mathbf{e}{}_{i}\right]\otimes\left[v^{j}\right]\left[\mathbf{e}{}_{j}\right]\right)=\left[u^{i}\right]\left[v^{j}\right]\boldsymbol{g}\!\left(\left[\mathbf{e}{}_{i}\right]\otimes\left[\mathbf{e}{}_{j}\right]\right)=\left[u^{i}\right]\left[v^{j}\right]\left[g_{ij}\right],
\]
where $g_{ij}=\boldsymbol{g}\left(\mathbf{e}{}_{i}\otimes\mathbf{e}{}_{j}\right)=\boldsymbol{\gamma}\left(\mathbf{e}_{i},\mathbf{e}_{j}\right)$,
so $\left[g_{ij}\right]$ represents $\boldsymbol{\gamma}$ and $\boldsymbol{g}$
relative to $\left[\mathbf{e}{}_{i}\right]$.

In terms of tensor map composition we have
\[
\boldsymbol{g}\circ\boldsymbol{u}\otimes\boldsymbol{v}=\boldsymbol{g}\circ\left[u^{i}\right]\left[\boldsymbol{e}_{i}\right]\otimes\left[v^{j}\right]\left[\boldsymbol{e}_{j}\right]=\left[u^{i}\right]\left[v^{j}\right]\left[\boldsymbol{g}\circ\boldsymbol{e}_{i}\otimes\boldsymbol{e}_{j}\right]=\left[u^{i}\right]\left[v^{j}\right]\left[\boldsymbol{g}_{ij}\right],
\]
or, in double-index notation, 
\[
g_{\mathrm{ab}}\circ u^{\mathrm{a}}\circ v^{\mathrm{b}}=g_{\mathrm{ab}}\circ\left[u^{i}\right]\left[e_{i}^{\mathrm{a}}\right]\circ\left[v^{j}\right]\left[e_{j}^{\mathrm{b}}\right]=\left[u^{i}\right]\left[v^{j}\right]\left[g_{\mathrm{ab}}\circ e_{i}^{\mathrm{a}}\circ e_{j}^{\mathrm{a}}\right]=\left[u^{i}\right]\left[v^{j}\right]\left[\boldsymbol{g}_{ij}\right],
\]
where $\left[\boldsymbol{g}_{ij}\right]=\left[\eta\mapsto\eta g_{ij}\right]=\left[g_{ij}\right]\left[\eta\mapsto\eta\right]$,
so $\left[g_{ij}\right]$ represents $\boldsymbol{g}$ relative to
$\left[\boldsymbol{e}_{i}\right]$, and indirectly relative to $\left[\mathbf{e}{}_{i}\right]$.

Thus, $\left(\mathbf{u},\mathbf{v}\right)\mapsto\boldsymbol{\gamma}\left(\mathbf{u},\mathbf{v}\right)$
can be represented in array form as 
\[
\left(\left[u^{i}\right],\left[v^{j}\right]\right)\mapsto\left[u^{i}\right]\left[v^{j}\right]\left[g_{ij}\right],
\]
while $\mathbf{u}\otimes\mathbf{v}\mapsto\boldsymbol{g}\left(\mathbf{u}\otimes\mathbf{v}\right)$
and $\boldsymbol{u}\otimes\boldsymbol{v}\mapsto\boldsymbol{g}\circ\boldsymbol{u}\otimes\boldsymbol{v}$
(or $u^{\mathrm{a}}\circ v^{\mathrm{b}}\mapsto g_{\mathrm{ab}}\circ u^{\mathrm{a}}\circ v^{\mathrm{b}}$)
can be represented in the same way or as
\[
\left[u^{i}\right]\left[v^{j}\right]\mapsto\left[u^{i}\right]\left[v^{j}\right]\left[g_{ij}\right].
\]
In index-free array notation we have $\left(\boldsymbol{\mathfrak{u}}^{\mathrm{T}},\boldsymbol{\mathfrak{v}}^{\mathrm{T}}\right)\mapsto\boldsymbol{\mathfrak{u}}^{\mathrm{T}}\otimes\boldsymbol{\mathfrak{v}}^{\mathrm{T}}\overset{0}{\underset{2}{\boldsymbol{\mathfrak{G}}}}$
in the first case and $\boldsymbol{\mathfrak{u}}^{\mathrm{T}}\otimes\boldsymbol{\mathfrak{v}}^{\mathrm{T}}\mapsto\boldsymbol{\mathfrak{u}}^{\mathrm{T}}\otimes\boldsymbol{\mathfrak{v}}^{\mathrm{T}}\overset{0}{\underset{2}{\boldsymbol{\mathfrak{G}}}}$
in the second. We can eliminate the outer product\linebreak{}
 $\boldsymbol{\mathfrak{u}}^{\mathrm{T}}\otimes\boldsymbol{\mathfrak{v}}^{\mathrm{T}}$
in the first case by writing this mapping in the more familiar matrix
form\linebreak{}
 $\left(\boldsymbol{\mathfrak{u}}^{\mathrm{T}},\boldsymbol{\mathfrak{v}}^{\mathrm{T}}\right)\mapsto\boldsymbol{\mathfrak{u}}^{\mathrm{T}}\boldsymbol{\mathfrak{G}}\boldsymbol{\mathfrak{v}}$,
where $\boldsymbol{\mathfrak{G}}=\left[g_{ij^{\mathrm{T}}}\right]$,
since $\left[u^{i}\right]\left[v^{j}\right]\left[g_{ij}\right]=\left[u^{i}\right]\left[g_{ij^{\mathrm{T}}}\right]\left[v^{j^{\mathrm{T}}}\right]$.

\emph{(5)}. Let $\boldsymbol{t}=t_{\mathrm{a}}^{\mathrm{b}}:V\rightarrow V$
be a linear operator on $V$. For each $\mathbf{e}_{j}$ and each
$\boldsymbol{e}_{j}$ there is a unique expansion
\[
\boldsymbol{t}\left(\mathbf{e}_{j}\right)=\left.\sum\right._{i}t_{j}^{i}\mathbf{e}_{i}=\left[t_{j}^{i}\right]\left[\mathbf{e}{}_{i}\right],\qquad\boldsymbol{t}\circ\boldsymbol{e}_{j}=\left.\sum\right._{i}t_{j}^{i}\boldsymbol{e}_{i}=\left[t_{j}^{i}\right]\left[\boldsymbol{e}_{i}\right],
\]
where $\boldsymbol{e}_{k}=\left(\eta\rightarrow\eta\mathbf{e}_{k}\right)$.
If $\mathbf{v}=\left[v^{j}\right]\left[\mathbf{e}_{j}\right]$ then
\[
\boldsymbol{t}\left(\mathbf{v}\right)=\boldsymbol{t}\left(\left[v^{j}\right]\left[\mathbf{e}_{j}\right]\right)=\left[v^{j}\right]\left[\boldsymbol{t}\left(\mathbf{e}_{j}\right)\right]=\left[v^{j}\right]\left[t_{j}^{i}\right]\left[\mathbf{e}_{i}\right],
\]
and if $\boldsymbol{v}=\left[v^{j}\right]\left[\boldsymbol{e}_{j}\right]$
then 
\begin{gather*}
\boldsymbol{t}\circ\boldsymbol{v}=\boldsymbol{t}\circ\left(\left[v^{j}\right]\left[\boldsymbol{e}_{j}\right]\right)=\left[v^{j}\right]\left[\boldsymbol{t}\circ\boldsymbol{e}_{j}\right]=\left[v^{j}\right]\left[t_{j}^{i}\right]\left[\boldsymbol{e}_{i}\right],
\end{gather*}
so $\mathbf{v}\mapsto\boldsymbol{t}\left(\mathbf{v}\right)$ and $\boldsymbol{v}\mapsto\boldsymbol{t}\circ\boldsymbol{v}$
(or $v{}^{\mathrm{a}}\mapsto t_{\mathrm{a}}^{\mathrm{b}}\circ v{}^{\mathrm{a}}$)
are both represented by the array function 
\[
\left[v^{j}\right]\mapsto\left[v^{j}\right]\left[t_{j}^{i}\right]
\]
relative to $\left[\mathbf{e}{}_{i}\right]$. This array function
has the matrix form $\boldsymbol{\mathfrak{v}}^{\mathrm{T}}\mapsto\boldsymbol{\mathfrak{v}}^{\mathrm{T}}\boldsymbol{\mathfrak{T}}$;
the array function $\left[v^{j}\right]^{\mathrm{T}}\mapsto\left[t_{j}^{i}\right]^{\mathrm{T}}\left[v^{j}\right]^{\mathrm{T}}$
obtained from $\left[v^{j}\right]\mapsto\left[v^{j}\right]\left[t_{j}^{i}\right]$
by transposition $ $has the more familiar matrix form $\boldsymbol{\mathfrak{v}}\mapsto\boldsymbol{\mathfrak{S}}\boldsymbol{\mathfrak{v}}$,
where $\boldsymbol{\mathfrak{S}}=\boldsymbol{\mathfrak{T}}^{\mathrm{T}}$. 
\begin{rem*}
It can be shown that if $\left[t_{j}^{i}\right]$ represents $t_{\mathrm{a}}^{\mathrm{b}}:V\rightarrow V$
relative to $\left\{ \mathbf{e}_{i}\right\} $ then $\left[t_{j^{\mathrm{T}}}^{i^{\mathrm{T}}}\right]=\left[t_{j}^{i}\right]^{\mathrm{T}}$
represents $t_{\mathrm{b}^{*}}^{\mathrm{a}^{*}}:V^{*}\rightarrow V^{*}$
relative to the dual basis $\left\{ \boldsymbol{e}^{i}\right\} $.
This connection suggests why the rules for placing transposed indices
(transposed subscripts are located to the left of 'old' superscripts,
transposed superscripts to the right of 'old' subscripts) are the
same for transposition of arrays and outer transposition of tensor
maps.
\end{rem*}
\emph{(6)}. In general, let 
\[
\boldsymbol{T}=T_{\mathrm{a_{1}\ldots a_{m}}}^{\mathrm{b_{1}\ldots b_{n}}}:V^{\otimes m}\rightarrow V^{\otimes n}
\]
be a linear tensor map and let $\left[\mathbf{e}_{i}\right]$ be a
basis for $V$. As $\left[\mathbf{e}{}_{i_{1}}\otimes\ldots\otimes\mathbf{e}{}_{i_{n}}\right]$
is a basis for $V^{\otimes n}$, there is a unique array $\left[t_{j_{1}\ldots j_{m}}^{i_{1}\ldots i_{n}}\right]$
such that
\[
\boldsymbol{T}\left(\mathbf{e}{}_{j_{1}}\otimes\ldots\otimes\mathbf{e}{}_{j_{m}}\right)=\left[t_{j_{1}\ldots j_{m}}^{i_{1}\ldots i_{n}}\right]\left[\mathbf{e}{}_{i_{1}}\otimes\ldots\otimes\mathbf{e}{}_{i_{n}}\right],
\]
 so 
\begin{gather*}
\boldsymbol{T}\left(\mathbf{v}{}_{1}\otimes\ldots\otimes\mathbf{v}{}_{m}\right)=\boldsymbol{T}\left(\left[v_{1}^{j_{1}}\right]\left[\mathbf{e}{}_{j_{1}}\right]\otimes\ldots\otimes\left[v_{m}^{j_{m}}\right]\left[\mathbf{e}{}_{j_{m}}\right]\right)=\\
\left[v_{1}^{j_{1}}\right]\ldots\left[v_{m}^{j_{m}}\right]\left[\boldsymbol{T}\left(\mathbf{e}{}_{j_{1}}\otimes\ldots\otimes\mathbf{e}{}_{j_{m}}\right)\right]=\left[v_{1}^{j_{1}}\right]\ldots\left[v_{m}^{j_{m}}\right]\left[t_{j_{1}\ldots j_{m}}^{i_{1}\ldots i_{n}}\right]\left[\mathbf{e}{}_{i_{1}}\otimes\ldots\otimes\mathbf{e}{}_{i_{n}}\right].
\end{gather*}

In terms of tensor map composition we have
\begin{gather*}
T_{\mathrm{a_{1}\ldots a}_{m}}^{\mathrm{b_{1}\ldots b}_{n}}\circ v_{1}^{\mathrm{a_{1}}}\circ\ldots\circ v_{m}^{\mathrm{a}_{m}}=T_{\mathrm{a_{1}\ldots a}_{m}}^{\mathrm{b_{1}\ldots b}_{n}}\circ\left[v_{1}^{j_{1}}\right]\left[e_{j_{1}}^{\mathrm{a_{1}}}\right]\circ\ldots\circ\left[v_{m}^{j_{m}}\right]\left[e_{j_{m}}^{\mathrm{a}_{m}}\right]=\\
\left[v_{1}^{j_{1}}\right]\ldots\left[v_{m}^{j_{m}}\right]\left[T_{\mathrm{a_{1}\ldots a}{}_{m}}^{\mathrm{b_{1}\ldots b}_{n}}\circ e_{j_{1}}^{\mathrm{a_{1}}}\circ\ldots\circ e_{j_{m}}^{\mathrm{a}_{m}}\right]=\left[v_{1}^{j_{1}}\right]\ldots\left[v_{m}^{j_{m}}\right]\left[t_{j_{1}\ldots j_{m}}^{i_{1}\ldots i_{n}}\right]\left[e_{i_{1}}^{\mathrm{b_{1}}}\circ\ldots\circ e_{i_{n}}^{\mathrm{b}_{n}}\right],
\end{gather*}
or, using index-free notation for tensor maps,
\begin{gather*}
\boldsymbol{T}\circ\boldsymbol{v}_{1}\otimes\ldots\otimes\boldsymbol{v}_{m}=\boldsymbol{T}\circ\left[v_{1}^{j_{1}}\right]\left[\boldsymbol{e}_{j_{1}}\right]\otimes\ldots\otimes\left[v_{m}^{j_{m}}\right]\left[\boldsymbol{e}_{j_{m}}\right]=\\
\left[v_{1}^{j_{1}}\right]\ldots\left[v_{m}^{j_{m}}\right]\left[\boldsymbol{T}\circ\boldsymbol{e}_{j_{1}}\otimes\ldots\otimes\boldsymbol{e}_{j_{m}}\right]=\left[v_{1}^{j_{1}}\right]\ldots\left[v_{m}^{j_{m}}\right]\left[t_{j_{1}\ldots j_{m}}^{i_{1}\ldots i_{n}}\right]\left[\boldsymbol{e}_{i_{1}}\otimes\ldots\otimes\boldsymbol{e}_{i_{n}}\right],
\end{gather*}
where $v_{k}^{\mathrm{a}_{k}}=\boldsymbol{v}_{k}=\left(\eta\mapsto\eta\mathbf{v}_{k}\right)$
and $e_{k}^{\mathrm{\mathrm{a}}_{k}}=e_{k}^{\mathrm{\mathrm{b}}_{k}}=\boldsymbol{e}_{k}=\left(\eta\mapsto\eta\mathbf{e}_{k}\right)$.
(If $n=0$ then $\left[\mathsf{e}_{i_{1}}\otimes\ldots\otimes\mathsf{e}_{i_{n}}\right]=\left[1\right]$
and $\left[\boldsymbol{e}_{i_{1}}\otimes\ldots\otimes\boldsymbol{e}_{i_{n}}\right]=\left[\eta\mapsto\eta\right]$
by convention.) Thus, the linear maps given by
\[
\mathbf{v}_{1}\otimes\ldots\otimes\mathbf{v}_{m}\mapsto\boldsymbol{T}\left(\mathbf{v}_{1}\otimes\ldots\otimes\mathbf{v}_{m}\right)
\]
and 
\begin{gather*}
\boldsymbol{v}_{1}\otimes\ldots\otimes\boldsymbol{v}_{m}\mapsto\boldsymbol{\boldsymbol{T}}\circ\boldsymbol{v}_{1}\otimes\ldots\otimes\boldsymbol{v}_{m}\quad\mathrm{or}\\
v_{1}^{\mathrm{a_{1}}}\circ\ldots\circ v_{m}^{\mathrm{a}_{m}}\mapsto T_{\mathrm{a_{1}\ldots a}_{m}}^{\mathrm{b_{1}\ldots b}_{n}}\circ v_{1}^{\mathrm{a_{1}}}\circ\ldots\circ v_{m}^{\mathrm{a}_{m}}
\end{gather*}
are represented relative to $\left[\mathbf{e}{}_{i}\right]$ by the
array map
\[
\left[v_{1}^{j_{1}}\right]\ldots\left[v_{m}^{j_{m}}\right]\mapsto\left[v_{1}^{j_{1}}\right]\ldots\left[v_{m}^{j_{m}}\right]\left[t_{j_{1}\ldots j_{m}}^{i_{1}\ldots i_{n}}\right].
\]
This means that these maps and the corresponding separately linear
maps can ultimately be represented by an array
\[
\left[t_{j_{1}\ldots j_{m}}^{i_{1}\ldots i_{n}}\right].
\]

\subsection{Representations of classical tensors in indicial notation}

$ $\\
Since tensor maps can be represented by arrays, it is clear that we
can also use conventional indicial notation to represent tensor maps.
Thus, $v^{i}$ can represent $v^{\mathrm{a}}$, $f_{i}$ can represent
$f_{\mathrm{a}}$, $g_{ij}$ can represent $g_{\mathrm{ab}}$, and
so forth, so there is an obvious analogy between indicial notation
and the notation used for tensor maps. This formal similarity does
of course belie a big conceptual difference, however. For example,
$v^{i}$ and $g_{ij}$ are scalars, or systems of scalars, and $i,j,\ldots$
are integers, but $v^{\mathrm{a}}$ and $g_{\mathrm{ab}}$ are functions,
and $\mathsf{a,b,\ldots}$ are formal symbols associated with 'inputs'
to or 'outputs' from such functions.

As we have seen, there is a subtle formal difference between tensor
map notation and indicial notation, too. In indicial notation, both
$s_{i}t^{i}$ and $t^{i}s_{i}$ are legitimate expressions. By contrast,
neither $\mathsf{t}{}^{\mathrm{a}}\mathsf{s}{}_{\mathrm{a}}$ (bilateral
tensors) nor $t^{\mathrm{a}}\circ s_{\mathrm{a}}$ (tensor maps) are
legitimate expressions. Only $\mathsf{s}{}_{\mathrm{a}}\mathsf{t}{}^{\mathrm{a}}$
and $s_{\mathrm{a}}\circ t^{\mathrm{a}}$ are legitimate expressions,
rendered in array notion as $\left[t^{i}\right]\left[s_{i}\right]$
and in indicial notation as either $t^{i}s_{i}$ or $s_{i}t^{i}$.
\begin{rem*}
As Penrose's abstract index notation \cite{key-5} is modeled on the
indicial notation, both $s_{a}t^{a}$ and $t^{a}s_{a}$ are legitimate
expressions in the notation he proposes. This is a point where the
difference between the present approach and Penrose's notation becomes
evident. 
\end{rem*}

\subsection{Multiplication of arrays and corresponding tensorial operations}

$ $\\
We have shown that any linear tensor map $t_{\mathrm{a}_{1}\ldots\mathrm{a}_{m}}^{\mathrm{b}_{1}\ldots\mathrm{b}_{n}}$
can be represented by a scalar array $\left[t_{i_{1}\ldots i_{m}}^{j_{1}\ldots j_{n}}\right]$,
and it is not difficult to show that if $V$ is an $N$-dimensional
vector space over $K$ then this mapping is a vector space isomorphism
\[
\Lambda_{\mathcal{E}}:\mathscr{L}\left[V^{\otimes m},V^{\otimes n}\right]\rightarrow\mathscr{A}_{m}^{n}\left[K^{N}\right]
\]
 for every basis $\mathcal{E}$ in $V$. 

It can also be shown that this isomorphism is compatible with the
multiplication operations in the two vector spaces. For example, we
may infer from the corresponding result in linear algebra that if
$\Lambda_{\mathcal{E}}\!\left(s_{\mathrm{a}}^{\mathrm{b}}\right)=\left[s_{i}^{j}\right]$
and $\Lambda_{\mathcal{E}}\!\left(t{}_{\mathrm{c}}^{\mathrm{d}}\right)=\left[t_{\ell}^{k}\right]$
then 
\[
\Lambda_{\mathcal{E}}\!\left(s_{\mathrm{a}}^{\mathrm{b}}\circ t{}_{\mathrm{c}}^{\mathrm{a}}\right)^{\mathrm{T}}=\Lambda_{\mathcal{E}}\!\left(s_{\mathrm{a}}^{\mathrm{b}}\right)^{\mathrm{T}}\Lambda_{\mathcal{E}}\!\left(t{}_{\mathrm{c}}^{\mathrm{a}}\right)^{\mathrm{T}}=\left[s_{i}^{j}\right]^{\mathrm{T}}\left[t_{j}^{k}\right]^{\mathrm{T}},
\]
so $\Lambda_{\mathcal{E}}^{\mathrm{T}}$, defined by $\Lambda_{\mathcal{E}}^{\mathrm{T}}\left(\boldsymbol{t}\right)=\Lambda_{\mathcal{E}}\left(\boldsymbol{t}\right)^{\mathrm{T}}$,
is an isomorphism between vector spaces mapping composition of tensor
maps onto multiplication of arrays.

In the general case, we have
\[
\Lambda_{\mathcal{E}}\!\left(\circ{}_{\mathrm{a}_{i_{1}}\ldots\mathrm{a}_{i_{\varrho}}}^{\mathrm{d}_{j_{1}}\ldots\mathrm{d}_{j_{\varrho}}}\left(s_{\mathrm{a}_{1}\ldots\mathrm{a}_{m}}^{\mathrm{b}_{1}\ldots\mathrm{b}_{n}},t_{\mathrm{c}_{1}\ldots\mathrm{d}_{p}}^{\mathrm{d}_{1}\ldots\mathrm{d}_{q}}\right)\right)^{\mathrm{T}}=\left[s_{I}^{J}\right]^{\mathrm{T}}\left[t_{J'}^{K}\right]^{\mathrm{T}},
\]
where $J$ and $J'$ contain exactly $\varrho$ matching indices,
but there are no indices matching other indices in $I$ and $K$.
There are again vector space isomorphisms which map composition of
tensor maps onto multiplication of arrays.

Recalling that 
\[
\Lambda_{\mathcal{B}}\left(\bullet{}_{\mathrm{a}_{i_{1}}\ldots\mathrm{a}_{i_{\varrho}}}^{\mathrm{d}_{j_{1}}\ldots\mathrm{d}_{j_{\varrho}}}\left(\mathsf{s}_{\mathrm{a}_{1}\ldots\mathrm{a}_{m}}^{\mathrm{b}_{1}\ldots\mathrm{b}_{n}},\mathsf{t}_{\mathrm{c}_{1}\ldots\mathrm{d}_{p}}^{\mathrm{d}_{1}\ldots\mathrm{d}_{q}}\right)\right)=\circ{}_{\mathrm{a}_{i_{1}}\ldots\mathrm{a}_{i_{\varrho}}}^{\mathrm{d}_{j_{1}}\ldots\mathrm{d}_{j_{\varrho}}}\left(\Lambda_{\mathcal{B}}\left(\mathsf{s}_{\mathrm{a}_{1}\ldots\mathrm{a}_{m}}^{\mathrm{b}_{1}\ldots\mathrm{b}_{n}}\right),\Lambda_{\mathcal{B}}\left(\mathsf{t}{}_{\mathrm{c}_{1}\ldots\mathrm{d}_{p}}^{\mathrm{d}_{1}\ldots\mathrm{d}_{q}}\right)\right),
\]
we conclude that the spaces of bilateral tensors, tensor maps and
scalar arrays are isomorphic not only as vector spaces but also as
vector spaces equipped with certain multiplication operations.

\section{Change of basis and corresponding change of representation}

\emph{(1)}. Relative to the standard basis $\left[1\right]$ for $K$
as a vector space, $\xi\in K$ is represented by $\xi=\left[\xi\right]$,
since $\xi=\xi1$. The scalar-like tensor map $\left(\eta\mapsto\eta\xi\right)$
corresponding to $\xi$ is also represented by $\xi=\left[\xi\right]$,
this time relative to the standard basis $\left[\left\langle 1\mapsto1\right\rangle \right]\!=\!\left[\eta\mapsto\eta\right]$
for $\mathscr{L}\left[K,K\right]$, since $\left(\eta\mapsto\eta\xi\right)=\xi\left(\eta\mapsto\eta\right)$.
For non-scalars and corresponding tensor maps, there is no canonical
basis, however, and the array representation depends on the basis
chosen. 

\emph{(2)}. Let $\left[\mathbf{e}_{i}\right]$, $\left[\bar{\mathbf{e}}_{j}\right]$
be indexed bases for $V$. As there is a unique expansion $\overline{\mathbf{e}}_{j}=\sum_{i}a_{j}^{i}\mathbf{e}_{i}$
for each $\overline{\mathbf{e}}_{j}$, there is a unique invertible
array $\left[a_{j}^{i}\right]$ such that 
\[
\left[\overline{\mathbf{e}}_{j}\right]=\left[a_{j}^{i}\right]\left[\mathbf{e}_{i}\right].
\]
Define $\left[\overline{v}^{j}\right]$ by $\left[\overline{v}^{j}\right]\left[\bar{\mathbf{e}}_{j}\right]=\mathbf{v}$
for every\textbf{ $\mathbf{v}\in V$}, so that 
\[
\left[v^{i}\right]\left[\mathbf{e}_{i}\right]=\left[\overline{v}^{j}\right]\left[\overline{\mathbf{e}}_{j}\right]=\left[\overline{v}^{j}\right]\left[a_{j}^{i}\right]\left[\mathbf{e}_{i}\right].
\]
As $\left[\mathbf{e}_{i}\right]$ is a basis this implies 
\begin{equation}
\left[v^{i}\right]=\left[\overline{v}^{j}\right]\left[a_{j}^{i}\right],\label{eq:contrav}
\end{equation}
which is equivalent to
\[
\left[\overline{v}^{j}\right]=\left[v^{i}\right]\left[a_{i}^{j}\right]^{-1}\quad\mathrm{and}\quad\left[\overline{v}^{j}\right]^{\mathrm{T}}=\left(\left[a_{i}^{j}\right]^{-1}\right)^{\mathrm{T}}\left[v^{i}\right]^{\mathrm{T}}=\left(\left[a_{i}^{j}\right]^{\mathrm{T}}\right)^{\mathrm{-1}}\left[v^{i}\right]^{\mathrm{T}}.
\]
In matrix form we have $\boldsymbol{\mathfrak{v}}{}^{\mathrm{T}}=\overline{\boldsymbol{\mathfrak{v}}}^{\mathrm{T}}\boldsymbol{\mathfrak{A}}$,
or equivalently $\overline{\boldsymbol{\mathfrak{v}}}^{\mathrm{T}}=\boldsymbol{\mathfrak{v}}^{\mathrm{T}}\boldsymbol{\mathfrak{A}}^{-1}$,
or equivalently $\overline{\boldsymbol{\mathfrak{v}}}=\mathbf{\left(\boldsymbol{\mathfrak{A}}^{\mathrm{-1}}\right)^{\mathrm{T}}\boldsymbol{\mathfrak{v}}}=\mathbf{\left(\boldsymbol{\mathfrak{A}}^{\mathrm{\mathrm{T}}}\right)^{\mathrm{-1}}\boldsymbol{\mathfrak{v}}}$. 

We obtain the same result if we consider the representation of the
vector-like tensor map $v^{\mathrm{a}}$ instead.

\emph{(3)}. Recall that $\left[v^{i}\right]\left[f_{i}\right]=f_{\mathrm{a}}\!\left(\mathbf{v}\right)$
and define $\left[\overline{f}_{j}\right]$ by $\left[\overline{v}^{j}\right]\left[\overline{f}_{j}\right]=f_{\mathrm{a}}\!\left(\mathbf{v}\right)$
for every $\left[\overline{v}^{j}\right]$. Using (\ref{eq:contrav}),
we get
\[
\left[\overline{v}^{j}\right]\left[\overline{f}_{j}\right]=\left[v^{i}\right]\left[f_{i}\right]=\left[\overline{v}^{j}\right]\left[a_{j}^{i}\right]\left[f_{i}\right].
\]
As this holds for every $\left[\overline{v}^{j}\right]$ we obtain
\begin{equation}
\left[\overline{f}_{j}\right]=\left[a_{j}^{i}\right]\left[f_{i}\right].\label{eq:cov}
\end{equation}
In matrix form this is simply $\overline{\boldsymbol{\mathfrak{f}}}=\boldsymbol{\mathfrak{A}}\mathbf{\boldsymbol{\mathfrak{f}}}$.

\emph{(4)}. Consider now a tensor map (linear transformation) 
\[
t_{\mathrm{a}}^{\mathrm{b}}:V\rightarrow V,\quad\left(\mathsf{\mathbf{u}\mapsto\mathbf{v}}\right)
\]
represented by the array map
\[
\left[u^{j}\right]\mapsto\left[v^{i}\right]=\left[u^{j}\right]\left[t_{j}^{i}\right]
\]
relative to $\left[\mathbf{e}_{i}\right]$ and by
\[
\left[\overline{u}^{j}\right]\mapsto\left[\overline{v}^{i}\right]=\left[\overline{u}^{j}\right]\left[\overline{t}_{j}^{i}\right]
\]
relative to $\left[\overline{\mathbf{e}}_{i}\right]$. Using (\ref{eq:contrav})
twice, we get
\[
\left[\overline{u}^{k}\right]\left[\overline{t}_{k}^{\ell}\right]\left[a_{\ell}^{i}\right]=\left[\overline{u}^{k}\right]\left[\overline{t}_{k}^{j}\right]\left[a_{j}^{i}\right]=\left[\overline{v}^{j}\right]\left[a_{j}^{i}\right]=\left[v^{i}\right]=\left[u^{j}\right]\left[t_{j}^{i}\right]=\left[\overline{u}^{k}\right]\left[a_{k}^{j}\right]\left[t_{j}^{i}\right],
\]
so
\[
\left[\overline{u}^{k}\right]\left[\overline{t}_{k}^{\ell}\right]=\left[\overline{u}^{k}\right]\left[a_{k}^{j}\right]\left[t_{j}^{i}\right]\left[a_{i}^{\ell}\right]^{-1},
\]
and as this holds for every $\left[\overline{u}^{k}\right]$ we conclude
that 
\begin{equation}
\left[\overline{t}_{k}^{\ell}\right]=\left[a_{k}^{j}\right]\left[t_{j}^{i}\right]\left[a_{i}^{\ell}\right]^{-1}.\label{eq:covcontrav}
\end{equation}
In matrix notation, we have $\boldsymbol{\mathfrak{v}}{}^{\mathrm{T}}=\boldsymbol{\mathfrak{u}}{}^{\mathrm{T}}\boldsymbol{\mathfrak{T}}$,
so $\overline{\boldsymbol{\mathfrak{v}}}^{\mathrm{T}}\boldsymbol{\mathfrak{A}}=\overline{\boldsymbol{\mathfrak{u}}}^{\mathrm{T}}\boldsymbol{\mathfrak{A}}\boldsymbol{\mathfrak{T}}$,
so $\overline{\mathbf{\boldsymbol{\mathfrak{v}}}}^{\mathrm{T}}=\overline{\boldsymbol{\mathfrak{u}}}^{\mathrm{T}}\boldsymbol{\mathfrak{A}}\boldsymbol{\mathfrak{T}}\boldsymbol{\mathfrak{A}}^{-1}$,
so with $\overline{\mathbf{\boldsymbol{\mathfrak{v}}}}^{\mathrm{T}}=\overline{\boldsymbol{\mathfrak{u}}}^{\mathrm{T}}\overline{\boldsymbol{\mathfrak{T}}}$
we have $\overline{\boldsymbol{\mathfrak{u}}}^{\mathrm{T}}\overline{\boldsymbol{\mathfrak{T}}}=\overline{\boldsymbol{\mathfrak{u}}}^{\mathrm{T}}\boldsymbol{\mathfrak{A}}\boldsymbol{\mathfrak{T}}\boldsymbol{\mathfrak{A}}^{-1}$
for every $\overline{\boldsymbol{\mathfrak{u}}}^{\mathrm{T}}$, so
$\overline{\mathbf{\boldsymbol{\mathfrak{T}}}}=\boldsymbol{\mathfrak{A}}\boldsymbol{\mathfrak{T}}\boldsymbol{\mathfrak{A}}^{-1}$. 
\begin{rem*}
 In the linear algebra literature, the change of coordinates equation
usually has the form $\overline{\mathbf{\boldsymbol{\mathfrak{T}}}}=\boldsymbol{\mathfrak{A}}^{-1}\boldsymbol{\mathfrak{T}}\boldsymbol{\mathfrak{A}}$
rather than $\overline{\mathbf{\boldsymbol{\mathfrak{T}}}}=\boldsymbol{\mathfrak{A}}\boldsymbol{\mathfrak{T}}\boldsymbol{\mathfrak{A}}^{-1}$.
This is typically due to the fact that the matrix $\boldsymbol{\mathfrak{T}}$
representing the linear transformation $\boldsymbol{t}$ is conventionally
defined by $\boldsymbol{t}\!\left(\mathbf{e}_{i}\right)=\sum_{j}t_{j}^{i}\mathbf{e}_{j}$
rather than by $\boldsymbol{t}\!\left(\mathbf{e}_{j}\right)=\sum_{i}t_{j}^{i}\mathbf{e}_{i}$,
where $i$ and $j$ are column and row indices, respectively, of entries
in $\boldsymbol{\mathfrak{T}}$. Then the linear transformation $\left[v^{j}\right]\left[\mathbf{e}_{j}\right]\mapsto\left[v^{j}\right]\left[t_{j}^{i}\right]\left[\mathbf{e}_{i}\right]$,
corresponding to a change of basis $\left[\mathbf{e}_{j}\right]\mapsto\left[t_{j}^{i}\right]\left[\mathbf{e}_{i}\right]$
if $\left[t_{j}^{i}\right]$ is invertible, is represented by $\left[t_{j}^{i}\right]^{\mathrm{T}}$
rather than $\left[t_{j}^{i}\right]$ as above. Thus, there is no
contradiction, since we have $\left(\boldsymbol{\mathfrak{A}}\boldsymbol{\mathfrak{T}}\boldsymbol{\mathfrak{A}}^{-1}\right)^{\mathrm{T}}=\left(\boldsymbol{\mathfrak{A}}^{-1}\right)^{\mathrm{T}}\boldsymbol{\mathfrak{T}}^{\mathrm{T}}\boldsymbol{\mathfrak{A}}^{\mathrm{T}}=\left(\boldsymbol{\mathfrak{A}}^{\mathrm{T}}\right)^{\mathrm{-1}}\boldsymbol{\mathfrak{T}}^{\mathrm{T}}\boldsymbol{\mathfrak{A}}^{\mathrm{T}}$.
\end{rem*}
\emph{(5)}. Relation (\ref{eq:covcontrav}) obviously generalizes
(\ref{eq:contrav}) and (\ref{eq:cov}), and (\ref{eq:covcontrav})
can in turn be generalized to 
\begin{equation}
\left[\overline{t}_{k_{1}\cdots k_{m}}^{\ell_{1}\cdots\ell_{n}}\right]=\left[a_{k_{1}}^{j_{1}}\right]\ldots\left[a_{k_{m}}^{j_{m}}\right]\left[t_{j_{1}\cdots j_{m}}^{i_{1}\cdots i_{n}}\right]\left[a_{i_{1}}^{\ell_{1}}\right]^{-1}\ldots\left[a_{i_{n}}^{\ell_{n}}\right]^{-1}.\label{eq:tenssim}
\end{equation}
Note that (\ref{eq:tenssim}) reduces to (\ref{eq:contrav}) for an
array $\left[t^{i}\right]$ representing a vector-like tensor map
and (\ref{eq:cov}) for an array $\left[t_{j}\right]$ representing
a tensor map which is a linear form.

Using index-free array notation, we can write (\ref{eq:tenssim})
as
\begin{gather}
\underset{m}{\overset{n}{\overline{\boldsymbol{\mathfrak{T}}}}}=\overset{\otimes m}{\boldsymbol{\boldsymbol{\mathfrak{A}}}}\:\underset{m}{\overset{n}{\mathbf{\boldsymbol{\mathfrak{T}}}}}\:\underset{\otimes n}{\boldsymbol{\mathfrak{A}}}^{\!-1},\qquad\mathrm{where}\label{eq:tenssimnoind}\\
\overset{\otimes m}{\boldsymbol{\boldsymbol{\mathfrak{A}}}}=\underset{m\geq1\;\mathrm{factors}}{\underbrace{\boldsymbol{\mathfrak{A}}\otimes\cdots\otimes\boldsymbol{\mathfrak{A}}}}\,,\quad\underset{\otimes n}{\boldsymbol{\mathfrak{A}}}^{\!-1}=\underset{n\geq1\;\mathrm{factors}}{\underbrace{\boldsymbol{\mathfrak{A}^{-1}}\otimes\cdots\otimes\boldsymbol{\mathfrak{A}^{-1}}}}\,,\qquad\overset{\otimes0}{\boldsymbol{\boldsymbol{\mathfrak{A}}}}=\underset{\otimes0}{\boldsymbol{\mathfrak{A}}}^{\!-1}=\left[\delta_{i}^{j}\right].\nonumber 
\end{gather}
As $\underset{0}{\overset{0}{\overline{\boldsymbol{\mathfrak{T}}}}}=\underset{0}{\overset{0}{\mathbf{\boldsymbol{\mathfrak{T}}}}}$,
(\ref{eq:tenssimnoind}) is consistent with the assertion that the
representation of a scalar-like tensor map does not depend on a choice
of basis.

We conclude that a classical tensor interpreted as a tensor map is
a tensor also in the traditional sense of 'something the coordinates
of which transform in accordance with (\ref{eq:tenssim}), or an equivalent
formula, under a change of basis'.
\begin{rem*}
Because conventional indicial notation allows $\left[a^{i}\right]\left[b_{i}\right]$
to be written as either $a^{i}b_{i}$ or $b_{i}a^{i}$, there are
many ways of expressing (\ref{eq:tenssim}) in indicial notation.
For example, \cite{key-10} gives the form
\[
\overline{T}_{j_{1}\cdots j_{m}}^{i_{1}\cdots i_{n}}=\breve{A}_{j_{1}}^{\ell_{1}}\cdots\breve{A}_{j_{m}}^{\ell_{m}}A_{k_{1}}^{i_{1}}\cdots A_{k_{n}}^{i_{n}}\,\, T_{\ell_{1}\cdots\ell_{m}}^{k_{1}\cdots k_{n}},\quad\mathrm{where}\quad A_{i}^{k}\breve{A}_{k}^{j}=\delta_{i}^{j}.
\]

\end{rem*}

\section{Final remarks}

\subsection{On notational conventions}

$ $\\
\emph{(1)}. The roles played by subscripts and superscripts of arrays
and tensors are determined by conventions. For example, a vector-like
tensor map could be denoted $v_{\mathrm{a}}$ as well as $v^{\mathrm{a}}$,
and in either case it could be represented by an array denoted $\left[v_{i}\right]$
as well as an array denoted $\left[v^{i}\right]$. Conventions such
as these must of course form an internally consistent system, and
what seems to be the most satisfactory system unfortunately breaks
with tradition with respect to the roles played by subscripts and
superscripts in systems of scalars. For example, a 'row vector' is
often denoted by $\left[\xi{}_{i}\right]=\left[\xi{}_{1}\cdots\xi{}_{N}\right]$
while a 'column vector' is denoted by 
\[
\left[\xi{}^{i}\right]=\left[\begin{array}{c}
\xi{}^{1}\\
\vdots\\
\xi{}^{N}
\end{array}\right],
\]
but here the opposite convention is used; a 'row vector' is written
as $\left[\xi{}^{i}\right]$, a 'column vector' as $\left[\xi{}_{i}\right]$.
This is because consistency with other conventions requires that superscripts
denote column numbers while subscripts denote row numbers.

\emph{(2)}. As we have seen, there is an asymmetry between shifting
of indices and raising/lowering of indices; for example, $t_{\mathrm{a}}^{\mathrm{b}^{\vee}}=t_{\mathrm{ab}^{*}}$
but $t_{\mathrm{a}}^{\mathrm{b}^{\wedge}}=t_{\mathrm{b'a}}$. This
is not a deep difference, however, but an effect of a notational convention
based on the rule that arguments of a function are written to the
right of the function symbol. But if we write $\left(\mathbf{u}\right)\! g$
instead of $g\!\left(\mathbf{u}\right)$, $\boldsymbol{u}\circ\boldsymbol{g}$
instead of $\boldsymbol{g}\circ\boldsymbol{u}$, $\left(\left(\mathbf{u}\right)\! g\right)\! h$
instead of $h\!\left(g\!\left(\mathbf{u}\right)\right)$, $t_{\mathrm{a}}^{\mathrm{b}}\circ g_{\mathrm{bb'}}$
instead of $g_{\mathrm{bb'}}\circ t_{\mathrm{a}}^{\mathrm{b}}$ and
so on, then, in particular, $t_{\mathrm{a}}^{\mathrm{b}^{\wedge}}=t_{\mathrm{a}}^{\mathrm{b}}\circ g_{\mathrm{bb'}}=t_{\mathrm{a}b'}$,
analogous to $t_{\mathrm{a}}^{\mathrm{b}^{\vee}}=t_{\mathrm{ab}^{*}}$. 

Reversing the notation for function application and function composition
in this way would obviously lead to far-reaching changes of formulas.
For example, the connection between tensor maps and arrays described
in Subsection 12.3 would look more natural, with $\boldsymbol{g}\circ\boldsymbol{h}$
represented by $\left[g_{i}^{j}\right]\left[h_{j}^{k}\right]$ instead
of $\left[h_{i}^{j}\right]\left[g_{j}^{k}\right]$. On the whole,
the reverse notation seems to be more natural, but the usual notation
is of course anchored to a very strong tradition.

\subsection{Arrays as classical tensors}

$ $\\
Recall that we have been discussing two main types of \emph{interpretations}
of classical tensors. There are, on the one hand, 'bilateral' interpretations,
where classical tensors are vectors in $V^{*\otimes m}\otimes V^{\otimes n}$,
$V^{\otimes n}\otimes V^{*\otimes m}$, $\mathscr{L}\left[V^{*m}\!\times\! V^{n},K\right]$
or $\mathscr{L}\left[V^{n}\!\times\! V^{*m},K\right]$, and on the
other hand interpretations of classical tensors as tensor maps in
\linebreak{}
$\mathscr{L}\left[V^{\otimes m},V^{\otimes n}\right]$ or $\mathscr{L}\left[V^{m},V^{\otimes n}\right]$.
These interpretations have been shown to be formally equivalent. We
also showed how these tensors could be \emph{represented} by arrays,
or in terms of classical indicial notation, in both cases relative
to a choice of basis.

It is possible, though, to adopt a more abstract point of view. The
formal equivalence of the two interpretations of classical tensors
suggest that both are \emph{realizations} of an abstract vector space
with additional structure, what we may loosely refer to as a\emph{
multialgebra.} The elements of such an algebraic structure are equipped
with strings of subscripts and superscripts, and the multialgebra
is endowed with multiplication operations such as those discussed
in some detail earlier. There are thus ascending sequences of subspaces
\[
\ldots,\mathcal{T}_{i-1}^{j}\subseteq\mathcal{T}_{i}^{j}\subseteq\mathcal{T}_{i+1}^{j},\ldots\qquad\qquad\ldots,\mathcal{T}_{i}^{j-1}\subseteq\mathcal{T}_{i}^{j}\subseteq\mathcal{T}_{i}^{j+1},\ldots
\]
such that if $s\in\mathcal{T}_{m}^{n}$, $t\in\mathcal{T}_{p}^{q}$
and $st$ is a product of $s$ and $t$ then $st\in\mathcal{T}_{m+p}^{n+q}$. 

In addition, arrays -- and corresponding indicial notation symbols
-- are also (potentially) double-indexed vectors in their own right.
Arrays, in particular, are equipped with (possibly empty) sequences
of subscripts and superscripts, and multiplication of arrays is entirely
analogous to these operations for the two other types of double-indexed
vectors considered here. 

At the same time, arrays continue to represent double-indexed vectors.
For example, an array represents itself relative to a unique standard
basis for the vector space of such arrays in the same way that an
element of $\mathbb{R}^{n}$ represents itself as a vector in $\mathbb{R}^{n}$
relative to the standard basis for $\mathbb{R}^{n}$.

\end{document}